\documentclass[a4paper,12pt]{article}

\usepackage{amsmath,amssymb,amsthm,amscd}

\newcommand{\Fg}{\mathfrak{g}}
\newcommand{\Fh}{\mathfrak{h}}

\newcommand{\BZ}{\mathbb{Z}}
\newcommand{\BC}{\mathbb{C}}

\newcommand{\CL}{\mathcal{L}}
\newcommand{\CB}{\mathcal{B}}
\newcommand{\CW}{\mathcal{W}}

\newcommand{\wt}{\mathop{\rm wt}\nolimits}
\newcommand{\cl}{\mathop{\rm cl}\nolimits}
\newcommand{\lev}{\mathop{\rm lev}\nolimits}
\newcommand{\hw}{\mathop{\rm h.w.}\nolimits}

\newcommand{\sw}{(\Fh_{\cl}^{\ast})^{0}}
\newcommand{\pos}{P_{\omega}^{\ast}}

\newcommand{\Hom}{\mathop{\rm Hom}\nolimits}
\newcommand{\Aut}{\mathop{\rm Aut}\nolimits}
\newcommand{\End}{\mathop{\rm End}\nolimits}
\newcommand{\GL}{\mathop{\rm GL}\nolimits}
\newcommand{\Bij}{\mathop{\rm Bij}\nolimits}

\newcommand{\vpi}{\varpi}
\newcommand{\ve}{\varepsilon}
\newcommand{\vp}{\varphi}

\newcommand{\ba}[1]{\overline{#1}}
\newcommand{\ti}[1]{\widetilde{#1}}
\newcommand{\ha}[1]{\widehat{#1}}

\makeatletter
\@addtoreset{equation}{subsection}

\renewcommand\section{\@startsection{section}{1}{0pt}
{-3.5ex plus -1ex minus -.2ex}{1.0ex plus .2ex}{\large\bf}}
\renewcommand\subsection{\@startsection{subsection}{1}{0pt}
{2.5ex plus 1ex minus .2ex}{-1em}{\bf}}
\makeatother

\theoremstyle{plain}
\newtheorem{thm}{Theorem}[subsection]
\newtheorem{lem}[thm]{Lemma}
\newtheorem{prop}[thm]{Proposition}

\newtheorem{conj}[thm]{Conjecture}
\newtheorem*{claim}{Claim}

\theoremstyle{definition}
\newtheorem{dfn}[thm]{Definition}
\newtheorem{ass}[thm]{Assumption}
\newtheorem*{obs}{Observations}

\theoremstyle{remark}
\newtheorem{rem}[thm]{Remark}

\begin{document}

\setlength{\baselineskip}{18pt}
\setcounter{section}{-1}

\title{{\Large\bf 
Construction of perfect crystals conjecturally \\[1.5mm]
corresponding to Kirillov-Reshetikhin modules \\
over twisted quantum affine algebras}}
\author{
 Satoshi Naito \\ 
 \small Institute of Mathematics, University of Tsukuba, \\
 \small Tsukuba, Ibaraki 305-8571, Japan \ 
 (e-mail: {\tt naito@math.tsukuba.ac.jp})
 \\[2mm] and \\[2mm]
 Daisuke Sagaki \\ 
 \small Institute of Mathematics, University of Tsukuba, \\
 \small Tsukuba, Ibaraki 305-8571, Japan \ 
 (e-mail: {\tt sagaki@math.tsukuba.ac.jp})
}
\date{}
\maketitle

%
%
%
\begin{abstract}
{\setlength{\baselineskip}{12pt}
Assuming the existence of the perfect crystal bases of 
Kirillov-Reshetikhin modules over simply-laced quantum affine algebras, 
we construct certain perfect crystals for twisted quantum affine algebras, 
and also provide compelling evidence that the constructed crystals are 
isomorphic to the conjectural crystal bases of Kirillov-Reshetikhin modules 
over twisted quantum affine algebras.}
\end{abstract}
%
%
%
%
\section{Introduction.}
\label{sec:intro}
The finite-dimensional irreducible modules 
over quantum affine algebras $U_{q}^{\prime}(\Fg)$ 
with (classical) weight lattice $P_{\cl}$ 
have been extensively studied from various viewpoints, 
but the study of these modules from the viewpoint of 
the crystal base theory still seems to be insufficient.
This is mainly because unlike (infinite-dimensional) 
integrable highest weight modules over quantum affine algebras 
$U_{q}(\Fg)$ with (affine) weight lattice $P$, 
finite-dimensional irreducible $U_{q}^{\prime}(\Fg)$-modules 
do not have a crystal base in general. 
(It is not even known which finite-dimensional 
irreducible $U_{q}^{\prime}(\Fg)$-modules have a crystal base.) 
However, it is conjectured (see Conjecture~\ref{conj:hkott} below, 
and also \cite[\S2.3]{HKOTY}, \cite[\S2.3]{HKOTT}) 
that a certain important class of finite-dimensional irreducible 
$U_{q}^{\prime}(\Fg)$-modules, 
called Kirillov-Reshetikhin modules (KR modules for short), 
do have a crystal base.
In this paper, assuming the existence of 
the perfect crystal bases $\CB^{i,s}$ of the KR modules 
over the simply-laced quantum affine algebras 
$U_{q}^{\prime}(\Fg)$, we construct certain perfect crystals 
$\ha{\CB}^{i,s}$ for twisted quantum affine algebras 
$U_{q}^{\prime}(\ha{\Fg})$.
Furthermore, we explicitly describe, in almost all cases, 
how the crystals $\ha{\CB}^{i,s}$ decompose, 
when regarded as a $U_{q}(\ha{\Fg}_{\ha{I}_{0}})$-module by restriction, 
into a direct sum of the crystal bases of 
irreducible highest weight $U_{q}(\ha{\Fg}_{\ha{I}_{0}})$-modules, 
where $U_{q}(\ha{\Fg}_{\ha{I}_{0}})$ is 
the quantized universal enveloping algebra of 
the (canonical) finite-dimensional, reductive Lie subalgebra 
$\ha{\Fg}_{\ha{I}_{0}}$ of $\ha{\Fg}$.
These results motivate a conjecture that 
the crystals $\ha{\CB}^{i,s}$ are isomorphic to 
the crystal bases of certain KR modules with 
specified Drinfeld polynomials (see \S\ref{subsec:conj} below) over 
the twisted quantum affine algebras $U_{q}^{\prime}(\ha{\Fg})$, 
since they agree with the conjectural branching rules 
(with $q = 1$) in \cite[Appendix~A]{HKOTT}.

We now describe our results more precisely.
Let $\Fg$ be a simply-laced affine Lie algebra over $\BC$, i.e., 
let $\Fg$ be the Kac-Moody algebra $\Fg(A)$ over $\BC$ 
associated to the generalized Cartan matrix (GCM for short) 
$A = (a_{ij})_{i,j \in I}$ of type 
$A_{2n-1}^{(1)} \, (n \ge 2)$, 
$A_{2n}^{(1)} \, (n \geq 1)$, 
$D_{n+1}^{(1)} \, (n \geq 3)$, 
$D_{4}^{(1)}$, or $E_{6}^{(1)}$, 
where $I$ is an index set for the simple roots 
(numbered as in \S\ref{subsec:diag-aff} below). 
We denote by $U_{q}(\Fg)$ (resp., $U_{q}^{\prime}(\Fg)$) 
the associated quantum affine algebra over 
$\BC(q)$ with $P$ (resp., $P_{\cl}$) as the weight lattice.
For each $i \in I_{0} := I \setminus \{0\}$, 
$s \in \BZ_{\ge 1}$, and 
$\zeta \in \BC(q)^{\times}:=\BC(q) \setminus \{0\}$, 
let $W^{(i)}_{s}(\zeta)$ be the finite-dimensional 
irreducible module over the quantum affine algebra 
$U_{q}^{\prime}(\Fg)$ whose Drinfeld polynomials 
$P_{j}(u) \in \BC(q)[u]$ for $j \in I_0$ are given by:
\begin{equation*}
P_{j}(u) = 
\begin{cases}
{\displaystyle\prod_{k=1}^{s}} 
(1 - \zeta q^{s+2-2k} u) & \text{if $j = i$,}\\[7mm]
1 & \text{otherwise}.
\end{cases}
\end{equation*}
(Here we are using the Drinfeld realization of 
$U_{q}^{\prime}(\Fg)$, and the classification of 
its finite-dimensional irreducible modules of ``type $1$'' by 
Drinfeld polynomials; see \cite{CP1} for details.)
We call this $U_{q}^{\prime}(\Fg)$-module $W^{(i)}_{s}(\zeta)$ 
a Kirillov-Reshetikhin module (KR module for short) 
over $U_{q}^{\prime}(\Fg)$. 
It is conjectured (see Conjecture~\ref{conj:hkott} below, and 
also \cite[\S2.3]{HKOTY}) that for every $i \in I$ and $s \in \BZ_{\ge 1}$, 
there exists some $\zeta^{(i)}_{s} \in \BC(q)^{\times}$ 
such that the KR module $W^{(i)}_{s}(\zeta^{(i)}_{s})$ 
has a crystal base.

Let $\omega : I \rightarrow I$ be a nontrivial diagram automorphism 
such that $\omega(0) = 0$, and $\ha{\Fg}$ the corresponding orbit 
Lie algebra of $\Fg$ (see \S\ref{subsec:orbit-aff} below for the definition).
Note that the orbit Lie algebra $\ha{\Fg}$ is 
a twisted affine Lie algebra, i.e., 
$\ha{\Fg}$ is the Kac-Moody algebra $\Fg(\ha{A})$ 
over $\BC$ associated to the GCM 
$\ha{A} = (\ha{a}_{ij})_{i,j \in \ha{I}}$ of type 
$D_{n+1}^{(2)} \, (n \ge 2)$, 
$A_{2n}^{(2)} \, (n \ge 1)$, 
$A_{2n-1}^{(2)} \, (n \ge 3)$, 
$D_{4}^{(3)}$, or $E_{6}^{(2)}$, 
where $\ha{I} \subset I$ is 
a certain complete set (containing $0 \in I$) of 
representatives of the $\omega$-orbits in $I$, and 
is also an index set for the simple roots of $\ha{\Fg}$.
In this paper, we assume that 
for (arbitrarily) fixed $i \in \ha{I}_{0} := \ha{I} \setminus \{0\}$ and 
$s \in \BZ_{\ge 1}$, there exists some $\zeta^{(i)}_{s} \in \BC(q)^{\times}$ 
such that for every $0 \le k \le N_{i} - 1$, the KR module 
$W^{(\omega^{k}(i))}_{s}(\zeta^{(i)}_{s})$ 
over $U_{q}^{\prime}(\Fg)$ has a crystal base, 
denoted by $\CB^{\omega^{k}(i), s}$, 
where $N_{i} \in \BZ_{\ge 1}$ is the number of elements of 
the $\omega$-orbit of $i \in \ha{I}_0$ in $I$; 
here we note that the $\zeta^{(i)}_{s} \in \BC(q)^{\times}$ is 
assumed to be independent of $0 \le k \le N_{i}-1$.
We further assume that 
the $\CB^{\omega^{k}(i), s}$, 
$0 \le k \le N_{i}-1$, are all perfect 
$U_{q}^{\prime}(\Fg)$-crystals of level $s$
(in the sense of Definition~\ref{dfn:perfect}).
Now, for the (fixed) $i \in \ha{I}_{0}$ and $s \in \BZ_{\geq 1}$, 
we define the tensor product $U_{q}^{\prime}(\Fg)$-crystal 
$\ti{\CB}^{i,s}$ equipped 
with the Kashiwara operators $e_{j}$ and $f_{j}$, $j \in I$, by:
\begin{equation*}
\ti{\CB}^{i,s} = 
\CB^{i,s} \otimes \CB^{\omega(i),s} \otimes \cdots \otimes 
\CB^{\omega^{N_{i}-1}(i), s}, 
\end{equation*}
on which the diagram automorphism 
$\omega :I \rightarrow I$ acts 
in a canonical way (see \S\ref{subsec:perfixed} below).
Also, we define $\omega$-Kashiwara operators 
$\ti{e}_{j}$ and $\ti{f}_{j}$ on $\ti{\CB}^{i,s}$ for 
$j \in \ha{I}$ by (see \eqref{eq:o-kas-op}, 
and also Remark~\ref{rem:link}):
\begin{equation*}
\ti{x}_{j} = 
\begin{cases}
x_{j}\,x_{\omega(j)}^{2}\,x_{j} & 
  \text{if $N_{j} = 2$ and $a_{j, \omega(j)} = a_{\omega(j), j} = -1$}, \\[3mm]
x_{j}\,x_{\omega(j)} \cdots x_{\omega^{N_{j}-1}(j)} & 
  \text{if $a_{j, \omega^{k}(j)} = 0$ for all $1 \le k \le N_{j}-1$},
\end{cases}
\end{equation*}
where $x$ is either $e$ or $f$.
Then, the $\omega$-Kashiwara operators 
$\ti{e}_{j}$ and $\ti{f}_{j}$, $j \in \ha{I}$, 
stabilize the fixed point subset $\ha{\CB}^{i,s}$ of 
$\ti{\CB}^{i,s}$ under the action of the diagram automorphism $\omega$,
and hence equip the $\ha{\CB}^{i,s}$ 
with a structure of $U_{q}^{\prime}(\ha{\Fg})$-crystal, 
where $U_{q}^{\prime}(\ha{\Fg})$ denotes 
the quantum affine algebra associated to $\ha{\Fg}$.
Furthermore, we prove that the $\ha{\CB}^{i,s}$ 
is a perfect $U_{q}^{\prime}(\ha{\Fg})$-crystal of level $s$ 
(in the sense of Definition~\ref{dfn:perfect}), 
thereby establishing Theorem~\ref{thm:main}.

Because the $U_{q}^{\prime}(\ha{\Fg})$-crystal 
$\ha{\CB}^{i, s}$ for the (fixed) 
$i \in \ha{I}_{0}=\ha{I} \setminus \{0\}$ and 
$s \in \BZ_{\ge 1}$ is perfect (and hence regular),
it decomposes, under restriction, 
into a direct sum of the crystal bases of 
integrable highest weight modules over 
the quantized universal enveloping algebra 
$U_{q}(\ha{\Fg}_{\ha{I}_{0}})$ of 
the finite-dimensional, reductive Lie subalgebra 
$\ha{\Fg}_{\ha{I}_{0}}$ of $\ha{\Fg}$ 
corresponding to $\ha{I}_{0} \subset \ha{I}$.
In fact, we can give, in almost all cases, 
an explicit description 
(see \S\ref{subsec:case-a} -- \S\ref{subsec:case-e}) 
of the branching rule 
with respect to the restriction to 
$U_{q}(\ha{\Fg}_{\ha{I}_{0}})$, i.e., 
how the $\ha{\CB}^{i,s}$ decomposes into 
a disjoint union of connected components 
as a $U_{q}(\ha{\Fg}_{\ha{I}_{0}})$-crystal.
This result deserves to be supporting evidence 
that the $\ha{\CB}^{i, s}$ is isomorphic 
as a $U_{q}^{\prime}(\ha{\Fg})$-crystal 
to the conjectural crystal base of 
a certain KR module with specified Drinfeld polynomials
(see \S\ref{subsec:conj} below), 
denoted by $\ha{W}^{(i)}_{s}(\ha{\zeta}^{(i)}_{s})$, 
over $U_{q}^{\prime}(\ha{\Fg})$, 
since our branching rule for the $\ha{\CB}^{i, s}$ 
indeed agrees with the ``branching rule'' (with $q = 1$) 
conjectured in \cite[Appendix~A]{HKOTT} 
for the KR module $\ha{W}^{(i)}_{s}(\ha{\zeta}^{(i)}_{s})$ 
over $U_{q}^{\prime}(\ha{\Fg})$, 
regarded as a $U_{q}(\ha{\Fg}_{\ha{I}_{0}})$-module 
by restriction.

Finally, we should mention the relationship 
between the $U_{q}^{\prime}(\ha{\Fg})$-crystals 
$\ha{\CB}^{i, s}$, $i \in \ha{I}_{0}$, $s \in \BZ_{\ge 1}$, and 
``virtual'' crystals defined in \cite{OSS1}, \cite{OSS2}. 
Since the $U_{q}^{\prime}(\ha{\Fg})$-crystals $\ha{\CB}^{i,s}$, 
$i \in \ha{I}_0$, $s \in \BZ_{\ge 1}$, are perfect 
(and hence simple), their crystal graphs are connected. 
Therefore, the $U_{q}^{\prime}(\ha{\Fg})$-crystals $\ha{\CB}^{i,s}$ 
coincide with virtual crystals at least in the cases 
where 
$\Fg$ is of type $A_{2n-1}^{(1)}$ and 
$\ha{\Fg}$ is of type $D_{n+1}^{(2)}$ for $n \ge 2$, 
$\Fg$ is of type $D_{n+1}^{(1)}$ and 
$\ha{\Fg}$ is of type $A_{2n-1}^{(2)}$ for $n \ge 3$, 
$\Fg$ is of type $D_{4}^{(1)}$ and 
$\ha{\Fg}$ is of type $D_{4}^{(3)}$, 
and where
$\Fg$ is of type $E_{6}^{(1)}$ and 
$\ha{\Fg}$ is of type $E_{6}^{(2)}$. 
This clarifies the representation-theoretical 
meaning of virtual crystals.

The organization of this paper is as follows. 
In \S\ref{sec:crystal}, 
we briefly review some basic notions 
in the theory of crystals for quantum affine algebras, and 
recall a conjecture about 
the existence of the crystal bases of KR modules.
In \S\ref{sec:main}, 
we first fix the notation for diagram automorphisms of 
simply-laced affine Lie algebras $\Fg$, and 
also for corresponding orbit Lie algebras $\ha{\Fg}$.
Then we define the $U_{q}^{\prime}(\ha{\Fg})$-crystals 
$\ha{\CB}^{i, s}$, $i \in \ha{I}_{0}$, $s \in \BZ_{\ge 1}$, 
and then state our main result (Theorem~\ref{thm:main}).
In \S\ref{sec:pre}, 
we show some technical propositions and 
lemmas about the fixed point subsets of regular crystals, 
or of tensor products of regular crystals 
under the action of the diagram automorphism $\omega$, 
which will be needed later.
In \S\ref{sec:prf-main}, 
we prove the perfectness of 
the $U_{q}^{\prime}(\ha{\Fg})$-crystals 
$\ha{\CB}^{i, s}$, $i \in \ha{I}_{0}$, $s \in \BZ_{\ge 1}$, 
thereby establishing Theorem~\ref{thm:main}.
In \S\ref{sec:branch}, 
we give explicit descriptions of 
the branching rules for 
the $U_{q}^{\prime}(\ha{\Fg})$-crystals 
$\ha{\CB}^{i, s}$, $i \in \ha{I}_{0}$, 
$s \in \BZ_{\ge 1}$, with a few exceptions, 
and then propose a conjecture that 
for each $i \in \ha{I}_{0}$ and $s \in \BZ_{\ge 1}$, 
the $\ha{\CB}^{i, s}$ is isomorphic 
as a $U_{q}^{\prime}(\ha{\Fg})$-crystal
to the conjectural crystal base of a certain KR module 
over the twisted quantum affine algebra 
$U_{q}^{\prime}(\ha{\Fg})$.
%
%
%
%
%
%
%
%
\section{Crystals for quantum affine algebras.}
\label{sec:crystal}
%
%
%
%
\subsection{Cartan data.}
\label{subsec:cd}
Let $A=(a_{ij})_{i,j \in I}$ be a symmetrizable generalized 
Cartan matrix (GCM for short). 
A Cartan datum for the GCM $A$ 
is, by definition, a quintuplet 
$(A, P, P^{\vee}, \Pi, \Pi^{\vee})$ consisting of 
the GCM $A$, 
a free $\BZ$-module $P^{\vee}$ of finite rank, 
its dual $P:=\Hom_{\BZ}(P^{\vee},\BZ)$, 
a subset $\Pi^{\vee}:=\bigl\{h_{j}\bigr\}_{j \in I}$ 
of $P^{\vee}$, and 
a subset $\Pi:=\bigl\{\alpha_{j}\bigr\}_{j \in I}$ 
of $P$ satisfying 
$\alpha_{k}(h_{j})=a_{jk}$ for $j,\,k \in I$. 
Further, 
we assume that the elements $h_{j}$, $j \in I$, 
of $\Pi^{\vee} \subset P^{\vee}$ are
linearly independent; 
however, we do not assume that 
the elements $\alpha_{j}$, $j \in I$, 
of $\Pi \subset P$ are
linearly independent. 
%
%
%
%
\subsection{Crystals.}
\label{subsec:crystal}
Let us briefly recall some basic notions 
in the theory of crystals from 
\cite[Chap.~4, \S4.5]{HK} (see also \cite[\S7]{Kas3}).
Let $A=(a_{ij})_{i,j \in I}$ be a symmetrizable GCM, and 
let $(A, P, P^{\vee}, \Pi, \Pi^{\vee})$ be a Cartan 
datum for the GCM $A$. 
A crystal associated to the Cartan datum 
$(A, P, P^{\vee}, \Pi, \Pi^{\vee})$ is 
a set $\CB$ equipped with maps 
 $\wt:\CB \rightarrow P$, 
 $e_{j},\,f_{j}:
  \CB \cup \{\theta\} \rightarrow 
  \CB \cup \{\theta\}$, $j \in I$, and 
 $\ve_{j},\vp_{j}:\CB \rightarrow 
 \BZ \cup \{-\infty\}$, $j \in I$, 
satisfying Conditions (1) -- (7) of 
\cite[Definition~4.5.1]{HK} 
(with $\ti{e}_{j}$, $\ti{f}_{j}$, $0$ replaced by 
$e_{j}$, $f_{j}$, $\theta$, respectively). 
We call the map $e_{j}$ (resp., $f_{j}$) the 
raising (resp., lowering) Kashiwara operator with respect 
to $\alpha_{j} \in \Pi$, and understand that 
$e_{j}\theta=f_{j}\theta=\theta$ for all $j \in I$. 
A crystal $\CB$ is said to be semiregular if 
\begin{equation}
\ve_{j}(b)=
  \max\bigl\{m \ge 0 \mid e_{j}^{m}b \ne \theta\bigr\}, \qquad
\vp_{j}(b)=
  \max\bigl\{m \ge 0 \mid f_{j}^{m}b \ne \theta\bigr\}, 
\end{equation}
for all $b \in \CB$ and $j \in I$; every 
crystal treated in this paper is semiregular. 

Let $\CB_{1}$, $\CB_{2}$ be crystals associated to 
the Cartan datum $(A, P, P^{\vee}, \Pi, \Pi^{\vee})$ above. 
We define the tensor product crystal 
$\CB_{1} \otimes \CB_{2}$ of 
$\CB_{1}$ and $\CB_{2}$ as in \cite[Definition~4.5.3]{HK} 
(see also \cite[\S7.3]{Kas3}). Note that 
the definition of tensor product crystals  
in \cite{OSS1} and \cite{OSS2} is different from 
\cite[Definition~4.5.3]{HK}, and hence from ours; 
the roles of $e_{j}$ and $f_{j}$ are 
interchanged for each $j \in I$
(see, for example, \cite[(2.10) and (2.11)]{OSS1}). 
%
%
%
%
\subsection{Quantum affine algebras.}
\label{subsec:qaa}
From now throughout this paper, we assume that 
a GCM $A=(a_{ij})_{i,j \in I}$ is of affine type. 
Take a special vertex $0 \in I$ as in 
\cite[\S4.8, Tables Aff~1 -- Aff~3]{Kac}, 
and set $I_{0}:=I \setminus \{0\}$. 
Let $\Fg=\Fg(A)$ be the affine Lie algebra 
over the field $\BC$ of complex numbers 
associated to the GCM $A$ of affine type. Then
\begin{equation}
\Fh=\left(\bigoplus_{j \in I} \BC h_{j}\right) \oplus \BC d
\end{equation}
is a Cartan subalgebra of $\Fg$, with 
$h_{j}$, $j \in I$, 
the simple coroots, 
and $d$ the scaling element. 
The simple roots $\alpha_{j} \in \Fh^{\ast}:=
\Hom_{\BC}(\Fh,\BC)$, $j \in I$, and 
fundamental weights $\Lambda_{j} \in \Fh^{\ast}$,
 $j \in I$, are defined by 
(see \cite[Chap.~10, \S10.1]{HK}):
%
%
\begin{equation} \label{eq:affine01}
\begin{array}{l}
\alpha_{k}(h_{j})=a_{jk}, \quad 
\alpha_{k}(d)=\delta_{k,0}, \\[2mm]
\Lambda_{k}(h_{j})=\delta_{k,j}, \quad
\Lambda_{k}(d)=0,
\end{array} 
\end{equation}
for $j,\,k \in I$. 
Let $E_{j}$, $F_{j}$, $j \in I$, be 
the Chevalley generators of $\Fg$, where 
$E_{j}$ (resp., $F_{j}$) corresponds to 
the simple root $\alpha_{j}$ (resp., $-\alpha_{j}$), 
and let 
\begin{equation}
\delta=\sum_{j \in I} a_{j}\alpha_{j} \in \Fh^{\ast}
\qquad \text{ and } \qquad 
c=\sum_{j \in I} a^{\vee}_{j} h_{j} \in \Fh
\end{equation}
be the null root and 
the canonical central element of $\Fg$, 
respectively. 
We take a dual weight lattice $P^{\vee}$ 
and a weight lattice $P$ as follows:
%
%
\begin{equation} \label{eq:lattices}
P^{\vee}=
\left(\bigoplus_{j \in I} \BZ h_{j}\right) \oplus \BZ d \, 
\subset \Fh
\quad \text{and} \quad 
P= 
\left(\bigoplus_{j \in I} \BZ \Lambda_{j}\right) \oplus 
   \BZ \left(\frac{1}{a_{0}}\delta\right)
   \subset \Fh^{\ast}.
\end{equation}
Clearly, we have $P \cong \Hom_{\BZ}(P^{\vee},\BZ)$. 
Here we should note that $a_{0} = 1$ except 
the case where $\Fg$ is of type $A_{2n}^{(2)}$, and 
$a_{0}=2$ in the case where $\Fg$ is of type 
$A_{2n}^{(2)}$. It is easily seen that the quintuplet
$(A,P,P^{\vee},\Pi,\Pi^{\vee})$ is a Cartan datum 
for the affine type GCM $A=(a_{ij})_{i,j \in I}$.
Let $U_{q}(\Fg)=\langle E_{j}, F_{j}, q^{h} \mid 
 j \in I,\,h \in P^{\vee}\rangle$ be 
the quantized universal enveloping algebra 
of $\Fg$ over the field $\BC(q)$ of 
rational functions in $q$ (with complex coefficients) 
with weight lattice $P$, and Chevalley generators 
$E_{j}$, $F_{j}$, $j \in I$. 

Now, we set 
%
%
\begin{equation} \label{eq:dual-cl}
\Fh_{\cl}:=\bigoplus_{j \in I} \BC h_{j} 
 \subset \Fh
\quad \text{and} \quad 
P_{\cl}^{\vee}:=\bigoplus_{j \in I} \BZ h_{j} 
 \subset P^{\vee}.
\end{equation}
For each $\lambda \in \Fh^{\ast}$, we define 
$\cl(\lambda) \in \Fh^{\ast}_{\cl}:=(\Fh_{\cl})^{\ast}$
to be the restriction $\lambda|_{\Fh_{\cl}}$ of 
$\lambda \in \Fh^{\ast}$ to $\Fh_{\cl}$ (we simply 
write $\lambda$ for $\cl(\lambda)$
if there is no fear of confusion). It follows 
that $\Fh_{\cl}^{\ast}=\cl(\Fh^{\ast}) 
\cong \Fh^{\ast}/\BC\delta$ as $\BC$-vector spaces, 
and 
%
%
\begin{equation} \label{eq:hcl}
\Fh_{\cl}^{\ast}=
 \bigoplus_{j \in I} \BC \cl(\Lambda_{j}). 
\end{equation}
Then we define the classical weight lattice 
$P_{\cl}$ to be $\cl(P) \subset \Fh_{\cl}^{\ast}$. 
We have
$P_{\cl}
  \cong \Hom_{\BZ}\bigl(P_{\cl}^{\vee},\BZ\bigr) 
  \cong P/(\BC\delta \cap P)$ 
as (free) $\BZ$-modules, and 
%
%
\begin{equation} \label{eq:lat-cl}
P_{\cl} = 
 \bigoplus_{j \in I} \BZ \Lambda_{j},
\end{equation}
where $\cl(\Lambda_{j}) \in \Fh_{\cl}^{\ast}$ is 
simply denoted by $\Lambda_{j}$ for $j \in I$. 
Further, we set
%
%
\begin{equation} \label{eq:subsets}
\begin{array}{c}
(P_{\cl})_{0}:=
  \bigl\{\mu \in P_{\cl} \mid \mu(c)=0\bigr\}, \qquad 
P_{\cl}^{+}:=
  \sum_{j \in I} \BZ_{\ge 0} \Lambda_{j}, \\[3mm]
(P_{\cl}^{+})_{s}:=
  \bigl\{\mu \in P_{\cl}^{+} \mid \mu(c) = s\bigr\} \quad 
  \text{for each $s \in \BZ_{\ge 0}$}.
\end{array}
\end{equation}
It is easily seen that the quintuplet 
$(A, P_{\cl}, P_{\cl}^{\vee}, \Pi, \Pi^{\vee})$ is 
also a Cartan datum for the affine type 
GCM $A=(a_{ij})_{i,j \in I}$. 
For simplicity, 
a crystal associated to the Cartan datum 
$(A, P_{\cl}, P_{\cl}^{\vee}, \Pi, \Pi^{\vee})$ is 
called a $U_{q}^{\prime}(\Fg)$-crystal, 
where $U_{q}^{\prime}(\Fg)$ denotes 
the $\BC(q)$-subalgebra of 
$U_{q}(\Fg)$ generated by $E_{j}$, $F_{j}$, $j \in I$, and 
$q^{h}$, $h \in P_{\cl}^{\vee}$ (which is 
the quantized universal enveloping algebra of 
$\Fg$ over $\BC(q)$ with weight lattice $P_{\cl}$). 
%
%
%
%
\subsection{Perfect crystals for quantum affine algebras.}
\label{subsec:cry-qaa}
We keep the notation of \S\ref{subsec:qaa}. 
Let us fix a proper subset $J$ of $I$. We set 
$A_{J}:=(a_{ij})_{i,j \in J}$, 
$\Pi_{J}:=
\bigl\{\alpha_{j}\bigr\}_{j \in J} \subset P_{\cl}$, 
$\Pi_{J}^{\vee}:=
\bigl\{h_{j}\bigr\}_{j \in J} \subset P_{\cl}^{\vee}$, 
and denote by $\Fg_{J}$ 
the Lie subalgebra of the affine Lie algebra $\Fg$ 
generated by $E_{j}$, $F_{j}$, $j \in J$, 
and $\Fh_{\cl}$. 
Then the quintuplet 
$(A_{J},P_{\cl},P_{\cl}^{\vee},\Pi_{J},\Pi^{\vee}_{J})$ 
is a Cartan datum for the GCM $A_{J}$. 
For simplicity, a crystal associated to this Cartan datum 
$(A_{J},P_{\cl},P_{\cl}^{\vee},\Pi_{J},\Pi^{\vee}_{J})$ 
is called a $U_{q}(\Fg_{J})$-crystal, where 
$U_{q}(\Fg_{J})$ denotes the $\BC(q)$-subalgebra of 
$U_{q}^{\prime}(\Fg)$ generated by $E_{j}$, $F_{j}$, $j \in J$, 
and $q^{h}$, $h \in P_{\cl}^{\vee}$ 
(which is the quantized universal 
enveloping algebra of $\Fg_{J}$ over $\BC(q)$ 
with weight lattice $P_{\cl}$). 
If $\CB$ is a $U_{q}^{\prime}(\Fg)$-crystal and 
$J$ is a (proper) subset of $I$, then 
the set $\CB$ equipped with 
the Kashiwara operators 
$e_{j}$, $f_{j}$, $j \in J$, and the maps 
$\wt:\CB \rightarrow P_{\cl}$, 
$\ve_{j},\,\vp_{j}:\CB \rightarrow \BZ \cup \{-\infty\}$, 
$j \in J$, is a $U_{q}(\Fg_{J})$-crystal. 
%
%
%
%
\begin{dfn}[{see \cite[\S1.4]{AK}}] \label{dfn:regular}
A $U_{q}^{\prime}(\Fg)$-crystal $\CB$ is said to be regular 
if for every proper subset $J \subsetneq I$, 
the $\CB$, regarded as a $U_{q}(\Fg_{J})$-crystal 
in the way above, is isomorphic to the crystal base of 
an integrable $U_{q}(\Fg_{J})$-module 
(for details about crystal bases, see, 
for example, \cite[Chap.~4, \S4.2]{HK} and \cite[\S4]{Kas3}). 
\end{dfn}

Let $W:=\bigl\langle r_{j} \mid j \in I \bigr\rangle 
\subset \GL(\Fh^{\ast})$ be the Weyl group of $\Fg$, where 
$r_{j} \in \GL(\Fh^{\ast})$ is the simple reflection in 
$\alpha_{j} \in \Fh^{\ast}$. Note that 
the weight lattice $P \subset \Fh^{\ast}$ is 
stable under the action of the Weyl group $W$, and that 
there exists an action of $W$ on $P_{\cl}$ 
induced from that on $P$, since $W\delta=\delta$. 

We can define an action of the Weyl group $W$ on 
a regular $U_{q}^{\prime}(\Fg)$-crystal $\CB$ 
as follows (see \cite[\S7]{Kas2}). For each $j \in I$, 
we define $S_{j}:\CB \rightarrow \CB$ by:
%
%
\begin{equation} \label{eq:si}
S_{j}b=\begin{cases}
f_{j}^{m}b & \text{if \ } m:=(\wt b)(h_{j}) \ge 0 \\[1.5mm]
e_{j}^{-m}b & \text{if \ } m:=(\wt b)(h_{j}) < 0
\end{cases} \qquad \text{for \,} b \in \CB.
\end{equation}
%
%
\begin{prop}
\label{prop:Weyl}
Let $\CB$ be a regular $U_{q}^{\prime}(\Fg)$-crystal. Then, 
there exists a unique action $S:W \rightarrow \Bij(\CB)$, 
$w \mapsto S_{w}$, of the Weyl group $W$ on the set $\CB$ 
such that $S_{r_{j}}=S_{j}$ for all $j \in I$, 
where $\Bij(\CB)$ is the group of all bijections from the set 
$\CB$ to itself. In addition, $\wt(S_{w}b)=w(\wt b)$ holds 
for $w \in W$ and $b \in \CB$. 
\end{prop}
%
%
\begin{dfn}[{see \cite[\S1.4]{AK}}] \label{dfn:extremal}
Let $\CB$ be a regular $U_{q}^{\prime}(\Fg)$-crystal. 
An element $b \in \CB$ is said to be extremal 
(or more accurately, $W$-extremal) 
if for every $w \in W$, either $e_{j}S_{w}b=\theta$ or 
$f_{j}S_{w}b=\theta$ holds for each $j \in I$. 
\end{dfn}
%
%
\begin{rem} \label{rem:extremal}
It immediately follows from the definition above that 
if $b \in \CB$ is an extremal element, 
then $S_{w}b \in \CB$ is an extremal element 
of weight $w\wt b$ for each $w \in W$. 
\end{rem}
We know the following lemma 
from \cite[Lemma~1.6\,(1)]{AK} and its proof 
(see also Lemma~\ref{lem:S2-01} below). 
%
%
\begin{lem} \label{lem:ak}
Let $\CB_{1},\,\CB_{2}$ be 
regular $U_{q}^{\prime}(\Fg)$-crystals of 
finite cardinality such that the weights of 
their elements are all contained 
in $(P_{\cl})_{0}$. 
Let $b_{1} \in \CB_{1}$ and $b_{2} \in \CB_{2}$ be 
extremal elements whose weights are 
contained in the same Weyl chamber with respect to 
the simple coroots $h_{j}$, $j \in I_{0}=I \setminus \{0\}$. 
Then, $b_{1} \otimes b_{2} \in \CB_{1} \otimes \CB_{2}$ 
is an extremal element.
Also, 
$S_{w}(b_{1} \otimes b_{2})=S_{w}b_{1} \otimes S_{w}b_{2}$ 
holds for all $w \in W$.
\end{lem}
%
%
%
%
\begin{dfn} \label{dfn:simple}
A regular $U_{q}^{\prime}(\Fg)$-crystal $\CB$ 
is said to be simple 
if it satisfies the following conditions: 

\vspace{1.5mm}

\noindent (S1) \, 
The set $\CB$ is of finite cardinality, and 
the weights of elements of $\CB$ are 
all contained in $(P_{\cl})_{0}$.

\noindent (S2) \, 
The set of all extremal elements of $\CB$ 
coincides with a Weyl group orbit in $\CB$. 

\noindent (S3) \, 
Let $b \in \CB$ be an extremal element, and set 
$\mu:=\wt b \in (P_{\cl})_{0}$. 
Then the subset $\CB_{\mu} \subset \CB$ of all elements of 
weight $\mu$ consists only of the element $b$, i.e., 
$\CB_{\mu}=\{b\}$.
\end{dfn}
\begin{rem}
Let $\CB$ be a regular $U_{q}^{\prime}(\Fg)$-crystal satisfying 
condition (S1) of Definition~\ref{dfn:simple}. Then 
we see that there exists at least 
one extremal element in $\CB$ (see 
the comment after the proof of 
\cite[Proposition~9.3.2]{Kas2}). 
\end{rem}
%
%
\begin{lem} \label{lem:dom}
Let $\CB$ be a simple $U_{q}^{\prime}(\Fg)$-crystal. Then there 
exists a unique extremal element $u \in \CB$ such that 
$(\wt u)(h_{j}) \ge 0$ for all $j \in I_{0}$.
\end{lem}

\begin{proof}
The existence of an extremal element 
with the desired property 
immediately follows from Remark~\ref{rem:extremal}. 
So, it remains to show the uniqueness. Let $u_{1},\,u_{2}$ 
be extremal elements whose weights are both dominant 
with respect to the simple coroots $h_{j}$, $j \in I_{0}$, 
and set $\mu_{1}:=\wt u_{1}$ and $\mu_{2}:=\wt u_{2}$. 
By condition (S2) of Definition~\ref{dfn:simple}, there exists 
some $w \in W$ such that $S_{w}u_{2}=u_{1}$. Then it follows that 
$\mu_{1}=w\mu_{2} \in W\mu_{2}$. 
We recall from \cite[Proposition~6.5]{Kac} that 
the Weyl group 
$W$ decomposes into 
the semidirect product $W_{I_{0}} 
\ltimes T$ of the Weyl group 
$W_{I_{0}}:=
\langle r_{j} \mid j \in I_{0}\rangle$ 
(of finite type) and the abelian group $T$ of translations. 
Since $\mu_{2} \in P_{\cl}$ is of level zero 
by condition (S1) of Definition~\ref{dfn:simple}, 
it follows from \cite[Chap.~6, formula (6.5.5)]{Kac} that 
$W\mu_{2}=W_{I_{0}}\mu_{2}$, 
and hence $\mu_{1} \in W_{I_{0}}\mu_{2}$. 
Since both of $\mu_{1}$ and $\mu_{2}$ are dominant 
with respect to the simple coroots $h_{j}$, $j \in I_{0}$, 
we deduce that $\mu_{1}=\mu_{2}$. 
Therefore, by condition (S3) of Definition~\ref{dfn:simple}, 
we obtain that $u_{1}=u_{2}$, which shows the uniqueness. 
\end{proof}
%
%
\begin{rem} \label{rem:dom}
Let $\CB$ be a simple $U_{q}^{\prime}(\Fg)$-crystal, and 
let $u \in \CB$ be the unique extremal element 
such that $(\wt u)(h_{j}) \ge 0$ for all $j \in I_{0}$.
We can show by the same argument as in the proof of
\cite[Corollary~5.2]{Kas4}, 
using \cite[Lemma~1.5]{AK}, that 
the weights of elements of $\CB$ are all contained in 
the convex hull of the $W$-orbit of $\wt u$. 
Hence, by Lemma~\ref{lem:dom}, it follows that 
the weights of elements of $\CB$ are 
all contained in the set 
$\wt u - \sum_{j \in I_{0}} \BZ_{\ge 0}\,\alpha_{j}$.
\end{rem}

We know from \cite[Remark~2.5.7]{NS3} that 
the definition of simple 
$U_{q}^{\prime}(\Fg)$-crystals above is 
equivalent to \cite[Definition~1.7]{AK} and 
\cite[Definition~4.9]{Kas4}. Thus we know 
the following from \cite[Lemmas~1.9 and 1.10]{AK}. 
%
%
\begin{prop} \label{prop:simple}
{\rm (1) }
The crystal graph of 
a simple $U_{q}^{\prime}(\Fg)$-crystal is connected. 

\vspace{1.5mm}

\noindent {\rm (2) }
Let $\CB_{1}$, $\CB_{2}$ be simple $U_{q}^{\prime}(\Fg)$-crystals. 
Then, the tensor product $\CB_{1} \otimes \CB_{2}$ 
is also a simple $U_{q}^{\prime}(\Fg)$-crystal.
In particular, the crystal graph of $\CB_{1} \otimes \CB_{2}$ is 
connected. 
\end{prop}
Let $\CB$ be a simple $U_{q}^{\prime}(\Fg)$-crystal. 
We define maps $\ve,\,\vp : \CB \rightarrow P_{\cl}^{+}$ by: 
%
%
\begin{equation} \label{eq:vevp-map}
\ve(b)=\sum_{j \in I} \ve_{j}(b) \Lambda_{j} 
\quad \text{and} \quad
\vp(b)=\sum_{j \in I} \vp_{j}(b) \Lambda_{j} \quad 
\text{for $b \in \CB$}.
\end{equation}
Further, we define a positive integer $\lev \CB$ 
(called the level of $\CB$) and a subset 
$\CB_{\min}$ of $\CB$ by:
%
%
\begin{align} 
& 
\lev \CB = 
 \min \bigl\{ 
 (\ve(b))(c) \mid b \in \CB \bigr\} 
 \in \BZ_{> 0}, \label{eq:level} \\
& \CB_{\min} =
 \bigl\{b \in \CB \mid (\ve(b))(c)=\lev \CB\bigr\} \subset \CB. 
 \label{eq:min-set}
\end{align}
%
%
\begin{rem} \label{rem:min}
It can easily be seen from 
the definition of crystals that 
$\vp(b)-\ve(b)=\wt b$ for every $b \in \CB$.
Since $\wt b \in (P_{\cl})_{0}$ by condition (S1) 
of Definition~\ref{dfn:simple}, 
we have $(\ve(b))(c)=(\vp(b))(c)$
for all $b \in \CB$, and hence $\lev \CB = 
 \min \bigl\{ (\vp(b))(c) \mid 
 b \in \CB \bigr\}$. 
\end{rem}
%
%
%
\begin{dfn} \label{dfn:perfect}
A simple $U_{q}^{\prime}(\Fg)$-crystal $\CB$ 
is said to be perfect 
if the restrictions of the maps 
$\ve,\,\vp: \CB \rightarrow P_{\cl}^{+}$ to 
$\CB_{\min}$ induce bijections 
$\CB_{\min} \rightarrow (P_{\cl}^{+})_{s}$, 
where $s:=\lev \CB$. 
\end{dfn}
%
%
%
\begin{rem} \label{rem:perfect}
(1) \, 
In the definition of 
perfect $U_{q}^{\prime}(\Fg)$-crystals $\CB$, 
it is often required that $\CB$ is isomorphic to 
the crystal base of a finite-dimensional 
$U_{q}^{\prime}(\Fg)$-module as a $U_{q}^{\prime}(\Fg)$-crystal 
(see, for example, Condition~(1) of 
\cite[Definition~10.5.1]{HK}); 
but we do not require it in this paper. 

\vspace{1.5mm}

\noindent (2) \, 
Our definition of 
perfect $U_{q}^{\prime}(\Fg)$-crystals
seems to be slightly different from the ones 
in \cite[Definition~10.5.1]{HK} and 
\cite[\S2.2]{HKOTT}: 
We can deduce from Remark~\ref{rem:dom}, 
condition (S3) of Definition~\ref{dfn:simple}, and 
Proposition~\ref{prop:simple} that 
if $\CB$ is a perfect $U_{q}^{\prime}(\Fg)$-crystal 
in the sense of Definition~\ref{dfn:perfect}, then 
$\CB$ satisfies Conditions (2) -- (5) of 
\cite[Definition~10.5.1]{HK}. But, 
in \cite[\S2.2]{HKOTT}, it is required that 
the perfect $U_{q}^{\prime}(\Fg)$-crystal $\CB$ is 
``finite'' (in the sense of \cite[Definition~2.5]{HKKOT}); 
we do not require this ``finiteness'' condition 
in our definition of perfect 
$U_{q}^{\prime}(\Fg)$-crystals, since 
it does not seem to be essential for our purposes. 

\vspace{1.5mm}

\noindent (3) \, 
If $\CB$ is a perfect $U_{q}^{\prime}(\Fg)$-crystal 
in the sense of Definition~\ref{dfn:perfect}, and 
is isomorphic to the crystal base of a finite-dimensional 
$U_{q}^{\prime}(\Fg)$-module as a $U_{q}^{\prime}(\Fg)$-crystal, 
then $\CB$ is perfect in the sense of \cite[\S2.10]{OSS1}. 
\end{rem}
%
%
\begin{lem} \label{lem:perfect}
Let $\CB_{1}$, $\CB_{2}$ be perfect $U_{q}^{\prime}(\Fg)$-crystals of 
the same level $s$. Then, the tensor product 
$\CB_{1} \otimes \CB_{2}$ is also 
a perfect $U_{q}^{\prime}(\Fg)$-crystal of level $s$. 
\end{lem}
\begin{proof}
We know from Proposition~\ref{prop:simple}\,(2) 
that the tensor product 
$\CB_{1} \otimes \CB_{2}$ is a simple 
$U_{q}^{\prime}(\Fg)$-crystal. In addition, 
by elementary arguments using 
the tensor product rule for crystals, 
we can easily show that 
the $\CB_{1} \otimes \CB_{2}$ is of level $s$, and 
the maps $\ve,\,\vp: 
(\CB_{1} \otimes \CB_{2})_{\min} \rightarrow 
(P_{\cl}^{+})_{s}$ are bijective. 
\end{proof}

We know the following proposition 
from \cite[Theorem~2.4]{OSS1}.
%
%
\begin{prop} \label{prop:perfect}
{\rm (1) }
Let $\CB_{1}$, $\CB_{2}$ be perfect $U_{q}^{\prime}(\Fg)$-crystals 
isomorphic to the crystal bases of finite-dimensional 
$U_{q}^{\prime}(\Fg)$-modules as a $U_{q}^{\prime}(\Fg)$-crystal.
Then, there exists a unique isomorphism 
{\rm(}called a combinatorial $R$-matrix\,{\rm)} 
$R:\CB_{1} \otimes \CB_{2} 
\stackrel{\sim}{\rightarrow} \CB_{2} \otimes \CB_{1}$
of $U_{q}^{\prime}(\Fg)$-crystals. 

\vspace{1.5mm}

\noindent
{\rm (2) }
Let $\CB$ be a perfect $U_{q}^{\prime}(\Fg)$-crystal 
isomorphic to the crystal base of a finite-dimensional 
$U_{q}^{\prime}(\Fg)$-module as a $U_{q}^{\prime}(\Fg)$-crystal.
Then, there exists a $\BZ$-valued function 
$H:\CB \otimes \CB \rightarrow \BZ$ 
{\rm(}called an energy function\,{\rm)} satisfying 
%
%
\begin{equation} \label{eq:energy-e}
H(e_{j}(b_{1} \otimes b_{2})) = 
\begin{cases}
H(b_{1} \otimes b_{2})+1 
 & \text{\rm if $j=0$ and 
   $\vp_{0}(b_{1}) \ge \ve_{0}(b_{2})$,} \\[1.5mm]
H(b_{1} \otimes b_{2})-1 
 & \text{\rm if $j=0$ and 
   $\vp_{0}(b_{1}) < \ve_{0}(b_{2})$,} \\[1.5mm]
H(b_{1} \otimes b_{2})
 & \text{\rm if $j \ne 0$}, 
\end{cases}
\end{equation}
for all $j \in I$ and 
$b_{1} \otimes b_{2} \in \CB \otimes \CB$ such that 
$e_{j}(b_{1} \otimes b_{2}) \ne \theta$. 
\end{prop}
%
%
\begin{rem} \label{rem:energy}
With the notation and assumption of 
Proposition~\ref{prop:perfect}\,(2), 
we have 
%
%
\begin{equation} \label{eq:energy-f}
H(f_{j}(b_{1} \otimes b_{2})) = 
\begin{cases}
H(b_{1} \otimes b_{2})-1 
 & \text{if $j=0$ and 
   $\vp_{0}(b_{1}) > \ve_{0}(b_{2})$,} \\[1.5mm]
H(b_{1} \otimes b_{2})+1 
 & \text{if $j=0$ and 
   $\vp_{0}(b_{1}) \le \ve_{0}(b_{2})$,} \\[1.5mm]
H(b_{1} \otimes b_{2})
 & \text{if $j \ne 0$}, 
\end{cases}
\end{equation}
for all $j \in I$ and 
$b_{1} \otimes b_{2} \in \CB \otimes \CB$ such that 
$f_{j}(b_{1} \otimes b_{2}) \ne \theta$. 
\end{rem}
%
%
%
%
\subsection{A conjectural family of perfect crystals.}
\label{subsec:conj}
In this subsection, let $\Fg$ be either 
a simply-laced affine Lie algebra, or a twisted affine Lie algebra. 
More specifically, let $\Fg$ be the affine Lie algebra of type 
$A_{n}^{(1)} \, (n \ge 2)$, 
$D_{n}^{(1)} \, (n \ge 4)$, 
$E_{6}^{(1)}$, or of type 
$D_{n+1}^{(2)} \, (n \ge 2)$, 
$A_{2n}^{(2)} \, (n \ge 1)$, 
$A_{2n-1}^{(2)} \, (n \ge 3)$, 
$D_{4}^{(3)}$, $E_{6}^{(2)}$, 
with the index set $I$ for the simple roots 
numbered as in \cite[\S4.8, Tables Aff~1 -- Aff~3]{Kac}. 
For each $i \in I_{0}=I \setminus \{0\}=\{1,\,2,\,\dots,\,n\}$, 
$s \in \BZ_{\ge 1}$, and $\zeta \in \BC(q)^{\times}:=
\BC(q) \setminus \{0\}$, we denote by $W^{(i)}_{s}(\zeta)$ the 
finite-dimensional irreducible module over the quantum affine 
algebra $U_{q}^{\prime}(\Fg)$ whose Drinfeld polynomials 
$P_{j}(u) \in \BC(q)[u]$ are specified as follows 
(see \cite[Definition 5.3]{KNT}):

\vspace{3mm}

\noindent
{\bf Case 1:} the case of type 
$A_{n}^{(1)} \, (n \ge 2)$, 
$D_{n}^{(1)} \, (n \ge 4)$, 
$E_{6}^{(1)}$. 
In this case, the Drinfeld polynomials $P_{j}(u)$, 
$j \in I_{0}=\{1,\,2,\,\dots,\,n\}$, are given by: 
\begin{equation*}
P_{j}(u) = 
\begin{cases}
{\displaystyle\prod_{k=1}^{s}} 
(1 - \zeta q^{s+2-2k} u) & \text{if $j = i$,}\\[7mm]
1 & \text{otherwise}.
\end{cases}
\end{equation*}

\vspace{3mm}

\noindent
{\bf Case 2:} the case of type 
$D_{n+1}^{(2)} \, (n \ge 2)$ 
(resp., $A_{2n-1}^{(2)} \, (n \ge 3)$). 
In this case, the Drinfeld polynomials $P_{j}(u)$, 
$j \in I_{0}=\{1,\,2,\,\dots,\,n-1,\,n\}$, are given by: 
\begin{equation*}
P_{j}(u) = 
\begin{cases}
{\displaystyle\prod_{k=1}^{s}} 
(1 - \zeta q^{d_{i}(s+2-2k)} u) & \text{if $j = i$,}\\[7mm]
1 & \text{otherwise}, 
\end{cases}
\end{equation*}
where $d_{i}=1$ if $i = n$ (resp., $i \ne n$), and 
$d_{i}=2$ otherwise. 

\vspace{3mm}

\noindent
{\bf Case 3:} the case of type 
$D_{4}^{(3)}$ (resp., $E_{6}^{(2)}$). 
In this case, the Drinfeld polynomials $P_{j}(u)$, 
$j \in I_{0}=\{1,\,2\}$ 
(resp., $j \in I_{0}=\{1,\,2,\,3,\,4\}$), 
are given by: 
\begin{equation*}
P_{j}(u) = 
\begin{cases}
{\displaystyle\prod_{k=1}^{s}} 
(1 - \zeta q^{d_{i}(s+2-2k)} u) & \text{if $j = i$,}\\[7mm]
1 & \text{otherwise}, 
\end{cases}
\end{equation*}
where $d_{i}=1$ if $i = 1$ (resp., $i=1,\,2$), and 
$d_{i}=3$ (resp., $d_{i}=2$) otherwise.

\vspace{3mm}

\noindent
{\bf Case 4:} the case of type 
$A_{2n}^{(2)} \, (n \ge 1)$. 
In this case, we should remark that 
the index set for the Drinfeld polynomials 
$P_{j}(u)$ is not $I_{0}=\{1,\,2,\,\dots,\,n\}$, 
but $\{0,\,1,\,\dots,\,n-1\}$, and that the Drinfeld 
polynomials $P_{j}(u)$, 
$j \in I_{0}=\{0,\,1,\,\dots,\,n-1\}$, 
are given by: 
\begin{equation*}
P_{j}(u) = 
\begin{cases}
{\displaystyle\prod_{k=1}^{s}} 
(1 - \zeta q^{2(s+2-2k)} u) & \text{if $j = n-i$,}\\[7mm]
1 & \text{otherwise}.
\end{cases}
\end{equation*}

\vspace{3mm}

\noindent
(In all cases above, we used the Drinfeld realization of 
$U_{q}^{\prime}(\Fg)$, and the classification of 
its finite-dimensional irreducible modules of type $1$ by 
Drinfeld polynomials; see \cite{CP1}, \cite{CP2} for details.)
We call this $U_{q}^{\prime}(\Fg)$-module $W^{(i)}_{s}(\zeta)$ 
a Kirillov-Reshetikhin module (KR module for short) over 
$U_{q}^{\prime}(\Fg)$. 

Now, let us fix (arbitrarily) $i \in I_{0}$ and 
$s \in \BZ_{\ge 1}$. 
%
%
\begin{conj}[{cf. \cite[Conjecture~2.1\,(1)]{HKOTT}}]
\label{conj:hkott}
For some $\zeta^{(i)}_{s} \in 
\BC(q)^{\times}=\BC(q) \setminus \{0\}$, 
the KR module $W^{(i)}_{s}(\zeta^{(i)}_{s})$ over $U_{q}^{\prime}(\Fg)$
has a crystal base $\CB^{i,s}$ that is a perfect 
$U_{q}^{\prime}(\Fg)$-crystal of level $s$. 
\end{conj}
%
%
\begin{rem}
Conjecture~\ref{conj:hkott} 
has already been proved in some cases 
(see \cite[Remark~2.3]{HKOTY} and 
\cite[Remark~2.6]{HKOTT}; see also Remark~\ref{rem:ass} below). 
\end{rem}
%
%
\begin{lem} \label{lem:unique}
Assume that Conjecture~\ref{conj:hkott} holds 
for the fixed $i \in I_{0}$ and $s \in \BZ_{\ge 1}$. Let 
$\CL^{i,s} \subset W^{(i)}_{s}(\zeta^{(i)}_{s})$ be 
the crystal lattice corresponding to the crystal base 
$\CB^{i,s}$. Suppose that $(\CL,\CB)$ is another crystal lattice and 
crystal base of $W^{(i)}_{s}(\zeta^{(i)}_{s})$ 
such that the crystal graph of $\CB$ is connected. 
Then, $\CL=f(q) \CL^{i,s}$ holds
for some $f(q) \in \BC(q) \setminus \{0\}$, and 
$\CB$ is isomorphic to 
$\CB^{i,s}$ as a $U_{q}^{\prime}(\Fg)$-crystal. 
Namely, the crystal base of $W^{(i)}_{s}(\zeta^{(i)}_{s})$ is unique, 
up to a nonzero constant multiple. 
\end{lem}

\begin{proof}
Let $b \in \CB^{i,s}$ be an extremal element, and 
set $\mu:=\wt b \in (P_{\cl})_{0}$. 
Note that the set $(\CB^{i,s})_{\mu}$ 
consists only of the element $b$
by condition (S3) of Definition~\ref{dfn:simple}, and 
hence that the $\mu$-weight space of 
$W^{(i)}_{s}(\zeta^{(i)}_{s})$ is one-dimensional.
Let $v \in \CL^{i,s}$ be an element of weight $\mu$ 
corresponding to the $b$ under the canonical projection 
$\CL^{i,s} \twoheadrightarrow 
 \CL^{i,s}/q\CL^{i,s}$. 
Here we remark that the crystal graph of 
$\CB^{i,s}$ is connected by 
Proposition~\ref{prop:simple}\,(1). 
By using Nakayama's lemma, we can show that 
$\CL^{i,s}$ is equal to the $A$-module generated by 
all elements of the form 
$x_{j_{1}} x_{j_{2}} \cdots x_{j_{k}} v$,  
$j_{1},\,j_{2},\,\dots,\,j_{k} \in I$, $k \ge 0$, 
where $A:=\bigl\{f(q) \in \BC(q) \mid 
\text{$f(q)$ is regular at $q=0$}\bigr\}$, and 
$x_{j}$ is either 
the raising Kashiwara operator $e_{j}$ or 
the lowering Kashiwara operator $f_{j}$ 
on $W^{(i)}_{s}(\zeta^{(i)}_{s})$
for each $j \in I$: 
\begin{equation} \label{eq:nakayama01}
\CL^{i,s} = 
 \sum_{
  j_{1},\,j_{2},\,\dots,\,j_{k} \in I; \, k \ge 0
  } 
A \, x_{j_{1}} x_{j_{2}} \cdots x_{j_{k}} v.
\end{equation}
Similarly, take an element 
$v^{\prime} \in \CL$ of weight $\mu$ 
corresponding to a unique element 
$b^{\prime} \in \CB$ of weight $\mu$. 
Because the crystal graph of $\CB$ is connected 
by the assumption, we obtain that 
\begin{equation} \label{eq:nakayama02}
\CL = 
 \sum_{
  j_{1},\,j_{2},\,\dots,\,j_{k} \in I; \, k \ge 0
  } 
A \, x_{j_{1}} x_{j_{2}} \cdots x_{j_{k}} v^{\prime}
\end{equation}
in the same way as above. 

Since the $\mu$-weight space of $W^{(i)}_{s}(\zeta^{(i)}_{s})$ is 
one-dimensional as mentioned above, it follows that 
$v^{\prime}=f(q)v$ for some $f(q) \in \BC(q) \setminus \{0\}$. 
Combining this fact with 
\eqref{eq:nakayama01} and \eqref{eq:nakayama02}, 
we have $\CL=f(q) \CL^{i,s}$. Hence, we have 
a $\BC$-linear isomorphism $\Psi: \CL^{i,s}/q \CL^{i,s}
\stackrel{\sim}{\rightarrow} \CL/q \CL$ induced from 
the transformation $v \mapsto f(q)v$ on 
$W^{(i)}_{s}(\zeta^{(i)}_{s})$. 
Furthermore, we can deduce 
from the connectedness of the crystal bases
$\CB^{i,s}$ and $\CB$ that 
the restriction $\Psi|_{\CB^{i,s}}$ of $\Psi$ to 
$\CB^{i,s}$ gives an isomorphism of 
$U_{q}^{\prime}(\Fg)$-crystals between 
$\CB^{i,s}$ and $\CB$. This proves the lemma. 
\end{proof}
%
%
%
%
\section{Construction of perfect crystals 
for twisted quantum affine algebras.}
\label{sec:main}
From now on throughout this paper, 
let $\Fg=\Fg(A)$ be the affine Lie algebra of 
type $A_{n}^{(1)} \ (n \ge 2)$, $D_{n}^{(1)} \ (n \ge 4)$, 
or $E_{6}^{(1)}$, and 
let $\omega:I \rightarrow I$ be 
a nontrivial diagram automorphism 
satisfying the (additional) condition that $\omega(0)=0$.

%
\subsection{Diagram automorphisms of 
simply-laced affine Lie algebras.}
\label{subsec:diag-aff}
Here we give all pairs $(\Fg,\omega)$ 
of an affine Lie algebra $\Fg$ and 
a nontrivial diagram automorphism 
$\omega:I \rightarrow I$ 
satisfying the condition that $\omega(0)=0$, 
after introducing our numbering of the index set $I$. 
(For the definition of the matrix 
$\ha{A}=(\ha{a}_{ij})_{i,j \in \ha{I}}$, 
see \S\ref{subsec:orbit-aff} below.)

\paragraph{Case (a).}
The affine Cartan matrix $A=(a_{ij})_{i,j \in I}$ 
is of type $A_{2n-1}^{(1)} \, (n \ge 2)$, 
and the diagram automorphism 
$\omega:I \rightarrow I$ is given by:
$\omega(0)=0$ and $\omega(j)=2n-j$ for 
$j \in I_{0}=I \setminus \{0\}$ 
(note that the order of $\omega$ is 2). 
Then the matrix $\ha{A}=(\ha{a}_{ij})_{i,j \in \ha{I}}$ is 
the affine Cartan matrix of type $D_{n+1}^{(2)}$: 

\vspace{1.5mm}

\hspace{10mm}
{\scriptsize
%
%
%
\unitlength 0.1in
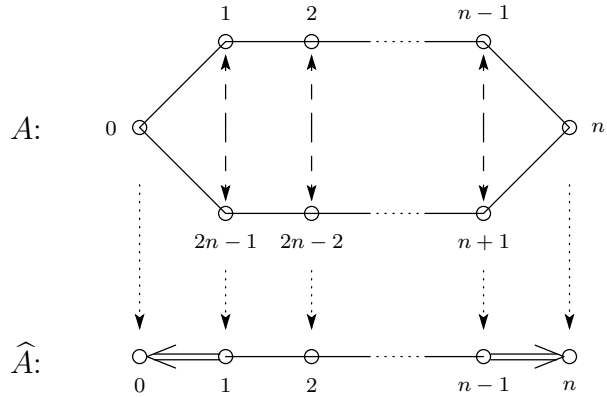
\begin{picture}( 39.1900, 19.5000)( -0.9400,-22.3900)
%
\special{pn 8}%
\special{ar 1988 524 38 38  0.0000000 6.2831853}%
%
\special{pn 8}%
\special{ar 1988 1424 38 38  0.0000000 6.2831853}%
%
\special{pn 8}%
\special{ar 2438 524 38 38  0.0000000 6.2831853}%
%
\special{pn 8}%
\special{ar 2438 1424 38 38  0.0000000 6.2831853}%
%
\special{pn 8}%
\special{ar 3338 524 38 38  0.0000000 6.2831853}%
%
\special{pn 8}%
\special{ar 3338 1424 38 38  0.0000000 6.2831853}%
%
\special{pn 8}%
\special{ar 3788 974 38 38  0.0000000 6.2831853}%
%
\special{pn 8}%
\special{pa 2738 524}%
\special{pa 3038 524}%
\special{dt 0.045}%
%
\special{pn 8}%
\special{pa 2738 1424}%
\special{pa 3038 1424}%
\special{dt 0.045}%
%
\special{pn 8}%
\special{ar 1538 2174 38 38  0.0000000 6.2831853}%
%
\special{pn 8}%
\special{ar 1988 2174 38 38  0.0000000 6.2831853}%
%
\special{pn 8}%
\special{ar 2438 2174 38 38  0.0000000 6.2831853}%
%
\special{pn 8}%
\special{ar 3338 2174 38 38  0.0000000 6.2831853}%
%
\special{pn 8}%
\special{ar 3788 2174 38 38  0.0000000 6.2831853}%
%
\special{pn 8}%
\special{ar 1538 974 38 38  0.0000000 6.2831853}%
%
\special{pn 8}%
\special{pa 3338 524}%
\special{pa 3788 974}%
\special{fp}%
\special{pa 3788 974}%
\special{pa 3338 1424}%
\special{fp}%
%
\special{pn 8}%
\special{pa 1988 524}%
\special{pa 1538 974}%
\special{fp}%
\special{pa 1538 974}%
\special{pa 1988 1424}%
\special{fp}%
%
\special{pn 8}%
\special{pa 1988 524}%
\special{pa 2738 524}%
\special{fp}%
\special{pa 3038 524}%
\special{pa 3338 524}%
\special{fp}%
\special{pa 3338 1424}%
\special{pa 3038 1424}%
\special{fp}%
\special{pa 2738 1424}%
\special{pa 1988 1424}%
\special{fp}%
%
\special{pn 8}%
\special{pa 2738 2174}%
\special{pa 3038 2174}%
\special{dt 0.045}%
%
\special{pn 8}%
\special{pa 1988 2174}%
\special{pa 2738 2174}%
\special{fp}%
\special{pa 3038 2174}%
\special{pa 3338 2174}%
\special{fp}%
%
\special{pn 8}%
\special{pa 1958 2160}%
\special{pa 1614 2160}%
\special{fp}%
%
\special{pn 8}%
\special{pa 1958 2190}%
\special{pa 1614 2190}%
\special{fp}%
%
\special{pn 8}%
\special{pa 1576 2174}%
\special{pa 1726 2114}%
\special{fp}%
%
\special{pn 8}%
\special{pa 1576 2174}%
\special{pa 1726 2234}%
\special{fp}%
%
\special{pn 8}%
\special{pa 3376 2190}%
\special{pa 3720 2190}%
\special{fp}%
%
\special{pn 8}%
\special{pa 3376 2160}%
\special{pa 3720 2160}%
\special{fp}%
%
\special{pn 8}%
\special{pa 3758 2174}%
\special{pa 3608 2234}%
\special{fp}%
%
\special{pn 8}%
\special{pa 3758 2174}%
\special{pa 3608 2114}%
\special{fp}%
%
\special{pn 8}%
\special{pa 1538 1274}%
\special{pa 1538 2024}%
\special{dt 0.045}%
\special{sh 1}%
\special{pa 1538 2024}%
\special{pa 1558 1958}%
\special{pa 1538 1972}%
\special{pa 1518 1958}%
\special{pa 1538 2024}%
\special{fp}%
%
\special{pn 8}%
\special{pa 1988 1724}%
\special{pa 1988 2024}%
\special{dt 0.045}%
\special{sh 1}%
\special{pa 1988 2024}%
\special{pa 2008 1958}%
\special{pa 1988 1972}%
\special{pa 1968 1958}%
\special{pa 1988 2024}%
\special{fp}%
\special{pa 2438 1724}%
\special{pa 2438 2024}%
\special{dt 0.045}%
\special{sh 1}%
\special{pa 2438 2024}%
\special{pa 2458 1958}%
\special{pa 2438 1972}%
\special{pa 2418 1958}%
\special{pa 2438 2024}%
\special{fp}%
\special{pa 3338 1724}%
\special{pa 3338 2024}%
\special{dt 0.045}%
\special{sh 1}%
\special{pa 3338 2024}%
\special{pa 3358 1958}%
\special{pa 3338 1972}%
\special{pa 3318 1958}%
\special{pa 3338 2024}%
\special{fp}%
\special{pa 3788 1274}%
\special{pa 3788 2024}%
\special{dt 0.045}%
\special{sh 1}%
\special{pa 3788 2024}%
\special{pa 3808 1958}%
\special{pa 3788 1972}%
\special{pa 3768 1958}%
\special{pa 3788 2024}%
\special{fp}%
\put(13.8700,-9.7400){\makebox(0,0){$0$}}%
\put(19.8700,-3.7400){\makebox(0,0){$1$}}%
\put(24.3700,-3.7400){\makebox(0,0){$2$}}%
\put(33.3700,-3.7400){\makebox(0,0){$n-1$}}%
\put(39.3700,-9.7400){\makebox(0,0){$n$}}%
\put(33.3700,-15.7400){\makebox(0,0){$n+1$}}%
\put(24.3700,-15.7400){\makebox(0,0){$2n-2$}}%
\put(19.8700,-15.7400){\makebox(0,0){$2n-1$}}%
\put(15.3700,-23.2400){\makebox(0,0){$0$}}%
\put(19.8700,-23.2400){\makebox(0,0){$1$}}%
\put(24.3700,-23.2400){\makebox(0,0){$2$}}%
\put(33.3700,-23.2400){\makebox(0,0){$n-1$}}%
\put(37.8700,-23.2400){\makebox(0,0){$n$}}%
%
\special{pn 8}%
\special{pa 1988 974}%
\special{pa 1988 600}%
\special{da 0.070}%
\special{sh 1}%
\special{pa 1988 600}%
\special{pa 1968 666}%
\special{pa 1988 652}%
\special{pa 2008 666}%
\special{pa 1988 600}%
\special{fp}%
%
\special{pn 8}%
\special{pa 2438 974}%
\special{pa 2438 600}%
\special{da 0.070}%
\special{sh 1}%
\special{pa 2438 600}%
\special{pa 2418 666}%
\special{pa 2438 652}%
\special{pa 2458 666}%
\special{pa 2438 600}%
\special{fp}%
%
\special{pn 8}%
\special{pa 3338 974}%
\special{pa 3338 600}%
\special{da 0.070}%
\special{sh 1}%
\special{pa 3338 600}%
\special{pa 3318 666}%
\special{pa 3338 652}%
\special{pa 3358 666}%
\special{pa 3338 600}%
\special{fp}%
%
\special{pn 8}%
\special{pa 3338 974}%
\special{pa 3338 1350}%
\special{da 0.070}%
\special{sh 1}%
\special{pa 3338 1350}%
\special{pa 3358 1282}%
\special{pa 3338 1296}%
\special{pa 3318 1282}%
\special{pa 3338 1350}%
\special{fp}%
%
\special{pn 8}%
\special{pa 2438 974}%
\special{pa 2438 1350}%
\special{da 0.070}%
\special{sh 1}%
\special{pa 2438 1350}%
\special{pa 2458 1282}%
\special{pa 2438 1296}%
\special{pa 2418 1282}%
\special{pa 2438 1350}%
\special{fp}%
%
\special{pn 8}%
\special{pa 1988 974}%
\special{pa 1988 1350}%
\special{da 0.070}%
\special{sh 1}%
\special{pa 1988 1350}%
\special{pa 2008 1282}%
\special{pa 1988 1296}%
\special{pa 1968 1282}%
\special{pa 1988 1350}%
\special{fp}%
\put(9.4100,-9.7800){\makebox(0,0){{\normalsize $A$:}}}%
\put(9.4100,-21.7800){\makebox(0,0){{\normalsize $\ha{A}$:}}}%
\end{picture}%
}

\paragraph{Case (b).}
The affine Cartan matrix $A=(a_{ij})_{i,j \in I}$ 
is of type $A_{2n}^{(1)} \, (n \ge 1)$, and 
the diagram automorphism 
$\omega:I \rightarrow I$ is given by:
$\omega(0)=0$ and $\omega(j)=2n+1-j$ 
for $j \in I_{0}=I \setminus \{0\}$ 
(note that the order of $\omega$ is 2). 
Then the matrix $\ha{A}=(\ha{a}_{ij})_{i,j \in \ha{I}}$ is 
the affine Cartan matrix of type $A_{2n}^{(2)}$: 

\vspace{4mm}

\hspace{-10mm}
{\scriptsize
%
%
%
\unitlength 0.1in
\begin{picture}( 56.9500, 23.3500)(  3.1100,-25.2300)
%
\special{pn 8}%
\special{ar 2392 804 38 38  0.0000000 6.2831853}%
%
\special{pn 8}%
\special{ar 2392 1704 38 38  0.0000000 6.2831853}%
%
\special{pn 8}%
\special{ar 2842 804 38 38  0.0000000 6.2831853}%
%
\special{pn 8}%
\special{ar 2842 1704 38 38  0.0000000 6.2831853}%
%
\special{pn 8}%
\special{ar 3742 804 38 38  0.0000000 6.2831853}%
%
\special{pn 8}%
\special{ar 3742 1704 38 38  0.0000000 6.2831853}%
%
\special{pn 8}%
\special{pa 3142 804}%
\special{pa 3442 804}%
\special{dt 0.045}%
%
\special{pn 8}%
\special{pa 3142 1704}%
\special{pa 3442 1704}%
\special{dt 0.045}%
%
\special{pn 8}%
\special{ar 1942 2454 38 38  0.0000000 6.2831853}%
%
\special{pn 8}%
\special{ar 2392 2454 38 38  0.0000000 6.2831853}%
%
\special{pn 8}%
\special{ar 2842 2454 38 38  0.0000000 6.2831853}%
%
\special{pn 8}%
\special{ar 3742 2454 38 38  0.0000000 6.2831853}%
%
\special{pn 8}%
\special{ar 4192 2454 38 38  0.0000000 6.2831853}%
%
\special{pn 8}%
\special{ar 1942 1254 38 38  0.0000000 6.2831853}%
%
\special{pn 8}%
\special{pa 2392 804}%
\special{pa 1942 1254}%
\special{fp}%
\special{pa 1942 1254}%
\special{pa 2392 1704}%
\special{fp}%
%
\special{pn 8}%
\special{pa 2392 804}%
\special{pa 3142 804}%
\special{fp}%
\special{pa 3442 804}%
\special{pa 3742 804}%
\special{fp}%
\special{pa 3742 1704}%
\special{pa 3442 1704}%
\special{fp}%
\special{pa 3142 1704}%
\special{pa 2392 1704}%
\special{fp}%
%
\special{pn 8}%
\special{pa 3142 2454}%
\special{pa 3442 2454}%
\special{dt 0.045}%
%
\special{pn 8}%
\special{pa 2392 2454}%
\special{pa 3142 2454}%
\special{fp}%
\special{pa 3442 2454}%
\special{pa 3742 2454}%
\special{fp}%
%
\special{pn 8}%
\special{pa 2362 2440}%
\special{pa 2018 2440}%
\special{fp}%
%
\special{pn 8}%
\special{pa 2362 2470}%
\special{pa 2018 2470}%
\special{fp}%
%
\special{pn 8}%
\special{pa 1980 2454}%
\special{pa 2130 2394}%
\special{fp}%
%
\special{pn 8}%
\special{pa 1980 2454}%
\special{pa 2130 2514}%
\special{fp}%
%
\special{pn 8}%
\special{pa 1942 1554}%
\special{pa 1942 2304}%
\special{dt 0.045}%
\special{sh 1}%
\special{pa 1942 2304}%
\special{pa 1962 2238}%
\special{pa 1942 2252}%
\special{pa 1922 2238}%
\special{pa 1942 2304}%
\special{fp}%
\put(17.9200,-12.5400){\makebox(0,0){$0$}}%
\put(23.9200,-6.5400){\makebox(0,0){$1$}}%
\put(28.4200,-6.5400){\makebox(0,0){$2$}}%
\put(37.4200,-6.5400){\makebox(0,0){$n-1$}}%
\put(37.4200,-18.5400){\makebox(0,0){$n+2$}}%
\put(28.4200,-18.5400){\makebox(0,0){$2n-1$}}%
\put(23.9200,-18.5400){\makebox(0,0){$2n$}}%
\put(19.4200,-26.0400){\makebox(0,0){$0$}}%
\put(23.9200,-26.0400){\makebox(0,0){$1$}}%
\put(28.4200,-26.0400){\makebox(0,0){$2$}}%
\put(37.4200,-26.0400){\makebox(0,0){$n-1$}}%
\put(41.9200,-26.0400){\makebox(0,0){$n$}}%
%
\special{pn 8}%
\special{pa 2392 1254}%
\special{pa 2392 880}%
\special{da 0.070}%
\special{sh 1}%
\special{pa 2392 880}%
\special{pa 2372 946}%
\special{pa 2392 932}%
\special{pa 2412 946}%
\special{pa 2392 880}%
\special{fp}%
%
\special{pn 8}%
\special{pa 2842 1254}%
\special{pa 2842 880}%
\special{da 0.070}%
\special{sh 1}%
\special{pa 2842 880}%
\special{pa 2822 946}%
\special{pa 2842 932}%
\special{pa 2862 946}%
\special{pa 2842 880}%
\special{fp}%
%
\special{pn 8}%
\special{pa 3742 1254}%
\special{pa 3742 880}%
\special{da 0.070}%
\special{sh 1}%
\special{pa 3742 880}%
\special{pa 3722 946}%
\special{pa 3742 932}%
\special{pa 3762 946}%
\special{pa 3742 880}%
\special{fp}%
%
\special{pn 8}%
\special{pa 3742 1254}%
\special{pa 3742 1630}%
\special{da 0.070}%
\special{sh 1}%
\special{pa 3742 1630}%
\special{pa 3762 1562}%
\special{pa 3742 1576}%
\special{pa 3722 1562}%
\special{pa 3742 1630}%
\special{fp}%
%
\special{pn 8}%
\special{pa 2842 1254}%
\special{pa 2842 1630}%
\special{da 0.070}%
\special{sh 1}%
\special{pa 2842 1630}%
\special{pa 2862 1562}%
\special{pa 2842 1576}%
\special{pa 2822 1562}%
\special{pa 2842 1630}%
\special{fp}%
%
\special{pn 8}%
\special{pa 2392 1254}%
\special{pa 2392 1630}%
\special{da 0.070}%
\special{sh 1}%
\special{pa 2392 1630}%
\special{pa 2412 1562}%
\special{pa 2392 1576}%
\special{pa 2372 1562}%
\special{pa 2392 1630}%
\special{fp}%
\put(13.4600,-12.5800){\makebox(0,0){{\normalsize $A$:}}}%
\put(13.4600,-24.5800){\makebox(0,0){{\normalsize $\ha{A}$:}}}%
%
\special{pn 8}%
\special{ar 4196 808 46 46  0.0000000 6.2831853}%
%
\special{pn 8}%
\special{ar 4196 1708 46 46  0.0000000 6.2831853}%
%
\special{pn 8}%
\special{pa 3746 808}%
\special{pa 4196 808}%
\special{fp}%
\special{pa 4196 808}%
\special{pa 4196 1708}%
\special{fp}%
\special{pa 4196 1708}%
\special{pa 3746 1708}%
\special{fp}%
%
\special{pn 8}%
\special{ar 3446 1258 906 906  5.8555537 5.9217789}%
\special{ar 3446 1258 906 906  5.9615140 6.0277391}%
\special{ar 3446 1258 906 906  6.0674742 6.1336994}%
\special{ar 3446 1258 906 906  6.1734345 6.2396597}%
\special{ar 3446 1258 906 906  6.2793948 6.3456199}%
\special{ar 3446 1258 906 906  6.3853550 6.4515802}%
\special{ar 3446 1258 906 906  6.4913153 6.5575405}%
\special{ar 3446 1258 906 906  6.5972756 6.6635007}%
\special{ar 3446 1258 906 906  6.7032358 6.7088078}%
%
\special{pn 8}%
\special{pa 4272 882}%
\special{pa 4250 838}%
\special{fp}%
\special{sh 1}%
\special{pa 4250 838}%
\special{pa 4262 908}%
\special{pa 4274 886}%
\special{pa 4298 890}%
\special{pa 4250 838}%
\special{fp}%
%
\special{pn 8}%
\special{pa 4272 1632}%
\special{pa 4250 1678}%
\special{fp}%
\special{sh 1}%
\special{pa 4250 1678}%
\special{pa 4296 1626}%
\special{pa 4272 1630}%
\special{pa 4260 1610}%
\special{pa 4250 1678}%
\special{fp}%
%
\special{pn 8}%
\special{pa 4166 2442}%
\special{pa 3822 2442}%
\special{fp}%
%
\special{pn 8}%
\special{pa 4166 2472}%
\special{pa 3822 2472}%
\special{fp}%
%
\special{pn 8}%
\special{pa 3784 2458}%
\special{pa 3934 2398}%
\special{fp}%
%
\special{pn 8}%
\special{pa 3784 2458}%
\special{pa 3934 2518}%
\special{fp}%
\put(41.9600,-6.5800){\makebox(0,0){$n$}}%
\put(41.9600,-18.5800){\makebox(0,0){$n+1$}}%
%
\special{pn 8}%
\special{pa 2396 2008}%
\special{pa 2396 2308}%
\special{dt 0.045}%
\special{sh 1}%
\special{pa 2396 2308}%
\special{pa 2416 2242}%
\special{pa 2396 2256}%
\special{pa 2376 2242}%
\special{pa 2396 2308}%
\special{fp}%
%
\special{pn 8}%
\special{pa 2846 2008}%
\special{pa 2846 2308}%
\special{dt 0.045}%
\special{sh 1}%
\special{pa 2846 2308}%
\special{pa 2866 2242}%
\special{pa 2846 2256}%
\special{pa 2826 2242}%
\special{pa 2846 2308}%
\special{fp}%
\special{pa 3746 2008}%
\special{pa 3746 2308}%
\special{dt 0.045}%
\special{sh 1}%
\special{pa 3746 2308}%
\special{pa 3766 2242}%
\special{pa 3746 2256}%
\special{pa 3726 2242}%
\special{pa 3746 2308}%
\special{fp}%
%
\special{pn 8}%
\special{pa 4196 2008}%
\special{pa 4196 2308}%
\special{dt 0.045}%
\special{sh 1}%
\special{pa 4196 2308}%
\special{pa 4216 2242}%
\special{pa 4196 2256}%
\special{pa 4176 2242}%
\special{pa 4196 2308}%
\special{fp}%
\put(14.9200,-3.5800){\makebox(0,0)[lb]{{\normalsize If $n \ge 2$, then}}}%
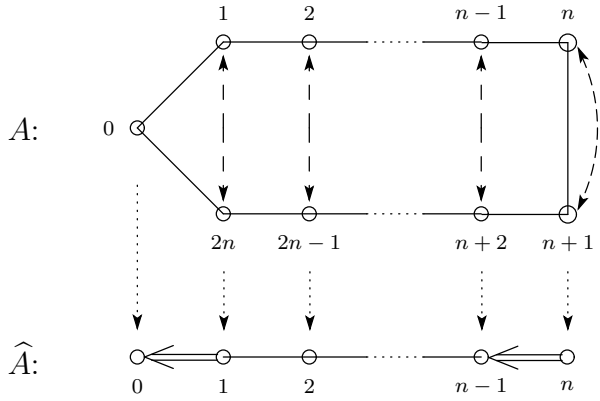
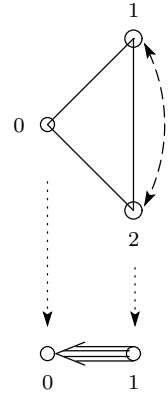
\put(49.4200,-3.5800){\makebox(0,0)[lb]{{\normalsize If $n=1$, then}}}%
%
\special{pn 8}%
\special{pa 5842 808}%
\special{pa 5392 1258}%
\special{fp}%
\special{pa 5392 1258}%
\special{pa 5842 1708}%
\special{fp}%
%
\special{pn 8}%
\special{ar 5392 1258 38 38  0.0000000 6.2831853}%
%
\special{pn 8}%
\special{ar 5100 1258 906 906  5.8555537 5.9217789}%
\special{ar 5100 1258 906 906  5.9615140 6.0277391}%
\special{ar 5100 1258 906 906  6.0674742 6.1336994}%
\special{ar 5100 1258 906 906  6.1734345 6.2396597}%
\special{ar 5100 1258 906 906  6.2793948 6.3456199}%
\special{ar 5100 1258 906 906  6.3853550 6.4515802}%
\special{ar 5100 1258 906 906  6.4913153 6.5575405}%
\special{ar 5100 1258 906 906  6.5972756 6.6635007}%
\special{ar 5100 1258 906 906  6.7032358 6.7088078}%
%
\special{pn 8}%
\special{ar 5842 808 46 46  0.0000000 6.2831853}%
%
\special{pn 8}%
\special{ar 5842 1708 46 46  0.0000000 6.2831853}%
%
\special{pn 8}%
\special{pa 5842 808}%
\special{pa 5842 808}%
\special{fp}%
\special{pa 5842 1708}%
\special{pa 5842 808}%
\special{fp}%
%
\special{pn 8}%
\special{pa 5926 882}%
\special{pa 5902 838}%
\special{fp}%
\special{sh 1}%
\special{pa 5902 838}%
\special{pa 5916 906}%
\special{pa 5928 886}%
\special{pa 5952 888}%
\special{pa 5902 838}%
\special{fp}%
%
\special{pn 8}%
\special{pa 5926 1632}%
\special{pa 5902 1678}%
\special{fp}%
\special{sh 1}%
\special{pa 5902 1678}%
\special{pa 5950 1628}%
\special{pa 5926 1630}%
\special{pa 5914 1610}%
\special{pa 5902 1678}%
\special{fp}%
%
\special{pn 8}%
\special{ar 5392 2458 38 38  0.0000000 6.2831853}%
%
\special{pn 8}%
\special{ar 5842 2458 38 38  0.0000000 6.2831853}%
\put(52.4200,-12.5800){\makebox(0,0){$0$}}%
\put(58.4200,-6.5800){\makebox(0,0){$1$}}%
\put(58.4200,-18.5800){\makebox(0,0){$2$}}%
%
\special{pn 8}%
\special{pa 5850 2008}%
\special{pa 5850 2308}%
\special{dt 0.045}%
\special{sh 1}%
\special{pa 5850 2308}%
\special{pa 5870 2242}%
\special{pa 5850 2256}%
\special{pa 5830 2242}%
\special{pa 5850 2308}%
\special{fp}%
%
\special{pn 8}%
\special{pa 5392 1558}%
\special{pa 5392 2308}%
\special{dt 0.045}%
\special{sh 1}%
\special{pa 5392 2308}%
\special{pa 5412 2242}%
\special{pa 5392 2256}%
\special{pa 5372 2242}%
\special{pa 5392 2308}%
\special{fp}%
%
\special{pn 8}%
\special{pa 5820 2442}%
\special{pa 5476 2442}%
\special{fp}%
%
\special{pn 8}%
\special{pa 5820 2472}%
\special{pa 5476 2472}%
\special{fp}%
%
\special{pn 8}%
\special{pa 5438 2458}%
\special{pa 5588 2398}%
\special{fp}%
%
\special{pn 8}%
\special{pa 5438 2458}%
\special{pa 5588 2518}%
\special{fp}%
%
\special{pn 8}%
\special{pa 5842 2420}%
\special{pa 5542 2420}%
\special{fp}%
%
\special{pn 8}%
\special{pa 5842 2496}%
\special{pa 5542 2496}%
\special{fp}%
\put(53.9200,-26.0800){\makebox(0,0){$0$}}%
\put(58.4200,-26.0800){\makebox(0,0){$1$}}%
\end{picture}%
}

\paragraph{Case (c).}
The affine Cartan matrix $A=(a_{ij})_{i,j \in I}$ 
is of type $D_{n+1}^{(1)} \, (n \ge 3)$, and 
the diagram automorphism 
$\omega:I \rightarrow I$ is given by:
$\omega(j)=j$ for $j \in I \setminus \{n,\,n+1\}$, and 
$\omega(n)=n+1$, $\omega(n+1)=n$ 
(note that the order of $\omega$ is 2). 
Then the matrix $\ha{A}=(\ha{a}_{ij})_{i,j \in \ha{I}}$ is 
the affine Cartan matrix of type $A_{2n-1}^{(2)}$: 

\vspace{3mm}

\hspace{10mm}
{\scriptsize
%
%
\unitlength 0.1in
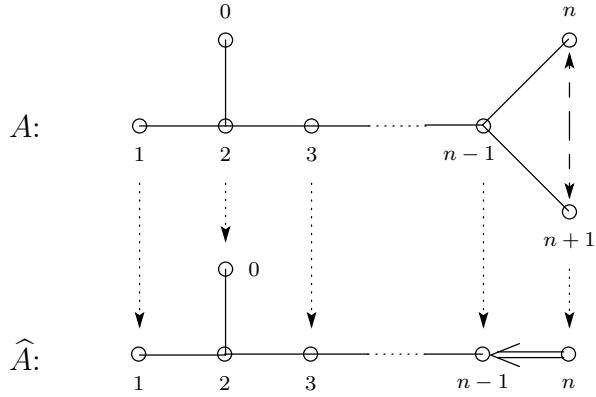
\begin{picture}( 39.2300, 19.4600)( -0.9400,-24.4000)
%
\special{pn 8}%
\special{ar 1538 2376 38 38  0.0000000 6.2831853}%
%
\special{pn 8}%
\special{ar 1988 2376 38 38  0.0000000 6.2831853}%
%
\special{pn 8}%
\special{ar 2438 2376 38 38  0.0000000 6.2831853}%
%
\special{pn 8}%
\special{ar 3338 2376 38 38  0.0000000 6.2831853}%
%
\special{pn 8}%
\special{ar 3788 2376 38 38  0.0000000 6.2831853}%
%
\special{pn 8}%
\special{pa 2738 2376}%
\special{pa 3038 2376}%
\special{dt 0.045}%
%
\special{pn 8}%
\special{pa 1988 2376}%
\special{pa 2738 2376}%
\special{fp}%
\special{pa 3038 2376}%
\special{pa 3338 2376}%
\special{fp}%
\put(15.3700,-25.2500){\makebox(0,0){$1$}}%
\put(19.8700,-25.2500){\makebox(0,0){$2$}}%
\put(24.3700,-25.2500){\makebox(0,0){$3$}}%
\put(33.3700,-25.2500){\makebox(0,0){$n-1$}}%
\put(37.8700,-25.2500){\makebox(0,0){$n$}}%
%
\special{pn 8}%
\special{pa 3762 2364}%
\special{pa 3416 2364}%
\special{fp}%
%
\special{pn 8}%
\special{pa 3762 2394}%
\special{pa 3416 2394}%
\special{fp}%
%
\special{pn 8}%
\special{pa 3380 2380}%
\special{pa 3530 2320}%
\special{fp}%
%
\special{pn 8}%
\special{pa 3380 2380}%
\special{pa 3530 2440}%
\special{fp}%
%
\special{pn 8}%
\special{ar 1542 1180 38 38  0.0000000 6.2831853}%
%
\special{pn 8}%
\special{ar 1992 1180 38 38  0.0000000 6.2831853}%
%
\special{pn 8}%
\special{ar 2442 1180 38 38  0.0000000 6.2831853}%
%
\special{pn 8}%
\special{ar 3342 1180 38 38  0.0000000 6.2831853}%
%
\special{pn 8}%
\special{ar 3792 1630 38 38  0.0000000 6.2831853}%
%
\special{pn 8}%
\special{ar 3792 730 38 38  0.0000000 6.2831853}%
%
\special{pn 8}%
\special{ar 1992 730 38 38  0.0000000 6.2831853}%
%
\special{pn 8}%
\special{pa 1992 2380}%
\special{pa 1542 2380}%
\special{fp}%
\special{pa 1992 1930}%
\special{pa 1992 2380}%
\special{fp}%
%
\special{pn 8}%
\special{ar 1992 1930 38 38  0.0000000 6.2831853}%
%
\special{pn 8}%
\special{pa 1542 1180}%
\special{pa 2742 1180}%
\special{fp}%
\special{pa 1992 1180}%
\special{pa 1992 730}%
\special{fp}%
%
\special{pn 8}%
\special{pa 2742 1180}%
\special{pa 3042 1180}%
\special{dt 0.045}%
%
\special{pn 8}%
\special{pa 3038 1176}%
\special{pa 3338 1176}%
\special{fp}%
\special{pa 3338 1176}%
\special{pa 3788 726}%
\special{fp}%
\special{pa 3338 1176}%
\special{pa 3788 1626}%
\special{fp}%
%
\special{pn 8}%
\special{pa 3792 1180}%
\special{pa 3792 804}%
\special{da 0.070}%
\special{sh 1}%
\special{pa 3792 804}%
\special{pa 3772 872}%
\special{pa 3792 858}%
\special{pa 3812 872}%
\special{pa 3792 804}%
\special{fp}%
\special{pa 3792 1180}%
\special{pa 3792 1554}%
\special{da 0.070}%
\special{sh 1}%
\special{pa 3792 1554}%
\special{pa 3812 1488}%
\special{pa 3792 1502}%
\special{pa 3772 1488}%
\special{pa 3792 1554}%
\special{fp}%
\put(19.9100,-5.7900){\makebox(0,0){$0$}}%
\put(15.4100,-13.2900){\makebox(0,0){$1$}}%
\put(19.9100,-13.2900){\makebox(0,0){$2$}}%
\put(24.4100,-13.2900){\makebox(0,0){$3$}}%
\put(32.6600,-13.2900){\makebox(0,0){$n-1$}}%
\put(37.9100,-5.7900){\makebox(0,0){$n$}}%
\put(37.9100,-17.7900){\makebox(0,0){$n+1$}}%
\put(21.4100,-19.2900){\makebox(0,0){$0$}}%
%
\special{pn 8}%
\special{pa 1542 1480}%
\special{pa 1542 2230}%
\special{dt 0.045}%
\special{sh 1}%
\special{pa 1542 2230}%
\special{pa 1562 2162}%
\special{pa 1542 2176}%
\special{pa 1522 2162}%
\special{pa 1542 2230}%
\special{fp}%
\special{pa 1992 1480}%
\special{pa 1992 1780}%
\special{dt 0.045}%
\special{sh 1}%
\special{pa 1992 1780}%
\special{pa 2012 1712}%
\special{pa 1992 1726}%
\special{pa 1972 1712}%
\special{pa 1992 1780}%
\special{fp}%
\special{pa 2442 1480}%
\special{pa 2442 2230}%
\special{dt 0.045}%
\special{sh 1}%
\special{pa 2442 2230}%
\special{pa 2462 2162}%
\special{pa 2442 2176}%
\special{pa 2422 2162}%
\special{pa 2442 2230}%
\special{fp}%
\special{pa 3342 1480}%
\special{pa 3342 2230}%
\special{dt 0.045}%
\special{sh 1}%
\special{pa 3342 2230}%
\special{pa 3362 2162}%
\special{pa 3342 2176}%
\special{pa 3322 2162}%
\special{pa 3342 2230}%
\special{fp}%
\special{pa 3792 1930}%
\special{pa 3792 2230}%
\special{dt 0.045}%
\special{sh 1}%
\special{pa 3792 2230}%
\special{pa 3812 2162}%
\special{pa 3792 2176}%
\special{pa 3772 2162}%
\special{pa 3792 2230}%
\special{fp}%
\put(9.4100,-11.7900){\makebox(0,0){{\normalsize $A$:}}}%
\put(9.4100,-23.7900){\makebox(0,0){{\normalsize $\ha{A}$:}}}%
\end{picture}%
}

\paragraph{Case (d).}
The affine Cartan matrix $A=(a_{ij})_{i,j \in I}$ 
is of type $D_{4}^{(1)}$, and 
the diagram automorphism 
$\omega:I \rightarrow I$ is given by:
$\omega(0)=0$, $\omega(1)=1$, $\omega(2)=3$, 
$\omega(3)=4$, and $\omega(4)=2$ 
(note that the order of $\omega$ is 3). 
Then the matrix $\ha{A}=(\ha{a}_{ij})_{i,j \in \ha{I}}$ is 
the affine Cartan matrix of type $D_{4}^{(3)}$: 

\vspace{-1.5mm}

\hspace{25mm}
{\scriptsize
%
%
\unitlength 0.1in
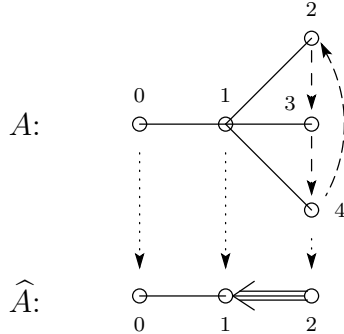
\begin{picture}( 27.0700, 16.6000)(  1.0600,-21.5400)
\put(11.4100,-11.7900){\makebox(0,0){{\normalsize $A$:}}}%
%
\special{pn 8}%
\special{ar 1742 1180 38 38  0.0000000 6.2831853}%
%
\special{pn 8}%
\special{ar 2192 1180 38 38  0.0000000 6.2831853}%
%
\special{pn 8}%
\special{ar 2642 1180 38 38  0.0000000 6.2831853}%
%
\special{pn 8}%
\special{ar 2642 730 38 38  0.0000000 6.2831853}%
%
\special{pn 8}%
\special{ar 2642 1630 38 38  0.0000000 6.2831853}%
%
\special{pn 8}%
\special{pa 1742 1180}%
\special{pa 2192 1180}%
\special{fp}%
\special{pa 2192 1180}%
\special{pa 2642 1180}%
\special{fp}%
\special{pa 2192 1180}%
\special{pa 2642 730}%
\special{fp}%
\special{pa 2192 1180}%
\special{pa 2642 1630}%
\special{fp}%
%
\special{pn 8}%
\special{pa 2642 1246}%
\special{pa 2642 1546}%
\special{da 0.070}%
\special{sh 1}%
\special{pa 2642 1546}%
\special{pa 2662 1480}%
\special{pa 2642 1494}%
\special{pa 2622 1480}%
\special{pa 2642 1546}%
\special{fp}%
%
\special{pn 8}%
\special{pa 2642 796}%
\special{pa 2642 1096}%
\special{da 0.070}%
\special{sh 1}%
\special{pa 2642 1096}%
\special{pa 2662 1030}%
\special{pa 2642 1044}%
\special{pa 2622 1030}%
\special{pa 2642 1096}%
\special{fp}%
%
\special{pn 8}%
\special{ar 2042 1180 772 772  5.7760868 5.8538070}%
\special{ar 2042 1180 772 772  5.9004391 5.9781593}%
\special{ar 2042 1180 772 772  6.0247915 6.1025117}%
\special{ar 2042 1180 772 772  6.1491438 6.2268640}%
\special{ar 2042 1180 772 772  6.2734961 6.3512163}%
\special{ar 2042 1180 772 772  6.3978485 6.4755687}%
\special{ar 2042 1180 772 772  6.5222008 6.5999210}%
\special{ar 2042 1180 772 772  6.6465531 6.7242733}%
\special{ar 2042 1180 772 772  6.7709055 6.7902838}%
%
\special{pn 8}%
\special{pa 2716 804}%
\special{pa 2694 760}%
\special{dt 0.045}%
\special{sh 1}%
\special{pa 2694 760}%
\special{pa 2706 828}%
\special{pa 2718 808}%
\special{pa 2742 810}%
\special{pa 2694 760}%
\special{fp}%
\put(17.4100,-10.2900){\makebox(0,0){$0$}}%
\put(21.9100,-10.2900){\makebox(0,0){$1$}}%
\put(26.4100,-5.7900){\makebox(0,0){$2$}}%
\put(25.2900,-10.6600){\makebox(0,0){$3$}}%
\put(27.9100,-16.2900){\makebox(0,0){$4$}}%
%
\special{pn 8}%
\special{ar 1742 2080 38 38  0.0000000 6.2831853}%
%
\special{pn 8}%
\special{ar 2192 2080 38 38  0.0000000 6.2831853}%
%
\special{pn 8}%
\special{ar 2642 2080 38 38  0.0000000 6.2831853}%
%
\special{pn 8}%
\special{pa 1742 2080}%
\special{pa 2192 2080}%
\special{fp}%
%
\special{pn 8}%
\special{pa 2604 2080}%
\special{pa 2230 2080}%
\special{fp}%
%
\special{pn 8}%
\special{pa 2620 2056}%
\special{pa 2266 2056}%
\special{fp}%
%
\special{pn 8}%
\special{pa 2620 2102}%
\special{pa 2266 2102}%
\special{fp}%
%
\special{pn 8}%
\special{pa 2230 2080}%
\special{pa 2342 2154}%
\special{fp}%
%
\special{pn 8}%
\special{pa 2230 2080}%
\special{pa 2342 2004}%
\special{fp}%
\put(17.4100,-22.2900){\makebox(0,0){$0$}}%
\put(21.9100,-22.2900){\makebox(0,0){$1$}}%
\put(26.4100,-22.2900){\makebox(0,0){$2$}}%
%
\special{pn 8}%
\special{pa 1742 1330}%
\special{pa 1742 1930}%
\special{dt 0.045}%
\special{sh 1}%
\special{pa 1742 1930}%
\special{pa 1762 1862}%
\special{pa 1742 1876}%
\special{pa 1722 1862}%
\special{pa 1742 1930}%
\special{fp}%
%
\special{pn 8}%
\special{pa 2192 1330}%
\special{pa 2192 1930}%
\special{dt 0.045}%
\special{sh 1}%
\special{pa 2192 1930}%
\special{pa 2212 1862}%
\special{pa 2192 1876}%
\special{pa 2172 1862}%
\special{pa 2192 1930}%
\special{fp}%
%
\special{pn 8}%
\special{pa 2642 1780}%
\special{pa 2642 1930}%
\special{dt 0.045}%
\special{sh 1}%
\special{pa 2642 1930}%
\special{pa 2662 1862}%
\special{pa 2642 1876}%
\special{pa 2622 1862}%
\special{pa 2642 1930}%
\special{fp}%
\put(11.4100,-20.7900){\makebox(0,0){{\normalsize $\ha{A}$:}}}%
\end{picture}%
}

\paragraph{Case (e).}
The affine Cartan matrix $A=(a_{ij})_{i,j \in I}$ 
is of type $E_{6}^{(1)}$, and 
the diagram automorphism 
$\omega:I \rightarrow I$ is given by:
$\omega(0)=0$, $\omega(1)=1$, $\omega(2)=2$, 
$\omega(3)=5$, $\omega(4)=6$, $\omega(5)=3$, 
and $\omega(6)=4$ (note that the order of $\omega$ is 2). 
Then the matrix $\ha{A}=(\ha{a}_{ij})_{i,j \in \ha{I}}$ is 
the affine Cartan matrix of type $E_{6}^{(2)}$: 

\vspace{3mm}

\hspace{15mm}
{\scriptsize
%
%
\unitlength 0.1in
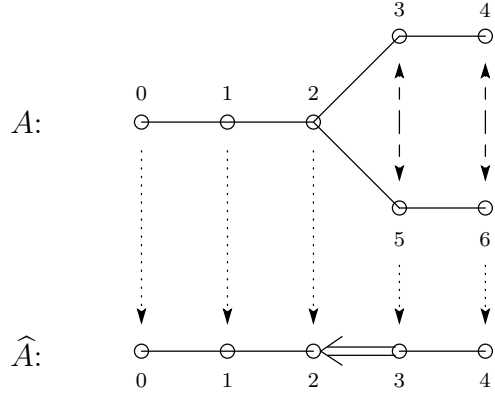
\begin{picture}( 34.7200, 19.5200)(  1.0600,-24.5600)
%
\special{pn 8}%
\special{ar 2192 2382 38 38  0.0000000 6.2831853}%
%
\special{pn 8}%
\special{ar 2642 2382 38 38  0.0000000 6.2831853}%
%
\special{pn 8}%
\special{ar 3092 2382 38 38  0.0000000 6.2831853}%
%
\special{pn 8}%
\special{pa 3070 2360}%
\special{pa 2716 2360}%
\special{fp}%
%
\special{pn 8}%
\special{pa 3070 2404}%
\special{pa 2716 2404}%
\special{fp}%
%
\special{pn 8}%
\special{pa 2680 2382}%
\special{pa 2792 2456}%
\special{fp}%
%
\special{pn 8}%
\special{pa 2680 2382}%
\special{pa 2792 2306}%
\special{fp}%
%
\special{pn 8}%
\special{ar 1742 1182 38 38  0.0000000 6.2831853}%
%
\special{pn 8}%
\special{ar 2192 1182 38 38  0.0000000 6.2831853}%
%
\special{pn 8}%
\special{ar 2642 1182 38 38  0.0000000 6.2831853}%
%
\special{pn 8}%
\special{ar 3092 732 38 38  0.0000000 6.2831853}%
%
\special{pn 8}%
\special{ar 3542 732 38 38  0.0000000 6.2831853}%
%
\special{pn 8}%
\special{ar 3092 1632 38 38  0.0000000 6.2831853}%
%
\special{pn 8}%
\special{ar 3542 1632 38 38  0.0000000 6.2831853}%
%
\special{pn 8}%
\special{pa 1742 1182}%
\special{pa 2642 1182}%
\special{fp}%
\special{pa 2642 1182}%
\special{pa 3092 732}%
\special{fp}%
\special{pa 3092 732}%
\special{pa 3542 732}%
\special{fp}%
\special{pa 2642 1182}%
\special{pa 3092 1632}%
\special{fp}%
\special{pa 3092 1632}%
\special{pa 3542 1632}%
\special{fp}%
%
\special{pn 8}%
\special{ar 1742 2382 38 38  0.0000000 6.2831853}%
%
\special{pn 8}%
\special{ar 3542 2382 38 38  0.0000000 6.2831853}%
%
\special{pn 8}%
\special{pa 1742 2382}%
\special{pa 2642 2382}%
\special{fp}%
\special{pa 3092 2382}%
\special{pa 3542 2382}%
\special{fp}%
\put(17.4100,-10.3100){\makebox(0,0){$0$}}%
\put(21.9100,-10.3100){\makebox(0,0){$1$}}%
\put(26.4100,-10.3100){\makebox(0,0){$2$}}%
\put(30.9100,-5.8900){\makebox(0,0){$3$}}%
\put(35.4100,-5.9600){\makebox(0,0){$4$}}%
\put(30.9100,-17.9600){\makebox(0,0){$5$}}%
\put(35.4100,-17.9600){\makebox(0,0){$6$}}%
%
\special{pn 8}%
\special{pa 3092 1182}%
\special{pa 3092 882}%
\special{da 0.070}%
\special{sh 1}%
\special{pa 3092 882}%
\special{pa 3072 948}%
\special{pa 3092 934}%
\special{pa 3112 948}%
\special{pa 3092 882}%
\special{fp}%
%
\special{pn 8}%
\special{pa 3092 1182}%
\special{pa 3092 1482}%
\special{da 0.070}%
\special{sh 1}%
\special{pa 3092 1482}%
\special{pa 3112 1414}%
\special{pa 3092 1428}%
\special{pa 3072 1414}%
\special{pa 3092 1482}%
\special{fp}%
%
\special{pn 8}%
\special{pa 3542 1182}%
\special{pa 3542 882}%
\special{da 0.070}%
\special{sh 1}%
\special{pa 3542 882}%
\special{pa 3522 948}%
\special{pa 3542 934}%
\special{pa 3562 948}%
\special{pa 3542 882}%
\special{fp}%
%
\special{pn 8}%
\special{pa 3542 1182}%
\special{pa 3542 1482}%
\special{da 0.070}%
\special{sh 1}%
\special{pa 3542 1482}%
\special{pa 3562 1414}%
\special{pa 3542 1428}%
\special{pa 3522 1414}%
\special{pa 3542 1482}%
\special{fp}%
\put(17.4100,-25.3100){\makebox(0,0){$0$}}%
\put(21.9100,-25.3100){\makebox(0,0){$1$}}%
\put(26.4100,-25.3100){\makebox(0,0){$2$}}%
\put(30.9100,-25.3100){\makebox(0,0){$3$}}%
\put(35.4100,-25.3100){\makebox(0,0){$4$}}%
%
\special{pn 8}%
\special{pa 1742 1332}%
\special{pa 1742 2232}%
\special{dt 0.045}%
\special{sh 1}%
\special{pa 1742 2232}%
\special{pa 1762 2164}%
\special{pa 1742 2178}%
\special{pa 1722 2164}%
\special{pa 1742 2232}%
\special{fp}%
%
\special{pn 8}%
\special{pa 2192 1332}%
\special{pa 2192 2232}%
\special{dt 0.045}%
\special{sh 1}%
\special{pa 2192 2232}%
\special{pa 2212 2164}%
\special{pa 2192 2178}%
\special{pa 2172 2164}%
\special{pa 2192 2232}%
\special{fp}%
%
\special{pn 8}%
\special{pa 2642 1332}%
\special{pa 2642 2232}%
\special{dt 0.045}%
\special{sh 1}%
\special{pa 2642 2232}%
\special{pa 2662 2164}%
\special{pa 2642 2178}%
\special{pa 2622 2164}%
\special{pa 2642 2232}%
\special{fp}%
%
\special{pn 8}%
\special{pa 3092 1932}%
\special{pa 3092 2232}%
\special{dt 0.045}%
\special{sh 1}%
\special{pa 3092 2232}%
\special{pa 3112 2164}%
\special{pa 3092 2178}%
\special{pa 3072 2164}%
\special{pa 3092 2232}%
\special{fp}%
%
\special{pn 8}%
\special{pa 3542 1932}%
\special{pa 3542 2232}%
\special{dt 0.045}%
\special{sh 1}%
\special{pa 3542 2232}%
\special{pa 3562 2164}%
\special{pa 3542 2178}%
\special{pa 3522 2164}%
\special{pa 3542 2232}%
\special{fp}%
\put(11.4100,-11.8100){\makebox(0,0){{\normalsize $A$:}}}%
\put(11.4100,-23.8100){\makebox(0,0){{\normalsize $\ha{A}$:}}}%
\end{picture}%
}

\vspace{10mm}

We define $\BC$-linear isomorphisms
$\omega:\Fh \rightarrow \Fh$ and 
$\omega^{\ast}:\Fh^{\ast} \rightarrow \Fh^{\ast}$ by:
\begin{equation}
\begin{array}{l}
\omega(h_{j})=h_{\omega(j)} \quad 
 \text{for $j \in I$}, \quad \omega(d)=d, \\[3mm]
(\omega^{\ast}(\lambda))(h)=\lambda(\omega^{-1}(h)) \quad 
 \text{for $\lambda \in \Fh^{\ast}$ and $h \in \Fh$}.
\end{array}
\end{equation}
It follows that 
the subsets $P^{\vee}$, $\Fh_{\cl}$, and $P_{\cl}^{\vee}$ 
of $\Fh$ are all stable under $\omega \in \GL(\Fh)$, and that 
$P \subset \Fh^{\ast}$ is stable under 
$\omega^{\ast} \in \GL(\Fh^{\ast})$. In addition, 
\begin{equation}
\begin{array}{l}
\omega(c)=c, \qquad 
\omega^{\ast}(\delta)=\delta, \\[2mm]
\omega^{\ast}(\alpha_{j})=\alpha_{\omega(j)}, \quad 
\omega^{\ast}(\Lambda_{j})=\Lambda_{\omega(j)} \quad 
\text{for $j \in I$}. 
\end{array}
\end{equation}
Also, we have a $\BC(q)$-algebra automorphism 
$\omega \in \Aut(U_{q}(\Fg))$ such that 
$\omega(E_{j})=E_{\omega(j)}$, 
$\omega(F_{j})=F_{\omega(j)}$ for $j \in I$, 
and $\omega(q^{h})=q^{\omega(h)}$ 
for $h \in P^{\vee}$. Since $P^{\vee}_{\cl}$ is 
stable under $\omega \in \GL(\Fh)$, 
we see that the $\BC(q)$-subalgebra 
$U_{q}^{\prime}(\Fg)$ is stable under 
$\omega \in \Aut(U_{q}(\Fg))$, thus obtaining 
a $\BC(q)$-algebra automorphism $\omega$ of 
$U_{q}^{\prime}(\Fg)$. 
Further, we define a $\BC$-linear automorphism 
$\omega^{\ast}:\Fh_{\cl}^{\ast} \rightarrow \Fh_{\cl}^{\ast}$ 
by: $\omega^{\ast}(\Lambda_{j})=\Lambda_{\omega(j)}$. 
Note that this $\BC$-linear automorphism of 
$\Fh_{\cl}^{\ast}=(\Fh_{\cl})^{\ast}$ can be thought of 
as the one induced from 
$\omega^{\ast} \in \GL(\Fh^{\ast})$ since $\Fh_{\cl}^{\ast} 
\cong \Fh^{\ast}/\BC\delta$, as well as 
the contragredient map of the restriction of 
$\omega \in \GL(\Fh)$ to $\Fh_{\cl}$. Then we set 
\begin{equation}
(\Fh_{\cl})^{0}:=
 \bigl\{h \in \Fh_{\cl} \mid 
  \omega(h)=h \bigr\}, \qquad 
(\Fh_{\cl}^{\ast})^{0}:=
 \bigl\{\lambda \in \Fh_{\cl}^{\ast} \mid 
  \omega^{\ast}(\lambda)=\lambda \bigr\}.
\end{equation}
%
%
%
%
\subsection{Orbit Lie algebras.}
\label{subsec:orbit-aff}

We choose (and fix) a complete set $\ha{I}$ 
(containing $0 \in I$) of 
representatives of the $\omega$-orbits in $I$ 
in such a way that 
if $j \in \ha{I}$, then $j \le \omega^{k}(j)$ 
for all $k \in \BZ_{\ge 0}$ 
(see the figures in \S\ref{subsec:diag-aff}). 
Now we set
%
%
\begin{equation} \label{eq:c}
c_{ij}:=\sum_{k=0}^{N_{j}-1}a_{i,\,\omega^{k}(j)} \quad 
\text{for $i,\,j \in \ha{I}$}, 
\quad \text{and} \quad 
c_{j}:=c_{jj} \quad \text{for $j \in \ha{I}$},
\end{equation}
where $N_{j}$ is the number of elements of the $\omega$-orbit 
of $j \in \ha{I}$ in $I$. 
(In fact, $N_{j}$ is equal to $1$, $2$, or $3$.)
%
%
\begin{rem}[{cf. \cite[\S2.2]{FSS}}] \label{rem:link} 
We see that $c_{j}=2$ except the case where 
the pair $(\Fg,\omega)$ is in Case (b) and 
$j=n$; if $c_{j}=2$, then the subdiagram 
of the Dynkin diagram of $A$ corresponding to 
the $\omega$-orbit of the $j$ is of type 
$A_{1} \times \dots \times A_{1}$ ($N_{j}$ times). 
On the other hand, if the pair $(\Fg,\omega)$ is 
in Case (b) and $j=n$, then $c_{j}=1$; in this case, 
the subdiagram of the Dynkin diagram of $A$ 
corresponding to 
the $\omega$-orbit of the $j$ is of type $A_{2}$. 
\end{rem}

Further, we set $\ha{a}_{ij}:=2c_{ij}/c_{j}$ 
for $i,\,j \in \ha{I}$. 

\begin{lem}[{see \cite[\S2.2]{FSS}}]
The matrix $\ha{A}:=(\ha{a}_{ij})_{i,j \in \ha{I}}$ 
is a generalized Cartan matrix of {\rm(}twisted\,{\rm)} affine type. 
Moreover, the explicit type of the {\rm GCM} $\ha{A}$ is 
as in \S\ref{subsec:diag-aff}. 
\end{lem}

Let $\ha{\Fg}:=\Fg(\ha{A})$ be the (affine) Kac-Moody algebra 
over $\BC$ associated to the GCM $\ha{A}$ above, 
which is called the orbit Lie algebra 
(corresponding to the diagram automorphism $\omega$). 
Then, $\ha{\Fh}=
\left(\bigoplus_{j \in \ha{I}} \BC\,\ha{h}_{j}\right) \oplus \BC\,\ha{d}$ 
is a Cartan subalgebra of $\ha{\Fg}$, with 
$\ha{\Pi}^{\vee}:=\{\ha{h}_{j}\}_{j \in \ha{I}}$ 
the set of simple coroots, 
and $\ha{d}$ 
the scaling element. 
Denote by 
$\ha{\Pi}:=\bigl\{\ha{\alpha}_{j}\bigr\}_{j \in \ha{I}} 
\subset \ha{\Fh}^{\ast}:=(\ha{\Fh})^{\ast}$ 
the set of simple roots, and 
$\ha{\Lambda}_{j} \in \ha{\Fh}^{\ast}$, $j \in \ha{I}$, 
the fundamental weights for the orbit Lie algebra $\ha{\Fg}$ 
(of affine type); 
note that $\ha{\alpha}_{j}(\ha{d})=\delta_{j,0}$ and 
$\ha{\Lambda}_{j}(\ha{d})=0$ for $j \in \ha{I}$. 
Let 
\begin{equation}
\ha{\delta}:=
\sum_{j \in \ha{I}} \ha{a}_{j}\,\ha{\alpha}_{j} 
\in \ha{\Fh}^{\ast}
\qquad \text{and} \qquad 
\ha{c}:=
\sum_{j \in \ha{I}} \ha{a}_{j}^{\vee}\,\ha{h}_{j} 
\in \ha{\Fh} 
\end{equation}
be  the null root and the canonical central element of
$\ha{\Fg}$, respectively. We take a dual weight lattice 
$\ha{P}^{\vee} \subset \ha{\Fh}$ and a weight lattice 
$\ha{P} \subset \ha{\Fh}^{\ast}$ as follows:
\begin{equation}
\ha{P}^{\vee}=
\left(\bigoplus_{j \in \ha{I}} 
\BZ\,\ha{h}_{j}\right) \oplus \BZ\,\ha{d} \, 
\subset \ha{\Fh}, 
\qquad 
\ha{P}=
\left(\bigoplus_{j \in \ha{I}} \BZ\,\ha{\Lambda}_{j}\right)
\oplus 
\BZ\left(\dfrac{1}{\ha{a}_{0}}\ha{\delta}\right) 
\subset \ha{\Fh}^{\ast}.
\end{equation}
Define $\ha{\Fh}_{\cl}, \, 
\ha{P}_{\cl}^{\vee} \subset \ha{\Fh}$, 
and $\ha{P}_{\cl} \subset 
\ha{\Fh}^{\ast}_{\cl}:=(\ha{\Fh}_{\cl})^{\ast}$ 
for the orbit Lie algebra $\ha{\Fg}$, and also define 
subsets $(\ha{P}_{\cl})_{0}$, $\ha{P}_{\cl}^{+}$, and 
$(\ha{P}_{\cl}^{+})_{s}$, $s \in \BZ_{\ge 0}$, of 
$\ha{P}_{\cl}$ as in \S\ref{subsec:qaa}. 
Note that $(\ha{A}, \ha{P}_{\cl}, \ha{P}_{\cl}^{\vee}, 
\ha{\Pi}, \ha{\Pi}^{\vee})$ is a Cartan datum 
for the GCM $\ha{A}$. Let $U_{q}(\ha{\Fg})$ be 
the quantized universal enveloping algebra 
of the orbit Lie algebra $\ha{\Fg}$ over $\BC(q)$ with 
weight lattice $\ha{P}$, and define 
its $\BC(q)$-subalgebra $U_{q}^{\prime}(\ha{\Fg})$ 
as in \S\ref{subsec:qaa} (which is the quantized 
universal enveloping algebra of $\ha{\Fg}$ over $\BC(q)$ 
with weight lattice $\ha{P}_{\cl}$). 
We call a crystal associated to the Cartan datum 
$(\ha{A}, \ha{P}_{\cl}, \ha{P}_{\cl}^{\vee}, 
\ha{\Pi}, \ha{\Pi}^{\vee})$
a $U_{q}^{\prime}(\ha{\Fg})$-crystal. 
Further, for a proper subset $\ha{J}$ of $\ha{I}$, 
let us define the Lie subalgebra 
$\ha{\Fg}_{\ha{J}}$ of $\ha{\Fg}$, and 
the $\BC(q)$-subalgebra $U_{q}(\ha{\Fg}_{\ha{J}})$ 
of $U_{q}^{\prime}(\ha{\Fg})$ corresponding to 
the subset $\ha{J}$ as in \S\ref{subsec:cry-qaa}. 
We call a crystal associated to the Cartan datum 
$(\ha{A}_{\ha{J}}, \ha{P}_{\cl}, \ha{P}_{\cl}^{\vee}, 
\ha{\Pi}_{\ha{J}}, \ha{\Pi}_{\ha{J}}^{\vee})$ 
for the GCM $\ha{A}_{\ha{J}}:=(\ha{a}_{ij})_{i,j \in \ha{J}}$ 
a $U_{q}(\ha{\Fg}_{\ha{J}})$-crystal, 
where 
$\ha{\Pi}_{\ha{J}}:=
 \bigl\{\ha{\alpha}_{j}\bigr\}_{j \in \ha{J}} 
 \subset \ha{P}_{\cl}$ and 
$\ha{\Pi}_{\ha{J}}^{\vee}:=
 \bigl\{\ha{h}_{j}\bigr\}_{j \in \ha{J}}
 \subset \ha{P}_{\cl}^{\vee}$. 

Now, let us define $\BC$-linear isomorphisms
$P_{\omega}:
 (\Fh_{\cl})^{0} \rightarrow \ha{\Fh}_{\cl}$ 
from the fixed point subspace $(\Fh_{\cl})^{0}$
onto $\ha{\Fh}_{\cl}=\bigoplus_{j \in \ha{I}}\BC \ha{h}_{j}$, 
and 
$P_{\omega}^{\ast}:
 \ha{\Fh}_{\cl}^{\ast} \rightarrow (\Fh_{\cl}^{\ast})^{0}$ 
from $\ha{\Fh}_{\cl}^{\ast}=
\bigoplus_{j \in \ha{I}}\BC \ha{\Lambda}_{j}$ 
onto the fixed point subspace $(\Fh_{\cl}^{\ast})^{0}$ by: 
%
%
\begin{equation} \label{eq:pos-aff} 
P_{\omega}
  \left(
    \dfrac{1}{N_{j}}
    {\displaystyle \sum_{k=0}^{N_{j}-1}} 
    h_{\omega^{k}(j)}
  \right)
 =\ha{h}_{j} 
\quad \text{and} \quad 
\pos(\ha{\Lambda}_{j})=
  {\displaystyle\sum_{k=0}^{N_{j}-1}} 
  \Lambda_{\omega^{k}(j)} 
\quad \text{for $j \in \ha{I}$}
\end{equation}
(here note that $(\Fh_{\cl}^{0})^{\ast}$ can be 
identified with $(\Fh^{\ast}_{\cl})^{0}$ 
in natural way). Then it is easily seen that 
$(P_{\omega}^{\ast}(\ha{\lambda}))(h)=
\ha{\lambda}(P_{\omega}(h))$ 
for $\ha{\lambda} \in \ha{\Fh}_{\cl}^{\ast}$ and 
$h \in (\Fh_{\cl})^{0}$, and that 
%
%
\begin{equation} \label{eq:pos-aff01}
P_{\omega}(c)=\ha{c},
\qquad 
\pos(\ha{\alpha}_{j})=
  \dfrac{2}{c_{j}}
  {\displaystyle\sum_{k=0}^{N_{j}-1}} 
  \alpha_{\omega^{k}(j)}
\quad \text{for $j \in \ha{I}$.}
\end{equation}
Furthermore, we can identify the Weyl group 
$\ha{W}:=\langle \ha{r}_{j} \mid j \in \ha{I} \rangle$ 
of the orbit Lie algebra $\ha{\Fg}$ with the subgroup 
$\ti{W}:=\bigl\{w \in W \mid 
\omega^{\ast}w=w\omega^{\ast}\bigr\}$ of the 
Weyl group $W=\langle r_{j} \mid j \in I \rangle$ of $\Fg$ 
as follows. Define $w_{j} \in W$ by:
%
%
\begin{equation} \label{eq:wj}
w_{j}=
\begin{cases}
  r_{j}\,r_{\omega(j)}\,r_{j}
  & \text{if $c_{j}=1$}, \\[3mm]
  r_{j}\,r_{\omega(j)} \cdots r_{\omega^{N_{j}-1}(j)}
  & \text{if $c_{j}=2$},
\end{cases}
\end{equation}
for each $j \in \ha{I}$ (see Remark~\ref{rem:link}). 
Then it follows that $w_{j} \in \ti{W}$ for all $j \in \ha{I}$. 
Also, we see from \cite[\S3]{FRS} that 
there exists a group isomorphism 
$\Theta: \ha{W} \rightarrow \ti{W}$ such that 
$\Theta(\ha{w})|_{\sw}=
\pos \circ \ha{w} \circ (\pos)^{-1}$ for 
each $\ha{w} \in \ha{W}$, and $\Theta(\ha{r}_{j})=w_{j}$ 
for all $j \in \ha{I}$.
%
%
%
%
\subsection{Fixed point subsets.}
\label{subsec:perfixed}
Let $\ha{I} \subset I$ be 
the index set (chosen as in \S\ref{subsec:orbit-aff}) 
for the orbit Lie algebra $\ha{\Fg}$  
corresponding to the $\omega$. 
Let us fix (arbitrarily) 
$i \in \ha{I}_{0}:=\ha{I} \setminus \{0\}$ and 
$s \in \BZ_{\ge 1}$. 
For the rest of this section, we make 
the following assumption (cf. Conjecture~\ref{conj:hkott}):
%
%
\begin{ass} \label{ass}
There exists some $\zeta^{(i)}_{s} \in \BC(q)^{\times}$ 
(independent of $0 \le k \le N_{i}-1$) 
such that for every $0 \le k \le N_{i}-1$, the KR module 
$W^{(\omega^{k}(i))}_{s}(\zeta^{(i)}_{s})$ 
over $U_{q}^{\prime}(\Fg)$ has 
a crystal base, denoted by $\CB^{\omega^{k}(i),s}$.
Further, the $\CB^{\omega^{k}(i),s}$, 
$0 \le k \le N_{i}-1$, are all perfect 
$U_{q}^{\prime}(\Fg)$-crystals of level $s$.
\end{ass}
%
%
\begin{rem} \label{rem:dom02}
Let $0 \le k \le N_{i}-1$. Then, 
since the $U_{q}^{\prime}(\Fg)$-crystal $\CB^{\omega^{k}(i),s}$ 
is perfect (and hence simple) by Assumption~\ref{ass}, 
it follows from Lemma~\ref{lem:dom} that there exists a unique 
extremal element of $\CB^{\omega^{k}(i),s}$, denoted by 
$u_{\omega^{k}(i),s}$, such that 
$(\wt u_{\omega^{k}(i),s})(h_{j}) \ge 0$ for all $j \in I_{0}$. 
In addition, we can show that 
$\wt u_{\omega^{k}(i),s}=s\vpi_{\omega^{k}(i)}$, where 
$\vpi_{\omega^{k}(i)}:=\Lambda_{\omega^{k}(i)}-
a_{\omega^{k}(i)}^{\vee}\Lambda_{0} \in P_{\cl}$. 
\end{rem}
%
%
%
\begin{rem}[{see \cite[Remark~2.3]{HKOTY}}] \label{rem:ass}
For Cases (a) and (b), we know from \cite{KMN} that 
Assumption~\ref{ass} is satisfied for 
all $i \in \ha{I}_{0}$ and $s \in \BZ_{\ge 1}$. 
For Case~(c) (resp., Case~(d)), 
we know from \cite{KMN} that 
Assumption~\ref{ass} is satisfied 
if $i=1,\,n$ (resp., $i=2$), 
and from \cite{Ko} that 
Assumption~\ref{ass} is satisfied 
if $i \ne n$ and $s=1$ (resp., if $s=1$). 
\end{rem}

First, we define a bijection
$\tau_{\omega}:
\CB^{i,s} \rightarrow 
\CB^{\omega(i),s}$ such that 
$\tau_{\omega} \circ e_{j}=
 e_{\omega(j)} \circ \tau_{\omega}$ 
and 
$\tau_{\omega} \circ f_{j}=
 f_{\omega(j)} \circ \tau_{\omega}$ for all $j \in I$
($\tau_{\omega}(\theta)$ is understood to be $\theta$), 
and such that 
$\wt(\tau_{\omega}(b))=\omega^{\ast}(\wt b)$
for each $b \in \CB^{i,s}$ 
as follows. 
Let $\rho:U_{q}^{\prime}(\Fg) \rightarrow 
\End_{\BC(q)}(W^{(i)}_{s}(\zeta^{(i)}_{s}))$ be 
the representation map affording 
the KR module $W^{(i)}_{s}(\zeta^{(i)}_{s})$ 
over $U_{q}^{\prime}(\Fg)$. 
It immediately follows that the 
the representation of $U_{q}^{\prime}(\Fg)$ 
on the (same) $\BC(q)$-vector space $W^{(i)}_{s}(\zeta^{(i)}_{s})$ 
given by $\rho \circ \omega^{-1}$, denoted by 
$(\rho \circ \omega^{-1}, W^{(i)}_{s}(\zeta^{(i)}_{s}))$, is 
finite-dimensional and irreducible. 
In addition, we can easily check that if 
$P_{j}(u) \in \BC(q)[u]$, $j \in I_{0}$, are 
the Drinfeld polynomials
of the KR module $W^{(i)}_{s}(\zeta^{(i)}_{s})$ 
over $U_{q}^{\prime}(\Fg)$ (see \S\ref{subsec:conj}), then 
the Drinfeld polynomials $P_{j}^{\omega}(u) \in \BC(q)[u]$, 
$j \in I_{0}$, of the representation 
$(\rho \circ \omega^{-1}, W^{(i)}_{s}(\zeta^{(i)}_{s}))$
of $U_{q}^{\prime}(\Fg)$ are given by: 
$P_{j}^{\omega}(u)=P_{\omega^{-1}(j)}(u)$ for each $j \in I_{0}$. 
(Here we have used Assumption~\ref{ass} that the 
$\zeta^{(i)}_{s} \in \BC(q)^{\times}$ is independent of 
$0 \le k \le N_{i}-1$.)
Because the finite-dimensional irreducible 
$U_{q}^{\prime}(\Fg)$-modules (of type $1$) are parametrized 
by their Drinfeld polynomials up to 
$U_{q}^{\prime}(\Fg)$-module isomorphism, 
it follows that the representation 
$(\rho \circ \omega^{-1}, W^{(i)}_{s}(\zeta^{(i)}_{s}))$ 
of $U_{q}^{\prime}(\Fg)$ is equivalent to 
the KR module $W^{(\omega(i))}_{s}(\zeta^{(i)}_{s})$ 
over $U_{q}^{\prime}(\Fg)$.
We denote by $\tau_{\omega}:
W^{(i)}_{s}(\zeta^{(i)}_{s}) \rightarrow 
W^{(\omega(i))}_{s}(\zeta^{(i)}_{s})$ an intertwining map 
between these two representations of $U_{q}^{\prime}(\Fg)$. 
Namely, $\tau_{\omega}:
W^{(i)}_{s}(\zeta^{(i)}_{s}) \rightarrow 
W^{(\omega(i))}_{s}(\zeta^{(i)}_{s})$ denotes 
a $\BC(q)$-linear isomorphism such that 
%
%
\begin{equation} \label{eq:tau-omega}
\tau_{\omega}(x v) = \omega(x) \tau_{\omega}(v) \quad 
\text{for $x \in U_{q}^{\prime}(\Fg)$ and 
$v \in W^{(i)}_{s}(\zeta^{(i)}_{s})$}. 
\end{equation}
It immediately follows from \eqref{eq:tau-omega} that 
for each $\mu \in P_{\cl}$, 
the $\mu$-weight space of 
$W^{(i)}_{s}(\zeta^{(i)}_{s})$ is sent to the 
$\omega^{\ast}(\mu)$-weight space of 
$W^{(\omega(i))}_{s}(\zeta^{(i)}_{s})$ under $\tau_{\omega}$. 
Also, we can easily deduce that 
%
%
\begin{equation} \label{eq:tau-omega01}
\tau_{\omega} \circ e_{j}=
e_{\omega(j)} \circ \tau_{\omega}, 
\quad \text{and} \quad
\tau_{\omega} \circ f_{j}=
f_{\omega(j)} \circ \tau_{\omega} \quad 
\text{for all $j \in I$}, 
\end{equation}
where $e_{j}$ (resp., $f_{j}$), $j \in I$, denote
the raising (resp., lowering) Kashiwara 
operators on $W^{(i)}_{s}(\zeta^{(i)}_{s})$, and 
also those on $W^{(\omega(i))}_{s}(\zeta^{(i)}_{s})$. 
Let us denote by $\CL^{i,s}$ the crystal lattice of 
$W^{(i)}_{s}(\zeta^{(i)}_{s})$. By \eqref{eq:tau-omega01}, 
we see that the image $\tau_{\omega}(\CL^{i,s})$ 
of the $\CL^{i,s}$ is stable under 
the action of the Kashiwara operators 
$e_{j}$ and $f_{j}$, $j \in I$, 
on $W^{(\omega(i))}_{s}(\zeta^{(i)}_{s})$. 
Hence these Kashiwara operators 
on $W^{(\omega(i))}_{s}(\zeta^{(i)}_{s})$ 
induces operators, denoted also 
by $e_{j}$ and $f_{j}$, $j \in I$, 
on the $\BC$-vector space 
$\tau_{\omega}(\CL^{i,s})/
 q\tau_{\omega}(\CL^{i,s})$. 
If $\tau_{\omega}:\CL^{i,s}/
 q\CL^{i,s} \rightarrow 
 \tau_{\omega}(\CL^{i,s})/
 q\tau_{\omega}(\CL^{i,s})$ denotes 
the induced $\BC$-linear map, then 
it follows from \eqref{eq:tau-omega01} that 
the set $\tau_{\omega}(\CB^{i,s}) \cup \{0\}$ is 
stable under the Kashiwara operators 
$e_{j}$ and $f_{j}$, $j \in I$, on 
$\tau_{\omega}(\CL^{i,s})/
 q\tau_{\omega}(\CL^{i,s})$, 
which means that 
$(\tau_{\omega}(\CL^{i,s}), 
  \tau_{\omega}(\CB^{i,s}))$ is 
a crystal base of $W^{(\omega(i))}_{s}(\zeta^{(i)}_{s})$. 
Therefore, it follows from Lemma~\ref{lem:unique} that 
$\tau_{\omega}(\CB^{i,s}) \cong 
 \CB^{\omega(i),s}$ as 
$U_{q}^{\prime}(\Fg)$-crystals. 
Thus we have obtained a bijection 
$\tau_{\omega}:
 \CB^{i,s} \rightarrow 
 \CB^{\omega(i),s}$ 
such that 
%
%
\begin{equation} \label{eq:tau-omega02}
\begin{array}{c}
\tau_{\omega} \circ e_{j}=
e_{\omega(j)} \circ \tau_{\omega}, 
\quad \text{and} \quad
\tau_{\omega} \circ f_{j}=
f_{\omega(j)} \circ \tau_{\omega} \quad 
 \text{for all $j \in I$}, \\[3mm]
\wt(\tau_{\omega}(b))=\omega^{\ast}(\wt b) \quad 
 \text{for each $b \in \CB^{i,s}$}.
\end{array}
\end{equation}
Here (and below) we understand that 
$\tau_{\omega}(\theta)=\theta$. 
Similarly, for each $1 \le k \le N_{i}-1$, 
we obtain a bijection 
$\tau_{\omega}:\CB^{\omega^{k}(i),s} 
 \rightarrow 
 \CB^{\omega^{k+1}(i),s}$ such that 
$\tau_{\omega} \circ e_{j}=
 e_{\omega(j)} \circ \tau_{\omega}$ and 
$\tau_{\omega} \circ f_{j}=
f_{\omega(j)} \circ \tau_{\omega}$ 
for all $j \in I$, and such that 
$\wt(\tau_{\omega}(b))=\omega^{\ast}(\wt b)$ 
for each $b \in \CB^{\omega^{k}(i),s}$.

Next, we set
%
%
\begin{equation} \label{eq:tiB}
\ti{\CB}^{i,s}:=
\begin{cases}
\CB^{i,s} 
 & \text{if $N_{i}=1$}, \\[1.5mm]
\CB^{i,s} \otimes 
\CB^{\omega(i),s} 
 & \text{if $N_{i}=2$}, \\[1.5mm]
\CB^{i,s} \otimes 
\CB^{\omega(i),s} \otimes 
\CB^{\omega^{2}(i),s}
 & \text{if $N_{i}=3$}.
\end{cases}
\end{equation}
Since the $\CB^{\omega^{k}(i),s}$, 
$0 \le k \le N_{i}-1$, are perfect 
$U_{q}^{\prime}(\Fg)$-crystals of (the same) level $s$ 
by Assumption~\ref{ass}, 
it follows from Lemma~\ref{lem:perfect} that 
$\ti{\CB}^{i,s}$ is a perfect 
$U_{q}^{\prime}(\Fg)$-crystal of level $s$.
We define an action of the diagram automorphism $\omega$ 
on $\ti{\CB}^{i,s}$ as follows. If $N_{i}=1$, then 
$\omega:\ti{\CB}^{i,s} \rightarrow \ti{\CB}^{i,s}$ 
is defined to be $\tau_{\omega}$. 
If $N_{i}=2$, then we first define a bijection from 
$\CB^{i,s} \otimes 
 \CB^{\omega(i),s}$
onto
$\CB^{\omega(i),s} \otimes 
 \CB^{i,s}$ by: 
$b_{1} \otimes b_{2} \mapsto 
 \tau_{\omega}(b_{1}) \otimes 
 \tau_{\omega}(b_{2})$ 
for $b_{1} \otimes b_{2} \in 
\CB^{i,s} \otimes \CB^{\omega(i),s}$. 
By Proposition~\ref{prop:perfect}\,(1), 
we have an isomorphism (a combinatorial $R$-matrix) 
$\CB^{\omega(i),s} \otimes 
 \CB^{i,s} 
 \stackrel{\sim}{\rightarrow}
 \CB^{i,s} \otimes 
 \CB^{\omega(i),s}$ 
of $U_{q}^{\prime}(\Fg)$-crystals. 
We now define $\omega:\ti{\CB}^{i,s} \rightarrow 
\ti{\CB}^{i,s}$ to be the composition of these maps:
\begin{equation}
\omega: \ti{\CB}^{i,s}=
\CB^{i,s} \otimes \CB^{\omega(i),s} 
  \stackrel{\tau_{\omega} \otimes \tau_{\omega}}
           {\longrightarrow}
   \CB^{\omega(i),s} \otimes \CB^{i,s} 
  \stackrel{\sim}{\rightarrow} 
   \CB^{i,s} \otimes \CB^{\omega(i),s}=
  \ti{\CB}^{i,s}.
\end{equation}
Similarly, if $N_{i}=3$, then we define 
an action of $\omega$ on $\ti{\CB}^{i,s}$ to be 
the composition of the map 
$\tau_{\omega} \otimes \tau_{\omega} \otimes \tau_{\omega}$ 
with combinatorial $R$-matrices:
\begin{align}
\omega: \ti{\CB}^{i,s}=
\CB^{i,s} \otimes 
\CB^{\omega(i),s} \otimes 
\CB^{\omega^{2}(i),s} 
& \stackrel{
     \tau_{\omega} \otimes 
     \tau_{\omega} \otimes 
     \tau_{\omega}}{\longrightarrow}
  \CB^{\omega(i),s} \otimes 
  \CB^{\omega^{2}(i),s} \otimes 
  \CB^{i,s} \nonumber \\
& \stackrel{\sim}{\rightarrow} 
  \CB^{i,s} \otimes 
  \CB^{\omega(i),s} \otimes 
  \CB^{\omega^{2}(i),s}=\ti{\CB}^{i,s}. 
\end{align}
In all cases above, we can deduce from 
the tensor product rule for crystals, 
\eqref{eq:tau-omega02}, and the comment after 
\eqref{eq:tau-omega02} that 
%
%
\begin{equation} \label{eq:omega}
\begin{array}{c}
\omega \circ e_{j}=
e_{\omega(j)} \circ \omega 
\quad \text{and} \quad
\omega \circ f_{j}=
f_{\omega(j)} \circ \omega \quad
 \text{ on $\ti{\CB}^{i,s}$ for all $j \in I$}, \\[3mm]
\wt(\omega(b))=\omega^{\ast}(\wt b) \quad
 \text{for each $b \in \ti{\CB}^{i,s}$}, 
\end{array}
\end{equation}
where $\omega(\theta)$ is understood to be $\theta$. 
Finally, we set 
\begin{equation}
\ha{\CB}^{i,s}:=
 \bigl\{ b \in \ti{\CB}^{i,s} \mid \omega(b)=b\bigr\}. 
\end{equation}
Note that the weights of elements of 
$\ha{\CB}^{i,s}$ are all contained in 
$(P_{\cl})_{0} \cap \sw$ by condition 
(S1) of Definition~\ref{dfn:simple} 
and the second equality of \eqref{eq:omega}. 
%
%
%
%
\subsection{Main result.}
\label{subsec:main}
For each $j \in \ha{I}$, 
we define $\omega$-Kashiwara operators 
$\ti{e}_{j}$ and $\ti{f}_{j}$ on 
$\ti{\CB}^{i,s} \cup \{\theta\}$ by: 
%
%
\begin{equation} \label{eq:o-kas-op}
\ti{x}_{j}=
  \begin{cases}
       x_{j} \, 
       x_{\omega(j)}^{2} \,
       x_{j}
    & \text{if \ } c_{j}=1, \\[3mm]
       x_{j}\,x_{\omega(j)} \cdots 
       x_{\omega^{N_{j}-1}(j)}
    & \text{if \ } c_{j}=2,
  \end{cases}
\end{equation}
where $x$ is either $e$ or $f$. 
The main result of this paper is 
the following theorem.
%
%
\begin{thm} \label{thm:main}
Let $i \in \ha{I}_{0}=\ha{I} \setminus \{0\}$ and 
$s \in \BZ_{\ge 1}$ 
{\rm(}fixed as in \S\ref{subsec:perfixed}\,{\rm)}. 
We keep Assumption~\ref{ass}. 
Then, the subset $\ha{\CB}^{i,s} \cup \{\theta\}$ 
of $\ti{\CB}^{i,s} \cup \{\theta\}$ is 
stable under the $\omega$-Kashiwara operators 
$\ti{e}_{j}$ and $\ti{f}_{j}$ on 
$\ti{\CB}^{i,s} \cup \{\theta\}$ 
for all $j \in \ha{I}$. Moreover, if we set 
%
%
\begin{equation} \label{eq:main}
\begin{cases}
\ha{\wt}\,b:=(\pos)^{-1}(\wt b) \in \ha{P}_{\cl}
 & \text{\rm for $b \in \ha{\CB}^{i,s}$}, \\[2mm]
\ha{\ve}_{j}(b):=
 \max\bigl\{m \ge 0 \mid 
 (\ti{e}_{j})^{m} b \ne \theta\bigr\} 
 & \text{\rm for $b \in \ha{\CB}^{i,s}$ and 
         $j \in \ha{I}$}, \\[2mm]
\ha{\vp}_{j}(b):=
 \max\bigl\{m \ge 0 \mid 
 (\ti{f}_{j})^{m} b \ne \theta \bigr\} 
 & \text{\rm for $b \in \ha{\CB}^{i,s}$ and 
                 $j \in \ha{I}$},
\end{cases}
\end{equation}
then the set $\ha{\CB}^{i,s}$ equipped with 
the $\omega$-Kashiwara operators $\ti{e}_{j}$ 
and $\ti{f}_{j}$, $j \in \ha{I}$, the maps 
$\ha{\wt}\,:\ha{\CB}^{i,s} \rightarrow \ha{P}_{\cl}$, 
and $\ha{\ve}_{j},\,\ha{\vp}_{j}:
\ha{\CB}^{i,s} \rightarrow \BZ_{\ge 0}$, $j \in \ha{I}$, 
becomes a perfect $U_{q}^{\prime}(\ha{\Fg})$-crystal of level $s$.  
\end{thm}

We will establish Theorem~\ref{thm:main} 
under the following plan. 
First, in \S\ref{subsec:fixed-reg}, we show 
that $\ha{\CB}^{i,s}$ is 
a regular $U_{q}^{\prime}(\ha{\Fg})$-crystal. 
Next, in \S\ref{subsec:prf-simple}, we prove that 
the $U_{q}^{\prime}(\ha{\Fg})$-crystal $\ha{\CB}^{i,s}$ is 
simple. 
Finally, in \S\ref{subsec:prf-bij}, 
we show that the level of 
$\ha{\CB}^{i,s}$ is equal to $s$, and that 
the restrictions of the maps $\ha{\ve}$, $\ha{\vp}$ 
to $(\ha{\CB}^{i,s})_{\min}$ induce bijections 
$(\ha{\CB}^{i,s})_{\min} \rightarrow (\ha{P}_{\cl}^{+})_{s}$, 
where the maps $\ha{\ve}, \ha{\vp} : \ha{\CB}^{i,s} 
\rightarrow \ha{P}_{\cl}^{+}$ are defined as 
in \eqref{eq:vevp-map}, and the set 
$(\ha{\CB}^{i,s})_{\min}$ is defined as in 
\eqref{eq:min-set}.
%
%
%
%
%
\section{Fixed point subsets of 
crystals under the action of $\omega$.}
\label{sec:pre}
Let $\Fg=\Fg(A)$ be the affine Lie algebra of 
type $A_{n}^{(1)} \ (n \ge 2)$, $D_{n}^{(1)} \ (n \ge 4)$, 
or $E_{6}^{(1)}$, and 
let $\omega:I \rightarrow I$ be a nontrivial diagram automorphism 
satisfying the condition that $\omega(0)=0$. 
%
%
\subsection{Fixed point subsets of crystal bases.}
\label{subsec:fixed}
Let us fix a proper subset $J$ of $I$ such that $\omega(J)=J$. 
For an integral weight $\lambda \in P_{\cl}$ 
that is dominant with respect to the simple coroots 
$h_{j}$, $j \in J$, which we call a $J$-dominant 
integral weight, we denote by $V_{J}(\lambda)$ 
the integrable highest weight $U_{q}(\Fg_{J})$-module 
of highest weight $\lambda$. 
Further, let us denote by $\CB_{J}(\lambda)$ 
the crystal base of $V_{J}(\lambda)$ with 
raising Kashiwara operators $e_{j}$, $j \in J$, and 
lowering Kashiwara operators $f_{j}$, $j \in J$. 

Let us take (and fix) a $J$-dominant integral weight 
$\lambda \in P_{\cl}$ such that 
$\omega^{\ast}(\lambda)=\lambda$. 
Then, as in \cite[\S3.2]{NS1}, we obtain an action 
$\omega:\CB_{J}(\lambda) \rightarrow \CB_{J}(\lambda)$ 
of the diagram automorphism $\omega$ on the crystal base 
$\CB_{J}(\lambda)$ satisfying the condition: 
%
%
\begin{equation} \label{eq:omega02}
\begin{array}{c}
\omega \circ e_{j}=
e_{\omega(j)} \circ \omega 
\quad \text{and} \quad
\omega \circ f_{j}=
f_{\omega(j)} \circ \omega \quad
 \text{for all $j \in I$}, \\[3mm]
\wt(\omega(b))=\omega^{\ast}(\wt b) \quad
 \text{for each $b \in \CB_{J}(\lambda)$}, 
\end{array}
\end{equation}
where $\omega(\theta)$ is understood to be $\theta$. 
We set 
\begin{equation}
\CB_{J}^{\omega}(\lambda):=\bigl\{b \in \CB_{J}(\lambda) \mid 
\omega(b)=b\bigr\}. 
\end{equation}
It immediately follows from \eqref{eq:omega02} that 
$\wt b \in P_{\cl} \cap \sw$ 
for all $b \in \CB_{J}^{\omega}(\lambda)$. 
Set $\ha{J}:=J \cap \ha{I} \subsetneq \ha{I}$.
For an integral weight 
$\ha{\lambda} \in \ha{P}_{\cl}$ 
that is dominant with 
respect to the simple coroots 
$\ha{h}_{j}$, $j \in \ha{J}$, 
which we call a $\ha{J}$-dominant integral weight, 
we denote by $\ha{V}_{\ha{J}}(\ha{\lambda})$ 
the integrable highest weight 
$U_{q}(\ha{\Fg}_{\ha{J}})$-module of 
highest weight $\ha{\lambda}$. 
Also, $\ha{\CB}_{\ha{J}}(\ha{\lambda})$ 
denotes the crystal base of 
$\ha{V}_{\ha{J}}(\ha{\lambda})$. 
Further, for each $j \in \ha{J}$, 
we define the $\omega$-Kashiwara operators 
$\ti{e}_{j}$ and $\ti{f}_{j}$ 
on $\CB_{J}(\lambda) \cup \{\theta\}$ by:
%
%
\begin{equation} \label{eq:o-kas-op01}
\ti{x}_{j}=
  \begin{cases}
       x_{j} \, 
       x_{\omega(j)}^{2} \,
       x_{j}
    & \text{if \ } c_{j}=1, \\[3mm]
       x_{j}\,x_{\omega(j)} \cdots 
       x_{\omega^{N_{j}-1}(j)}
    & \text{if \ } c_{j}=2,
  \end{cases}
\end{equation}
where $x$ is either $e$ or $f$.

We know the following theorem from 
\cite[Theorem~2.2.1\,(1) -- (3)]{NS2}; 
note that the restriction $\omega|_{J}$ of $\omega$ to the 
subset $J$ of $I$ is a diagram automorphism for 
the finite-dimensional, reductive 
Lie subalgebra $\Fg_{J}$ of $\Fg$, and 
the Lie subalgebra $\ha{\Fg}_{\ha{J}}$ of $\ha{\Fg}$ 
can be thought of as the orbit Lie algebra of 
$\Fg_{J}$ corresponding to the $\omega|_{J}$. 
%
%
%
\begin{thm} \label{thm:fixed}
Let $\lambda \in P_{\cl}$ be a $J$-dominant integral weight 
such that $\omega^{\ast}(\lambda)=\lambda$, and set 
$\ha{\lambda}:=(\pos)^{-1}(\lambda)$. Then, 
the subset $\CB_{J}^{\omega}(\lambda) \cup \{\theta\}$ 
of $\CB_{J}(\lambda) \cup \{\theta\}$ is 
stable under the $\omega$-Kashiwara operators 
$\ti{e}_{j}$ and $\ti{f}_{j}$ 
on $\CB_{J}(\lambda) \cup \{\theta\}$ 
for all $j \in \ha{J}$.
Moreover, the set $\CB_{J}^{\omega}(\lambda)$ equipped
with the $\omega$-Kashiwara operators 
$\ti{e}_{j}$, $\ti{f}_{j}$, $j \in \ha{J}$, and 
the maps
%
%
\begin{equation} \label{eq:fixed}
\begin{cases}
\ha{\wt}\,b:=(\pos)^{-1}(\wt b) 
 \in \ha{P}_{\cl}
 & \text{\rm for $b \in \CB_{J}^{\omega}(\lambda)$}, \\[2mm]
\ha{\ve}_{j}(b):=
 \max\bigl\{m \ge 0 \mid 
  (\ti{e}_{j})^{m} b \ne \theta\bigr\} 
 \in \BZ_{\ge 0}
 & \text{\rm for $b \in \CB_{J}^{\omega}(\lambda)$ and 
         $j \in \ha{J}$}, \\[2mm]
\ha{\vp}_{j}(b):=
 \max\bigl\{m \ge 0 \mid 
  (\ti{f}_{j})^{m} b \ne \theta\bigr\} 
 \in \BZ_{\ge 0}
 & \text{\rm for $b \in \CB_{J}^{\omega}(\lambda)$ and 
         $j \in \ha{J}$},
\end{cases}
\end{equation}
becomes a $U_{q}(\ha{\Fg}_{\ha{J}})$-crystal 
isomorphic to the crystal base $\ha{\CB}_{\ha{J}}(\ha{\lambda})$ of 
the integrable highest weight $U_{q}(\ha{\Fg}_{\ha{J}})$-module 
$\ha{V}_{\ha{J}}(\ha{\lambda})$ of highest weight $\ha{\lambda}$. 
\end{thm}

For each $m \in \BZ_{\ge 1}$ and $j \in \ha{I}$, 
we define operators $\ti{e}(m)_{j}$ and $\ti{f}(m)_{j}$ 
on $\CB(\lambda) \cup \{\theta\}$ by: 
%
%
\begin{equation} \label{eq:o-kas-op-n}
\ti{x}(m)_{j}=
  \begin{cases}
       x_{j}^{m} \, 
       x_{\omega(j)}^{2m} \,
       x_{j}^{m}
    & \text{if \ } c_{j}=1, \\[3mm]
       x_{j}^{m}\,
       x_{\omega(j)}^{m} \cdots 
       x_{\omega^{N_{j}-1}(j)}^{m}
    & \text{if \ } c_{j}=2,
  \end{cases}
\end{equation}
where $x$ is either $e$ or $f$. 
We know the following from 
\cite[Theorem~2.2.1\,(4)]{NS2}. 
%
%
\begin{prop} \label{prop:o-kas-op-n}
Let $\lambda \in P_{\cl}$ be a $J$-dominant integral weight 
such that $\omega^{\ast}(\lambda)=\lambda$. 
Then, for every $m \in \BZ_{\ge 1}$ and $j \in \ha{J}$, 
we have $\ti{e}(m)_{j}=(\ti{e}_{j})^{m}$ and 
$\ti{f}(m)_{j}=(\ti{f}_{j})^{m}$ 
on $\CB^{\omega}_{J}(\lambda) \cup \{\theta\}$. 
\end{prop}
%
%
\begin{prop} \label{prop:vevp}
Let $\lambda \in P_{\cl}$ be a $J$-dominant integral weight 
such that $\omega^{\ast}(\lambda)=\lambda$. 
Then, for each $b \in \CB^{\omega}_{J}(\lambda)$, 
we have 
$\ha{\ve}_{j}(b)=\ve_{\omega^{k}(j)}(b)$ and 
$\ha{\vp}_{j}(b)=\vp_{\omega^{k}(j)}(b)$ 
for all $j \in \ha{J}$ and $0 \le k \le N_{j}-1$. 
\end{prop}

\begin{proof}
Let $b \in \CB^{\omega}_{J}(\lambda)$. Then 
we see from \eqref{eq:omega02} that 
$\ve_{\omega^{k}(j)}(b)=\ve_{j}(b)$ and 
$\vp_{\omega^{k}(j)}(b)=\vp_{j}(b)$ 
for all $j \in \ha{J}$ and $0 \le k \le N_{j}-1$. 
So, we need only show that 
$\ha{\ve}_{j}(b)=\ve_{j}(b)$ and 
$\ha{\vp}_{j}(b)=\vp_{j}(b)$. But, these equalities 
follow from \cite[Lemma~2.1.3 and Theorem~2.2.1\,(2)]{NS2}. 
\end{proof}
%
%
%
%
\subsection{Fixed point subsets of regular crystals.}
\label{subsec:fixed-reg}
Let $\CB$ be a regular $U_{q}^{\prime}(\Fg)$-crystal 
with an action $\omega:\CB \rightarrow \CB$ of 
the diagram automorphism $\omega$ satisfying the condition: 
%
%
\begin{equation} \label{eq:omega03}
\begin{array}{c}
\omega \circ e_{j}=
e_{\omega(j)} \circ \omega 
\quad \text{and} \quad
\omega \circ f_{j}=
f_{\omega(j)} \circ \omega \quad
 \text{for all $j \in I$}, \\[3mm]
\wt(\omega(b))=\omega^{\ast}(\wt b) \quad
 \text{for each $b \in \CB$}. 
\end{array}
\end{equation}
Here (and below) we understand that $\omega(\theta)=\theta$. 
We set 
\begin{equation}
\CB^{\omega}:=\bigl\{b \in \CB \mid \omega(b)=b\bigr\}, 
\end{equation}
and assume that $\CB^{\omega} \ne \emptyset$. 
Note that $\wt b \in P_{\cl} \cap \sw$ 
for all $b \in \CB^{\omega}$ 
by the second equality of \eqref{eq:omega03}. 
For each $j \in \ha{I}$, 
we define $\omega$-Kashiwara operators 
$\ti{e}_{j}$ and $\ti{f}_{j}$ 
on $\CB \cup \{\theta\}$ by: 
%
%
\begin{equation} \label{eq:o-kas-op02}
\ti{x}_{j}=
  \begin{cases}
       x_{j} \, 
       x_{\omega(j)}^{2} \,
       x_{j}
    & \text{if \ } c_{j}=1, \\[3mm]
       x_{j}\,x_{\omega(j)} \cdots 
       x_{\omega^{N_{j}-1}(j)}
    & \text{if \ } c_{j}=2,
  \end{cases}
\end{equation}
where $x$ is either $e$ or $f$. 
Further, we define maps 
$\ha{\wt}:\CB^{\omega} \rightarrow P_{\cl}$ and 
$\ha{\ve}_{j},\,\ha{\vp}_{j}: 
\CB^{\omega} \rightarrow \BZ_{\ge 0}$, 
$j \in \ha{I}$, by: 
%
%
\begin{equation} \label{eq:wtvevp}
\begin{cases}
\ha{\wt}\,b:=(\pos)^{-1}(\wt b) 
 \in \ha{P}_{\cl}
 & \text{for $b \in \CB^{\omega}$}, \\[2mm]
\ha{\ve}_{j}(b):=
 \max\bigl\{m \ge 0 \mid 
  (\ti{e}_{j})^{m} b \ne \theta\bigr\} 
 \in \BZ_{\ge 0}
 & \text{for $b \in \CB^{\omega}$ and 
         $j \in \ha{I}$}, \\[2mm]
\ha{\vp}_{j}(b):=
 \max\bigl\{m \ge 0 \mid 
  (\ti{f}_{j})^{m} b \ne \theta\bigr\} 
 \in \BZ_{\ge 0}
 & \text{for $b \in \CB^{\omega}$ and 
         $j \in \ha{I}$}.
\end{cases}
\end{equation}
%
%
\begin{prop} \label{prop:fixed-reg}
Let $\CB$ be a regular $U_{q}^{\prime}(\Fg)$-crystal 
with an action $\omega:\CB \rightarrow \CB$ 
of the diagram automorphism $\omega$ 
satisfying \eqref{eq:omega03}. 
Then, the subset $\CB^{\omega} \cup \{\theta\}$ 
of $\CB \cup \{\theta\}$ is stable 
under the $\omega$-Kashiwara operators $\ti{e}_{j}$
and $\ti{f}_{j}$ on $\CB \cup \{\theta\}$ 
for all $j \in \ha{I}$. 
Moreover, the the fixed point subset $\CB^{\omega}$ 
equipped with the maps 
$\ha{\wt}$, $\ti{e}_{j}$, $\ti{f}_{j}$, 
$j \in \ha{I}$, and $\ha{\ve}_{j}$, $\ha{\vp}_{j}$, 
$j \in \ha{I}$, becomes a regular 
$U_{q}^{\prime}(\ha{\Fg})$-crystal. 
\end{prop} 
\begin{proof}
{\bf Step 1.}
For a proper subset $\ha{J}$ of $\ha{I}$, we set 
$J:=\bigl\{ \omega^{k}(j) \mid 
 0 \le k \le N_{j}-1,\, j \in \ha{J}
 \bigr\} \subsetneq I$. 
Since $\CB$ is 
a regular $U_{q}^{\prime}(\Fg)$-crystal, it follows that 
$\CB$ is isomorphic as a $U_{q}(\Fg_{J})$-crystal to 
the crystal base of 
an integrable $U_{q}(\Fg_{J})$-module, and hence 
to a direct sum of the crystal bases of integrable 
highest weight $U_{q}(\Fg_{J})$-modules. 
Namely, there exists an isomorphism of 
$U_{q}(\Fg_{J})$-crystals:
%
%
\begin{equation} \label{eq:reg01}
\Psi_{J}:\CB \stackrel{\sim}{\rightarrow}
 \CB_{J}(\lambda_{1}) \sqcup \CB_{J}(\lambda_{2}) 
 \sqcup \cdots \sqcup 
 \CB_{J}(\lambda_{p}), 
\end{equation}
for some $J$-dominant integral weights 
$\lambda_{1},\,\lambda_{2},\,\dots,\,\lambda_{p} \in P_{\cl}$. 
Put $\CB_{t}:=\Psi_{J}^{-1}(\CB_{J}(\lambda_{t}))$, and 
$b_{t}:=\Psi_{J}^{-1}(v_{\lambda_{t}})$ for $1 \le t \le p$, 
where $v_{\lambda_{t}}$ is the highest weight element of 
$\CB_{J}(\lambda_{t})$. 
Assume that 
$\CB^{\omega} \cap \CB_{t} \ne \emptyset$
for any $1 \le t \le p^{\prime}$, and 
$\CB^{\omega} \cap \CB_{t} = \emptyset$
for all $p^{\prime}+1 \le t \le p$ ($1 \le p^{\prime} \le p$). 
Then the highest weight elements $b_{t} \in \CB_{t}$, 
$1 \le t \le p^{\prime}$, are all fixed by 
$\omega:\CB \rightarrow \CB$, i.e.,  
$\omega(b_{t})=b_{t}$ 
for all $1 \le t \le p^{\prime}$. 
Indeed, if $1 \le t \le p^{\prime}$ and 
$b \in \CB^{\omega} \cap \CB_{t} \ne \emptyset$, 
then there exist $j_{1},\,j_{2},\,\dots,\,j_{l} \in J$ 
such that $b_{t}=e_{j_{1}} e_{j_{2}} \cdots e_{j_{l}} b$. 
So, it follows from \eqref{eq:omega03} that 
$\omega(b_{t})=
e_{\omega(j_{1})}
e_{\omega(j_{2})} \cdots 
e_{\omega(j_{l})} b$, 
since $\omega(b)=b$. 
Here we note that $\omega(j_{1}),\,\omega(j_{2}),\,\dots,\,
\omega(j_{l}) \in J$, since $J$ is stable under $\omega$. 
Thus, because $\CB_{t}$ is a connected component of 
$\CB$ regarded as a $U_{q}(\Fg_{J})$-crystal, 
it follows that $\omega(b_{t})=
e_{\omega(j_{1})}
e_{\omega(j_{2})} \cdots 
e_{\omega(j_{l})} b$ 
is also contained in $\CB_{t}$. 
In addition, we see from \eqref{eq:omega03} 
that $\omega(b_{t}) \in \CB_{t}$ is the highest weight 
element with respect to $e_{j}$, $j \in J$. 
Therefore, we conclude that $\omega(b_{t})=b_{t}$ by 
the uniqueness of the highest weight element in $\CB_{t}$, 
which is isomorphic to $\CB_{J}(\lambda_{t})$ as a 
$U_{q}(\Fg_{J})$-crystal, and hence that 
$\omega^{\ast}(\lambda_{t})=\lambda_{t}$. 
Note also that $\omega(\CB_{t})=\CB_{t}$ since 
$\CB_{t}$ is a connected $U_{q}(\Fg_{J})$-crystal. 
Moreover, since $\CB_{t}$ is connected 
as a $U_{q}(\Fg_{J})$-crystal, it follows from 
\eqref{eq:omega02} and \eqref{eq:omega03} that 
the following diagram commutes: 
\begin{equation*}
\begin{CD}
\CB_{t} @>{\sim}>{\Psi_{J}|_{\CB_{t}}}> 
\CB_{J}(\lambda_{t}) \\
@V{\omega}VV @VV{\omega}V \\
\CB_{t} @>{\sim}>{\Psi_{J}|_{\CB_{t}}}> 
\CB_{J}(\lambda_{t}). \\
\end{CD}
\end{equation*}
Hence we deduce that 
%
%
\begin{equation} \label{eq:reg02}
\Psi_{J}(\CB^{\omega})= 
 \CB_{J}^{\omega}(\lambda_{1}) \sqcup 
 \CB_{J}^{\omega}(\lambda_{2}) 
 \sqcup \cdots \sqcup 
 \CB_{J}^{\omega}(\lambda_{p^{\prime}}).
\end{equation}
Note that by Theorem~\ref{thm:fixed}, 
the set on the right-hand side of \eqref{eq:reg02}, 
equipped with the $\omega$-Kashiwara operators, 
becomes a $U_{q}(\ha{\Fg}_{\ha{J}})$-crystal 
isomorphic to the crystal base 
$\ha{\CB}_{\ha{J}}(\ha{\lambda}_{1}) \sqcup 
 \ha{\CB}_{\ha{J}}(\ha{\lambda}_{2}) 
 \sqcup \cdots \sqcup 
 \ha{\CB}_{\ha{J}}(\ha{\lambda}_{p^{\prime}})$ of 
 the integrable $U_{q}(\ha{\Fg}_{\ha{J}})$-module
$\ha{V}_{\ha{J}}(\ha{\lambda}_{1}) \oplus 
  \ha{V}_{\ha{J}}(\ha{\lambda}_{2}) 
  \oplus \cdots \oplus 
  \ha{V}_{\ha{J}}(\ha{\lambda}_{p^{\prime}})$, 
where $\ha{\lambda}_{t}:=(\pos)^{-1}(\lambda_{t})$ 
for $1 \le t \le p^{\prime}$. 

\vspace{5mm}

\noindent {\bf Step 2.} 
First we show that 
the set $\CB^{\omega} \cup \{\theta\}$ is stable under 
the $\omega$-Kashiwara operators $\ti{e}_{j}$, 
$j \in \ha{I}$. 
Let us fix (arbitrarily) $j \in \ha{I}$, and 
let $b \in \CB^{\omega}$ be such that $\ti{e}_{j}b \ne \theta$.
Set $\ha{J}:=\bigl\{j\bigr\} \subsetneq \ha{I}$, 
$J:=\bigl\{\omega^{k}(j) \mid 0 \le k \le N_{j}-1\bigr\} 
\subsetneq I$, and let $\Psi_{J}$ be 
the isomorphism \eqref{eq:reg01} 
of $U_{q}(\Fg_{J})$-crystals. 
Then, from the definitions 
\eqref{eq:o-kas-op01} and \eqref{eq:o-kas-op02} of 
the $\omega$-Kashiwara operator $\ti{e}_{j}$ on 
$\bigl(
 \CB_{J}(\lambda_{1}) \sqcup 
 \CB_{J}(\lambda_{2}) 
 \sqcup \cdots \sqcup 
 \CB_{J}(\lambda_{p})
 \bigr) \cup \{\theta\}$
and the one on $\CB \cup \{\theta\}$, respectively, we see 
that $\ti{e}_{j} \circ \Psi_{J}=\Psi_{J} \circ \ti{e}_{j}$. 
Namely, we have 
$\ti{e}_{j}b = \Psi_{J}^{-1} (\ti{e}_{j} (\Psi_{J}(b)))$.
Also, we know from \eqref{eq:reg02} that 
the image $\Psi_{J}(\CB^{\omega})$ is 
a disjoint union of the fixed point subsets 
$\CB_{J}^{\omega}(\lambda_{t})$ 
under the action of $\omega$ 
of the crystal bases $\CB_{J}(\lambda_{t})$ 
of integrable highest weight $U_{q}(\Fg_{J})$-modules. 
Hence, we see 
by Theorem~\ref{thm:fixed} that 
$\ti{e}_{j} (\Psi_{J}(b))$ is contained in 
$\Psi_{J}(\CB^{\omega})$. 
Consequently, $\ti{e}_{j}b=
\Psi_{J}^{-1} (\ti{e}_{j} (\Psi_{J}(b)))$ is 
contained in $\CB^{\omega}$.
Similarly, we can show that 
the set $\CB^{\omega} \cup \{\theta\}$ is stable under 
the $\omega$-Kashiwara operators $\ti{f}_{j}$, 
$j \in \ha{I}$. This proves the first assertion.

Next, let us prove the second assertion. 
We show only the equality 
$\ha{\vp}_{j}(b)=\ha{\ve}_{j}(b)+
(\ha{\wt}\,b)(\ha{h}_{j})$ for each $b \in \CB^{\omega}$ 
and $j \in \ha{I}$ (i.e., 
Condition (1) of \cite[Definition~4.5.1]{HK}); 
the other conditions immediately follow from 
the definition \eqref{eq:o-kas-op02} of the $\omega$-Kashiwara operators 
$\ti{e}_{j}$, $\ti{f}_{j}$, $j \in \ha{I}$, 
the definitions \eqref{eq:wtvevp}
of the maps $\ha{\wt}$, 
$\ha{\ve}_{j}$, $\ha{\vp}_{j}$, $j \in \ha{I}$, 
for $\CB^{\omega}$, and equality \eqref{eq:pos-aff01}. 
Let us fix (arbitrarily) $j \in \ha{I}$, 
and set $\ha{J}:=\{j\} \subsetneq \ha{I}$, 
$J:=\bigl\{\omega^{k}(j) \mid 0 \le k \le N_{j}-1\bigr\} 
\subsetneq I$ as above.
Then we deduce from \eqref{eq:reg02} and 
Theorem~\ref{thm:fixed} that 
for each $b \in \CB^{\omega}$, 
%
%
\begin{equation} \label{eq:reg03}
\ha{\vp}_{j}(\Psi_{J}(b))=
\ha{\ve}_{j}(\Psi_{J}(b))+
\bigl(\ha{\wt}(\Psi_{J}(b))\bigr)(\ha{h}_{j}).
\end{equation}
In addition, since $\Psi_{J}$ is 
an isomorphism of $U_{q}(\Fg_{J})$-crystals 
(see \eqref{eq:reg01}), we see from 
the definitions \eqref{eq:fixed} and 
\eqref{eq:wtvevp} of the maps 
$\ha{\wt}$, 
$\ha{\ve}_{j}$, $\ha{\vp}_{j}$ for 
$\CB_{J}^{\omega}(\lambda_{1}) \sqcup 
 \CB_{J}^{\omega}(\lambda_{2}) 
 \sqcup \cdots \sqcup 
 \CB_{J}^{\omega}(\lambda_{p^{\prime}})$ 
and for $\CB^{\omega}$ that 
%
%
\begin{equation} \label{eq:reg04} 
\ha{\wt}\,(\Psi_{J}(b))=\ha{\wt}\,b, \quad 
\ha{\ve}_{j}(\Psi_{J}(b))=\ha{\ve}_{j}(b), \quad 
\ha{\vp}_{j}(\Psi_{J}(b))=\ha{\vp}_{j}(b) 
\end{equation}
for each $b \in \CB^{\omega}$. 
Therefore, by combining \eqref{eq:reg03} 
and \eqref{eq:reg04}, 
we obtain that $\ha{\vp}_{j}(b)=\ha{\ve}_{j}(b)+
(\ha{\wt}\,b)(\ha{h}_{j})$, as desired.

Finally, the regularity of 
the $U_{q}^{\prime}(\ha{\Fg})$-crystal 
$\CB^{\omega}$ 
follows from \eqref{eq:reg02} and the comment after it.
This proves the second assertion. 
\end{proof}

By the same argument as in 
the proof of Proposition~\ref{prop:fixed-reg}, 
we also obtain the following 
Lemmas~\ref{lem:o-kas-op-n} and \ref{lem:vevp}, 
using Propositions~\ref{prop:o-kas-op-n} 
and \ref{prop:vevp}, respectively. 
%
%
\begin{lem} \label{lem:o-kas-op-n}
Let $\CB$ be a regular $U_{q}^{\prime}(\Fg)$-crystal 
with an action $\omega:\CB \rightarrow \CB$ of 
the diagram automorphism $\omega$ 
satisfying \eqref{eq:omega03}. 
Then, we have $\ti{e}(m)_{j}=(\ti{e}_{j})^{m}$ and 
$\ti{f}(m)_{j}=(\ti{f}_{j})^{m}$ 
on the fixed point subset $\CB^{\omega}$ 
for every $m \ge 0$ and $j \in \ha{I}$, 
where $\ti{e}(m)_{j}$ and $\ti{f}(m)_{j}$ are 
defined by\,{\rm:}
\begin{equation}
\ti{x}(m)_{j}=
  \begin{cases}
       x_{j}^{m} \, 
       x_{\omega(j)}^{2m} \,
       x_{j}^{m}
    & \text{if \ } c_{j}=1, \\[3mm]
       x_{j}^{m}\,
       x_{\omega(j)}^{m} \cdots 
       x_{\omega^{N_{j}-1}(j)}^{m}
    & \text{if \ } c_{j}=2,
  \end{cases}
\end{equation}
where $x$ is either $e$ or $f$. 
\end{lem}
%
%
%
\begin{lem} \label{lem:vevp}
Let $\CB$ be a regular $U_{q}^{\prime}(\Fg)$-crystal 
with an action $\omega:\CB \rightarrow \CB$ of 
the diagram automorphism $\omega$ 
satisfying \eqref{eq:omega03}. 
Then, for every element $b$ 
of the fixed point subset $\CB^{\omega}$, 
we have $\ha{\ve}_{j}(b)=\ve_{\omega^{k}(j)}(b)$ and 
$\ha{\vp}_{j}(b)=\vp_{\omega^{k}(j)}(b)$ 
for all $j \in \ha{I}$ and $k \ge 0$.
\end{lem}

Let $\CB$ be a regular $U_{q}^{\prime}(\Fg)$-crystal 
with an action $\omega:\CB \rightarrow \CB$ of 
the diagram automorphism $\omega$ 
satisfying \eqref{eq:omega03}. 
Then it follows from 
Proposition~\ref{prop:fixed-reg} that 
the fixed point subset $\CB^{\omega}$ becomes 
a regular $U_{q}^{\prime}(\ha{\Fg})$-crystal. 
Hence there exists a unique action 
$\ha{S}:\ha{W} \rightarrow 
\Bij(\CB^{\omega})$, $\ha{w} \mapsto \ha{S}_{\ha{w}}$, 
of the Weyl group $\ha{W}$ of the orbit Lie algebra 
$\ha{\Fg}$ on the set $\CB^{\omega}$ such that 
$\ha{S}_{\ha{r}_{j}}=\ha{S}_{j}$ for all $j \in \ha{I}$, 
where $\ha{S}_{j}$ is defined as in \eqref{eq:si} 
(see Proposition~\ref{prop:Weyl}). 
%
%
\begin{lem} \label{lem:Weyl}
We have $\ha{S}_{\ha{w}}=S_{\Theta(\ha{w})}$ on 
$\CB^{\omega}$ for all $\ha{w} \in \ha{W}$, 
where $\Theta:\ha{W} \rightarrow \ti{W}$ 
is the isomorphism of groups introduced at the end of 
\S\ref{subsec:orbit-aff}. In particular, 
the equality $\ha{S}_{j}=S_{w_{j}}$ 
holds on $\CB^{\omega}$ for all $j \in \ha{I}$. 
\end{lem}

\begin{proof}
We need only show that the equality 
$\ha{S}_{j}=S_{w_{j}}$ 
holds on $\CB^{\omega}$ 
for all $j \in \ha{I}$. 
Fix an arbitrary $j \in \ha{I}$. 
Let $b \in \CB^{\omega}$, and set 
$m:=(\ha{\wt}\, b)(\ha{h}_{j})$. 
Here we note that 
%
%
\begin{align} \label{eq:wt01}
(\wt b)(h_{j}) & = 
\frac{1}{N_{j}}\sum_{k=0}^{N_{j}-1}
\bigl((\omega^{\ast})^{k}(\wt b)\bigr)(h_{j})
\quad \text{since $\wt b \in \sw$} \nonumber \\[1.5mm]
& = 
\frac{1}{N_{j}}\sum_{k=0}^{N_{j}-1}
(\wt b)(\omega^{-k}(h_{j})) = 
\frac{1}{N_{j}}\sum_{k=0}^{N_{j}-1}
(\wt b)(h_{\omega^{-k}(j)}) \nonumber \\[1.5mm]
& = (\wt b)(P_{\omega}^{-1}(\ha{h}_{j})) 
\quad \text{by the definition 
\eqref{eq:pos-aff} of $P_{\omega}$} \nonumber \\[1.5mm]
& = \bigl((\pos)^{-1}(\wt b)\bigr)(\ha{h}_{j}) 
  = (\ha{\wt}\,b)(\ha{h}_{j}). 
\end{align}
So, we have $(\wt b)(h_{j})=m$. In addition, 
since $\wt b \in \sw$, it follows that 
%
%
\begin{equation} \label{eq:wt02}
(\wt b)(h_{\omega^{k}(j)})=m \quad 
\text{for all $0 \le k \le N_{j}-1$}.
\end{equation}

Now, assume that $m \ge 0$. 
Then we obtain that 
\begin{align*}
\ha{S}_{j}b= (\ti{f}_{j})^{m}b = \ti{f}(m)_{j} b = 
\begin{cases}
f_{j}^{m} f_{\omega(j)}^{2m} f_{j}^{m} b 
 & \text{if $c_{j}=1$}, \\[3mm]
f_{j}^{m} f_{\omega(j)}^{m} \cdots 
f_{\omega^{N_{j}-1}(j)}^{m} b
 & \text{if $c_{j}=2$},
\end{cases}
\end{align*}
from the definition of $\ha{S}_{j}$ and 
Lemma~\ref{lem:o-kas-op-n}. 
On the other hand, it follows that 
\begin{equation*}
S_{w_{j}}b=
\begin{cases}
f_{j}^{m} f_{\omega(j)}^{2m} f_{j}^{m} b 
 & \text{if $c_{j}=1$}, \\[3mm]
f_{j}^{m} f_{\omega(j)}^{m} \cdots 
f_{\omega^{N_{j}-1}(j)}^{m} b
 & \text{if $c_{j}=2$}.
\end{cases}
\end{equation*}
Indeed, if $c_{j}=1$, then 
we have $w_{j}=r_{j}r_{\omega(j)}r_{j}$ by 
Remark~\ref{rem:link} and 
the definition \eqref{eq:wj} of $w_{j}$, and hence 
$S_{w_{j}}b = S_{j}S_{\omega(j)}S_{j}b 
   = S_{j}S_{\omega(j)}f_{j}^{m}b$ by definition. 
Since $c_{j}=1$, we deduce from 
Remark~\ref{rem:link} and \eqref{eq:wt02} that 
\begin{align*}
(\wt (f_{j}^{m}b))(h_{\omega(j)}) & =
(\wt b - m \alpha_{j})(h_{\omega(j)})=
(\wt b)(h_{\omega(j)}) - m \alpha_{j}(h_{\omega(j)}) \\ 
& =m-m \times (-1)= 2m.
\end{align*}
Therefore, it follows that 
$S_{j}S_{\omega(j)}f_{j}^{m}b
   = S_{j}f_{\omega(j)}^{2m}f_{j}^{m}b$. 
Similarly, it follows that 
$S_{j}f_{\omega(j)}^{2m}f_{j}^{m}b = 
f_{j}^{m}f_{\omega(j)}^{2m}f_{j}^{m}b$. 
Thus, we obtain that 
$S_{w_{j}}b=f_{j}^{m} f_{\omega(j)}^{2m} f_{j}^{m} b$. 
The proof for the case where $c_{j}=2$ is easier.

Consequently, 
we obtain that $\ha{S}_{j}b=S_{w_{j}}b$ if $m \ge 0$. 
The proof for the case $m \le 0$ is similar. 
This proves the lemma. 
\end{proof}
%
%
%
%
\subsection{Fixed point subsets of tensor products of 
regular crystals.}
\label{subsec:fixed-tensor}

Let $\CB_{1}$ (resp., $\CB_{2}$) be a regular 
$U_{q}^{\prime}(\Fg)$-crystal
with an action of the diagram automorphism $\omega$ 
satisfying \eqref{eq:omega03}, with $\CB=\CB_{1}$ 
(resp., $\CB=\CB_{2}$). 
Define an action 
$\omega:\CB_{1} \otimes \CB_{2} \rightarrow 
\CB_{1} \otimes \CB_{2}$ of 
the diagram automorphism $\omega$ by: 
$\omega(b_{1} \otimes b_{2})= 
 \omega(b_{1}) \otimes \omega(b_{2})$ 
for $b_{1} \otimes b_{2} \in \CB_{1} \otimes \CB_{2}$. 
Then, we can easily check by the tensor product rule 
for $U_{q}^{\prime}(\Fg)$-crystals that this action 
$\omega:\CB_{1} \otimes \CB_{2} \rightarrow 
\CB_{1} \otimes \CB_{2}$ satisfies condition 
\eqref{eq:omega03}, with $\CB=\CB_{1} \otimes \CB_{2}$.
Hence it follows from 
Proposition~\ref{prop:fixed-reg} that 
the fixed point subset 
$(\CB_{1} \otimes \CB_{2})^{\omega}$ of 
$\CB_{1} \otimes \CB_{2}$ under this action 
becomes a regular 
$U_{q}^{\prime}(\ha{\Fg})$-crystal, when equipped with 
the $\omega$-Kashiwara operators $\ti{e}_{j}$ and 
$\ti{f}_{j}$, $j \in \ha{I}$. 

On the other hand, we obviously have
%
%
\begin{equation} \label{eq:fixed-ten01}
(\CB_{1} \otimes \CB_{2})^{\omega} = 
 \bigl\{b_{1} \otimes b_{2} \mid 
 b_{1} \in \CB_{1}^{\omega},\,
 b_{2} \in \CB_{2}^{\omega}\bigr\},
\end{equation}
where $\CB_{1}^{\omega}$ and $\CB_{2}^{\omega}$ 
are the fixed point subsets of $\CB_{1}$ and 
$\CB_{2}$ under the action of $\omega$, respectively. 
We know from Proposition~\ref{prop:fixed-reg} that 
both of the fixed point subsets 
$\CB_{1}^{\omega}$ and $\CB_{2}^{\omega}$ 
become regular $U_{q}^{\prime}(\ha{\Fg})$-crystals, 
when equipped with 
the $\omega$-Kashiwara operators $\ti{e}_{j}$ and 
$\ti{f}_{j}$, $j \in \ha{I}$. 
Denote by 
$\CB_{1}^{\omega} \ha{\otimes}\, \CB_{2}^{\omega}$ 
the tensor product $U_{q}^{\prime}(\ha{\Fg})$-crystal 
of $\CB_{1}^{\omega}$ and $\CB_{2}^{\omega}$; 
we use the symbol $\ha{\otimes}\,$ 
instead of $\otimes$ to emphasize that 
it means the tensor product of 
$U_{q}^{\prime}(\ha{\Fg})$-crystals. 
Also, we denote by 
$\ha{e}_{j}$, $j \in \ha{I}$, and 
$\ha{f}_{j}$, $j \in \ha{I}$, 
the raising Kashiwara operators and 
lowering Kashiwara operators, respectively, 
on the tensor product 
$U_{q}^{\prime}(\ha{\Fg})$-crystal 
$\CB_{1}^{\omega} \ha{\otimes}\, \CB_{2}^{\omega}$. 
%
%
\begin{lem} \label{lem:fixed-tensor}
Let $\Phi:
\CB_{1}^{\omega} \ha{\otimes}\, \CB_{2}^{\omega} 
\rightarrow (\CB_{1} \otimes \CB_{2})^{\omega}$ 
be the map defined by\,{\rm:}
$\Phi(b_{1} \ha{\otimes}\, b_{2}) = b_{1} \otimes b_{2}$ 
for $b_{1} \ha{\otimes}\, b_{2} \in 
\CB_{1}^{\omega} \ha{\otimes}\, \CB_{2}^{\omega}$. 
Then, $\Phi$ is an isomorphism of 
$U_{q}^{\prime}(\ha{\Fg})$-crystals between 
$\CB_{1}^{\omega} \ha{\otimes}\, \CB_{2}^{\omega}$
and $(\CB_{1} \otimes \CB_{2})^{\omega}$. 
\end{lem}
\begin{proof}
It is obvious by \eqref{eq:fixed-ten01} that 
$\Phi$ is bijective. In addition, it immediately follows that 
$\Phi$ preserves weights. So, it remains to show show that 
$\Phi \circ \ha{e}_{j}=\ti{e}_{j} \circ \Phi$ and 
$\Phi \circ \ha{f}_{j}=\ti{f}_{j} \circ \Phi$ 
for all $j \in \ha{I}$ since both of 
$\CB_{1}^{\omega} \ha{\otimes}\, \CB_{2}^{\omega}$
and $(\CB_{1} \otimes \CB_{2})^{\omega}$ are 
semiregular $U_{q}^{\prime}(\ha{\Fg})$-crystals. 

We show only the equality  
$\Phi \circ \ha{e}_{j}=
 \ti{e}_{j} \circ \Phi$ 
in the case where 
 $c_{j}=1$ (see Remark~\ref{rem:link}), since
the proof of the equality 
$\Phi \circ \ha{f}_{j}=
 \ti{f}_{j} \circ \Phi$ is similar, and 
the proof for the case where $c_{j}=2$ is easier
(cf. \cite[Proposition~6.4]{OSS1} in the case where $c_{j}=2$). 
Let $b_{1} \ha{\otimes}\, b_{2} \in 
\CB_{1}^{\omega} \ha{\otimes}\, \CB_{2}^{\omega}$. 
We deduce from Lemma~\ref{lem:vevp} and 
the tensor product rules for 
the $U_{q}^{\prime}(\Fg)$-crystal 
$\CB_{1} \otimes \CB_{2}$ and for 
the $U_{q}^{\prime}(\ha{\Fg})$-crystal 
$\CB_{1}^{\omega} \ha{\otimes}\, \CB_{2}^{\omega}$, 
together with a computation similar to \eqref{eq:wt01}, 
that 
$\ha{\ve}_{j}(b_{1} \ha{\otimes}\, b_{2})= 
\ve_{j}(b_{1} \otimes b_{2})$ 
(note that $b_{1} \otimes b_{2} = 
\Phi(b_{1} \ha{\otimes}\, b_{2})$ 
by definition), 
where 
$\ha{\ve}_{j}(b_{1} \ha{\otimes}\, b_{2}) 
  :=\max \bigl\{n \ge 0 \mid 
    \ha{e}_{j}^{n} (b_{1} \ha{\otimes}\, b_{2}) \ne \theta
    \bigr\}$. 
Also, we see 
from Lemma~\ref{lem:vevp} that 
$\ha{\ve}_{j}(b_{1} \otimes b_{2})= 
\ve_{j}(b_{1} \otimes b_{2})$. 
Combining these equalities,
we obtain that 
for $b_{1} \ha{\otimes}\, b_{2} \in 
\CB_{1}^{\omega} \ha{\otimes}\, \CB_{2}^{\omega}$, 
%
%
\begin{equation} \label{eq:fixed-ten02}
\ha{\ve}_{j}(b_{1} \ha{\otimes}\, b_{2})= 
\ha{\ve}_{j}(b_{1} \otimes b_{2}).
\end{equation}
Therefore, we conclude that 
$\ha{e}_{j}(b_{1} \ha{\otimes}\, b_{2}) = \theta$ 
if and only if $\ti{e}_{j}(\Phi(b_{1} \ha{\otimes}\, b_{2})) = 
\ti{e}_{j}(b_{1} \otimes b_{2}) = \theta$.

Now, we assume that 
$\ha{e}_{j}(b_{1} \ha{\otimes}\, b_{2}) \ne \theta$; 
note that $\ti{e}_{j}(b_{1} \otimes b_{2}) \ne \theta$, 
as seen above.
We further assume that 
$\ha{\vp}_{j}(b_{1}) \ge \ha{\ve}_{j}(b_{2})$, 
since the proof for the case where 
$\ha{\vp}_{j}(b_{1}) < \ha{\ve}_{j}(b_{2})$ is similar. 
Then we have $\ha{e}_{j}(b_{1} \ha{\otimes}\, b_{2}) = 
\ti{e}_{j} b_{1} \ha{\otimes}\, b_{2}$ 
by the definition of $\ha{e}_{j}$. 
Set $b_{1}^{\prime} \otimes b_{2}^{\prime}:=
\ti{e}_{j}(b_{1} \otimes b_{2})$. 
Because $\ha{\vp}_{j}(b_{1}) \ge \ha{\ve}_{j}(b_{2})$ 
by the assumption, it follows from Lemma~\ref{lem:vevp} that 
$\vp_{j}(b_{1}) \ge \ve_{j}(b_{2})$. Hence we deduce that
\begin{equation*}
b_{1}^{\prime} \otimes b_{2}^{\prime} 
= \ti{e}_{j}(b_{1} \otimes b_{2})
= (e_{j}e_{\omega(j)}^{2}e_{j})
  (b_{1} \otimes b_{2})
= e_{j}e_{\omega(j)}^{2}(e_{j}b_{1} \otimes b_{2}).
\end{equation*}
Here, we obviously have $b_{1} \in \CB_{1}^{\omega}$, 
$b_{2} \in \CB_{2}^{\omega}$ since 
$b_{1} \ha{\otimes}\, b_{2} \in 
\CB_{1}^{\omega} \ha{\otimes}\, \CB_{2}^{\omega}$. 
In addition, 
since $b_{1}^{\prime} \otimes b_{2}^{\prime} 
= \ti{e}_{j}(b_{1} \otimes b_{2}) \in 
(\CB_{1} \otimes \CB_{2})^{\omega}$ 
by Proposition~\ref{prop:fixed-reg}, 
it follows from \eqref{eq:fixed-ten01} that
$b_{1}^{\prime} \in \CB_{1}^{\omega}$, 
$b_{2}^{\prime} \in \CB_{2}^{\omega}$. 
Thus, we have
$\wt b_{1},\,\wt b_{1}^{\prime},\, 
 \wt b_{2},\,\wt b_{2}^{\prime} \in \sw$. 
Consequently, by the tensor product rule for crystals, 
we obtain that 
\begin{align*}
b_{1}^{\prime} \otimes b_{2}^{\prime} & =
e_{j}e_{\omega(j)}^{2}(e_{j}b_{1} \otimes b_{2}) \\ 
& = 
(e_{j}e_{\omega(j)}^{2}e_{j}b_{1}) \otimes b_{2}
\quad \text{or} \quad 
(e_{\omega(j)}e_{j}b_{1}) \otimes (e_{j}e_{\omega(j)}b_{2}).
\end{align*}
Indeed, 
the $b_{1}^{\prime} \otimes b_{2}^{\prime}$ 
cannot be equal to 
$(e_{\omega(j)}^{2}e_{j}b_{1}) \otimes (e_{j}b_{2})$, 
for example, since 
$\wt (e_{j}b_{2}) = \wt b_{2} + \alpha_{j}$ is 
not contained in $\sw$. 
Moreover, we have the following claim. 

\begin{claim}
The $b_{1}^{\prime} \otimes b_{2}^{\prime}$ 
cannot be equal to 
$\bigl(e_{\omega(j)}e_{j}b_{1}\bigr) \otimes 
\bigl(e_{j}e_{\omega(j)}b_{2}\bigr)$. 
\end{claim}

\noindent 
{\it Proof of Claim.} 
Suppose that $b_{1}^{\prime} \otimes b_{2}^{\prime}=
\bigl(e_{\omega(j)}e_{j}b_{1}\bigr) \otimes 
\bigl(e_{j}e_{\omega(j)}b_{2}\bigr)$. Then we have 
$b_{1}^{\prime}=e_{\omega(j)}e_{j}b_{1}$. 
Let $J:=\bigl\{j,\,\omega(j)\bigr\} \subsetneq I$, and 
denote by $\CB_{1}(b_{1})$ the connected component of
$\CB_{1}$, regarded as a $U_{q}(\Fg_{J})$-crystal,
containing the $b_{1} \in \CB_{1}$. 
Note that $b_{1}^{\prime}$ is also contained in 
$\CB_{1}(b_{1})$. 
We see that 
$\CB_{1}(b_{1})$ is isomorphic 
as a $U_{q}(\Fg_{J})$-crystal to 
the crystal base $\CB_{J}(\lambda)$ for some 
$J$-dominant integral weight 
$\lambda \in P_{\cl}$ (see \eqref{eq:reg01}). 
Because $\CB_{1}^{\omega} \cap \CB_{1}(b_{1})$ contains 
$b_{1}$ and $b_{1}^{\prime}$ (and hence is nonempty), 
we deduce as in Step 1 of the proof of 
Proposition~\ref{prop:fixed-reg} that
$\omega^{\ast}(\lambda)=\lambda$, and that 
$\CB_{1}^{\omega} \cap \CB_{1}(b_{1})$ 
equipped with $\omega$-Kashiwara operators becomes
a $U_{q}(\ha{\Fg}_{\ha{J}})$-crystal isomorphic to 
the crystal base $\ha{\CB}_{\ha{J}}(\ha{\lambda})$, 
where $\ha{J}:=\bigl\{j\}$ and 
$\ha{\lambda}:=(\pos)^{-1}(\lambda)$. 
Since $\ha{\CB}_{\ha{J}}(\ha{\lambda})$ is connected, 
it follows that both of $\ha{\wt}\,b_{1}$ and 
$\ha{\wt}\,b_{1}^{\prime}$ are contained 
in the set 
$\ha{\lambda} + \BZ \ha{\alpha}_{j}$, 
and hence that $\ha{\wt}\,b_{1}-
\ha{\wt}\,b_{1}^{\prime} \in \BZ\ha{\alpha}_{j}$. 
But, we obtain by \eqref{eq:pos-aff01} that 
\begin{equation*}
\ha{\wt}\,b_{1}^{\prime}=
\ha{\wt}\bigl(e_{\omega(j)}e_{j}b_{1}\bigr) = 
\ha{\wt}\,b_{1} + \frac{1}{2}\ha{\alpha}_{j} 
\not\in \ha{\wt}\,b_{1} + \BZ \ha{\alpha}_{j}, 
\end{equation*}
which is a contradiction. This proves the claim. 

\vspace{3mm}

Therefore, we conclude that 
$\ti{e}_{j}(b_{1} \otimes b_{2})= 
 (e_{j}e_{\omega(j)}^{2}e_{j}b_{1}) \otimes b_{2} =
 \ti{e}_{j}b_{1} \otimes b_{2}$, from which 
the equality $\Phi (\ha{e}_{j}(b_{1} \ha{\otimes}\, b_{2})) = 
\ti{e}_{j}(\Phi (b_{1} \ha{\otimes}\, b_{2}))$ 
immediately follows. 
This completes the proof of the lemma. 
\end{proof}

Because the fixed point subsets 
$\CB_{1}^{\omega}$, $\CB_{2}^{\omega}$, 
and $(\CB_{1} \otimes \CB_{2})^{\omega}$, 
equipped with the $\omega$-Kashiwara operators 
$\ti{e}_{j}$ and $\ti{f}_{j}$, $j \in I$, 
are all regular $U_{q}^{\prime}(\ha{\Fg})$-crystals 
by Propositioin~\ref{prop:fixed-reg}, 
there exist actions of the Weyl group $\ha{W}$ of $\ha{\Fg}$ 
on them by Proposition~\ref{prop:Weyl}, all of which are denoted 
by $\ha{S}_{\ha{w}}$, $\ha{w} \in \ha{W}$. 
With this notation, we have the following lemma.
%
%
\begin{lem} \label{lem:ak02}
Let $b_{1} \in \CB_{1}^{\omega}$ and 
$b_{2} \in \CB_{2}^{\omega}$ be $\ha{W}$-extremal elements
whose weights are contained in the same Weyl chamber with 
respect to the simple coroots $\ha{h}_{j}$, $j \in \ha{I}_{0}$. 
Then, $b_{1} \otimes b_{2} \in (\CB_{1} \otimes \CB_{2})^{\omega}$ 
is also a $\ha{W}$-extremal element, and the equality 
$\ha{S}_{\ha{w}}(b_{1} \otimes b_{2}) = \ha{S}_{\ha{w}}b_{1}
\otimes \ha{S}_{\ha{w}}b_{2}$ holds 
for all $\ha{w} \in \ha{W}$. 
\end{lem}

\begin{proof}
Since the tensor product 
$U_{q}^{\prime}(\Fg)$-crystal 
$\CB_{1}^{\omega} \ha{\otimes}\, \CB_{2}^{\omega}$ is 
regular, there exists an action 
of the Weyl group $\ha{W}$ on 
$\CB_{1}^{\omega} \ha{\otimes}\, \CB_{2}^{\omega}$ 
(see Proposition~\ref{prop:Weyl}). 
In addition, since $\Phi:
\CB_{1}^{\omega} \ha{\otimes}\,\CB_{2}^{\omega}
\rightarrow (\CB_{1} \otimes \CB_{2})^{\omega}$ is 
an isomorphism of $U_{q}^{\prime}(\ha{\Fg})$-crystals 
by Lemma~\ref{lem:fixed-tensor}, 
it follows from the definitions of the actions of 
the Weyl group $\ha{W}$ 
on $\CB_{1}^{\omega} \ha{\otimes}\, \CB_{2}^{\omega}$ and 
on $(\CB_{1} \otimes \CB_{2})^{\omega}$ that 
$\Phi \circ \ha{S}_{\ha{w}} = 
\ha{S}_{\ha{w}} \circ \Phi$ 
for all $\ha{w} \in \ha{W}$. 
Hence, we see that 
$b_{1} \otimes b_{2} \in (\CB_{1} \otimes \CB_{2})^{\omega}$ 
is a $\ha{W}$-extremal element 
since $\Phi^{-1}(b_{1} \otimes b_{2})=
b_{1} \ha{\otimes}\, b_{2}$ 
is a $\ha{W}$-extremal element of the tensor product 
$U_{q}^{\prime}(\ha{\Fg})$-crystal 
$\CB_{1}^{\omega} \ha{\otimes}\,\CB_{2}^{\omega}$
by Lemma~\ref{lem:ak}. 
Also, the equality 
$\ha{S}_{\ha{w}}(b_{1} \otimes b_{2}) = \ha{S}_{\ha{w}}b_{1}
\otimes \ha{S}_{\ha{w}}b_{2}$ immediately follows from 
Lemma~\ref{lem:ak} and the equality 
$\Phi \circ \ha{S}_{\ha{w}} = 
\ha{S}_{\ha{w}} \circ \Phi$. 
\end{proof}
%
%
%
%
\section{Proof of the main result.}
\label{sec:prf-main}
In this section, we use the notation of 
\S\ref{subsec:perfixed} and keep 
Assumption~\ref{ass}. 
Recall 
that we fixed (arbitrarily) 
$i \in \ha{I}_{0}=\ha{I} \setminus \{0\}$ and 
$s \in \BZ_{\ge 1}$. 
By \eqref{eq:omega}, 
$\ti{\CB}^{i,s}$ is a regular $U_{q}^{\prime}(\Fg)$-crystal
with an action $\omega:\ti{\CB}^{i,s} \rightarrow 
\ti{\CB}^{i,s}$ of the diagram automorphism $\omega$ 
satisfying condition \eqref{eq:omega03}, 
with $\CB=\ti{\CB}^{i,s}$.
Hence we can apply the results in 
\S\ref{subsec:fixed-reg} to the case where 
$\CB=\ti{\CB}^{i,s}$, 
$\CB^{\omega}=\ha{\CB}^{i,s}$, and 
those in \S\ref{subsec:fixed-tensor} to the case where 
$\CB_{1}=\CB_{2}=\ti{\CB}^{i,s}$, 
$\CB_{1}^{\omega}=\CB_{2}^{\omega}=\ha{\CB}^{i,s}$. 
%
%
\subsection{Proof of regularity.}
\label{subsec:prf-reg}
The following proposition immediately follows from 
Proposition~\ref{prop:fixed-reg} applied to 
the case $\CB=\ti{\CB}^{i,s}$ 
(note that $\CB^{\omega}=\ha{\CB}^{i,s}$).
%
%

\begin{prop} \label{prop:reg}
The subset $\ha{\CB}^{i,s} \cup \{\theta\}$ of 
$\ti{\CB}^{i,s} \cup \{\theta\}$ is stable 
under the $\omega$-Kashiwara operators $\ti{e}_{j}$ and 
$\ti{f}_{j}$ on $\ti{\CB}^{i,s} \cup \{\theta\}$
for all $j \in \ha{I}$. 
In addition, the $\ha{\CB}^{i,s}$ 
equipped with the maps 
$\ha{\wt}$, $\ti{e}_{j}$, $\ti{f}_{j}$, 
$j \in \ha{I}$, and $\ha{\ve}_{j}$, $\ha{\vp}_{j}$, 
$j \in \ha{I}$, becomes a regular 
$U_{q}^{\prime}(\ha{\Fg})$-crystal. 
\end{prop}
%
%
%
%
\subsection{Proof of simplicity.}
\label{subsec:prf-simple}
In this subsection, we prove that 
the $U_{q}^{\prime}(\ha{\Fg})$-crystal 
$\ha{\CB}^{i,s}$ is simple. First we show the 
following proposition. 
%
%
\begin{prop} \label{prop:S1}
The set $\ha{\CB}^{i,s}$ is of finite cardinality, and 
the weights of elements of 
$\ha{\CB}^{i,s}$ are all 
contained in $(\ha{P}_{\cl})_{0}$. 
Namely, the $U_{q}^{\prime}(\ha{\Fg})$-crystal 
$\ha{\CB}^{i,s}$ 
satisfies condition {\rm (S1)} 
of Definition~\ref{dfn:simple}. 
\end{prop}

\begin{proof}
Since $\ti{\CB}^{i,s}$ is the tensor product of 
the crystal bases $\CB^{\omega^{k}(i),s}$, 
$0 \le k \le N_{i}-1$, of finite cardinality, 
it follows that 
$\ti{\CB}^{i,s}$ is a finite set, and hence so is 
$\ha{\CB}^{i,s}$. 
Therefore, it suffices to 
show that $\ha{\wt}\,b \in (\ha{P}_{\cl})_{0}$ for all 
$b \in \ha{\CB}^{i,s}$. 
Let $b \in \ha{\CB}^{i,s}$. Then we deduce that 
%
%
\begin{align} \label{eq:wt03}
(\ha{\wt}\, b)(\ha{c}\,) 
 & = ((\pos)^{-1}(\wt b))(\ha{c}\,) \nonumber \\
 & = (\wt b)(P_{\omega}^{-1}(\ha{c}\,)) 
   \quad \text{(see the comment after \eqref{eq:pos-aff})} 
   \nonumber \\
 & = (\wt b)(c) 
   \quad \text{by \eqref{eq:pos-aff01}}. 
\end{align}
Because the $U_{q}^{\prime}(\Fg)$-crystal 
$\ti{\CB}^{i,s}$ is perfect (and hence simple), 
we have $(\wt b)(c)=0$, 
and hence $(\ha{\wt}\, b)(\ha{c}\,)=0$. 
This proves the proposition. 
\end{proof}
The rest of this subsection is devoted to proving that 
the $U_{q}^{\prime}(\ha{\Fg})$-crystal 
$\ha{\CB}^{i,s}$ satisfies condition (S2), 
and also condition (S3), of 
Definition~\ref{dfn:simple}. 
Since $\ti{\CB}^{i,s}$ is a perfect 
(and hence simple) $U_{q}^{\prime}(\Fg)$-crystal, 
it follows from Lemma~\ref{lem:dom} that 
there exists a unique $W$-extremal element of $\ti{\CB}^{i,s}$, 
denoted by $\ti{u}$, such that $\wt \ti{u}$ is $I_{0}$-dominant, 
i.e., $(\wt \ti{u})(h_{j}) \ge 0$ for all 
$j \in I_{0}=I \setminus \{0\}$ (this $\ti{u}$ is equal to the element 
$\ti{u}_{i,s} \in \ti{\CB}^{i,s}$ in Remark~\ref{rem:dom03}). 
%
%
%
%
\begin{lem} \label{lem:tiu}
The $W$-extremal element $\ti{u} \in \ti{\CB}^{i,s}$ is 
contained in the fixed point subset $\ha{\CB}^{i,s}$. 
Moreover, the element $\ti{u}$ is also 
a $\ha{W}$-extremal element of 
the $U_{q}^{\prime}(\ha{\Fg})$-crystal 
$\ha{\CB}^{i,s}$. Namely, for every $\ha{w} \in \ha{W}$, either 
$\ti{e}_{j} \ha{S}_{\ha{w}}\ti{u} = \theta$ or 
$\ti{f}_{j} \ha{S}_{\ha{w}}\ti{u} = \theta$ holds 
for each $j \in \ha{I}$. 
\end{lem}

\begin{proof}
We first show that $\omega(\ti{u})=\ti{u}$. Note that 
$\wt (\omega(\ti{u}))$ is $I_{0}$-dominant. In addition, 
for each $w \in W$, we have
\begin{equation} \label{eq:sim01}
S_{w}(\omega(\ti{u}))=
\omega(S_{w^{\prime}}\ti{u})
\end{equation}
for some $w^{\prime} \in W$. 
Indeed, we see from the definition of the action of $W$ 
on $\ti{\CB}^{i,s}$ and \eqref{eq:omega} that 
$S_{j}(\omega(b))=\omega(S_{\omega^{-1}(j)}b)$ for 
all $j \in I$ and $b \in \ti{\CB}^{i,s}$. Therefore, 
if $w=r_{j_{1}}r_{j_{2}} \cdots r_{j_{p}}$ is 
a reduced expression of $w \in W$, then we have 
\begin{align*}
S_{w}(\omega(\ti{u})) & =
S_{j_{1}}S_{j_{2}} \cdots S_{j_{p}}(\omega(\ti{u})) \\
& =
\omega(S_{\omega^{-1}(j_{1})}
S_{\omega^{-1}(j_{2})} \cdots 
S_{\omega^{-1}(j_{p})}\ti{u})
= \omega(S_{w^{\prime}}\ti{u}),
\end{align*}
where $w^{\prime}:=r_{\omega^{-1}(j_{1})}
r_{\omega^{-1}(j_{2})} \cdots 
r_{\omega^{-1}(j_{p})} \in W$. 
Since $S_{w^{\prime}}\ti{u} \in \ti{\CB}^{i,s}$ is 
a $W$-extremal element, 
we deduce from \eqref{eq:omega} and \eqref{eq:sim01}
that $\omega(\ti{u})$ is also a $W$-extremal element, 
whose weight $\wt (\omega(\ti{u}))$ is $I_{0}$-dominant. 
Consequently, 
it follows from the uniqueness of $\ti{u}$ 
(see Lemma~\ref{lem:dom}) that  
$\omega(\ti{u})=\ti{u}$. 
Thus we obtain that $\ti{u} \in \ha{\CB}^{i,s}$. 

Now, let $\ha{w} \in \ha{W}$. We know from 
Lemma~\ref{lem:Weyl} that 
$\ha{S}_{\ha{w}}\ti{u} = 
S_{\Theta(\ha{w})}\ti{u}$. 
Since $\ti{u}$ is a $W$-extremal element 
of the $U_{q}^{\prime}(\Fg)$-crystal 
$\ti{\CB}^{i,s}$, either 
$e_{j} \ha{S}_{\ha{w}}\ti{u} = 
 e_{j}S_{\Theta(\ha{w})}\ti{u} = 
 \theta$ or 
$f_{j} \ha{S}_{\ha{w}}\ti{u} = 
 f_{j}S_{\Theta(\ha{w})}\ti{u} = 
 \theta$ holds 
for each $j \in I$. Therefore, 
it follows from the definition of the 
$\omega$-Kashiwara operators 
$\ti{e}_{j}$, $\ti{f}_{j}$, $j \in \ha{I}$, that 
$\ti{e}_{j} \ha{S}_{\ha{w}}\ti{u} = \theta$ or 
$\ti{f}_{j} \ha{S}_{\ha{w}}\ti{u} = \theta$ 
for each $j \in \ha{I}$. This proves the lemma. 
\end{proof}

Let $b^{\prime} \in \ha{\CB}^{i,s}$ be 
a $\ha{W}$-extremal element of 
the regular $U_{q}^{\prime}(\ha{\Fg})$-crystal 
$\ha{\CB}^{i,s}$. 
Then there exists some
$\ha{w} \in \ha{W}$ such that 
$\ha{\wt}\,(\ha{S}_{\ha{w}}b^{\prime})$ is $\ha{I}_{0}$-dominant, 
i.e., $\bigl(\ha{\wt}\,(\ha{S}_{\ha{w}}b^{\prime})\bigr)
(\ha{h}_{j}) \ge 0$ for all 
$j \in \ha{I}_{0}=\ha{I} \setminus \{0\}$. 
Hence, in order to prove that 
the $U_{q}^{\prime}(\ha{\Fg})$-crystal 
$\ha{\CB}^{i,s}$ satisfies condition (S2) of 
Definition~\ref{dfn:simple}, 
it suffices to show the next proposition.  
\begin{prop} \label{prop:S2}
Let $b \in \ha{\CB}^{i,s}$ be 
a $\ha{W}$-extremal element 
of the regular $U_{q}^{\prime}(\ha{\Fg})$-crystal 
$\ha{\CB}^{i,s}$ such that 
$\ha{\wt}\,b$ is $\ha{I}_{0}$-dominant. 
Then we have $b=\ti{u}$. 
\end{prop}
In the proof below of Proposition~\ref{prop:S2}, 
we need some lemmas. 
%
%
\begin{lem} \label{lem:S2-01}
Let $\ha{\mu},\,\ha{\nu}$ be elements of $(\ha{P}_{\cl})_{0}$ 
that are contained in the same Weyl chamber with respect to 
$\ha{h}_{j}$, $j \in \ha{I}_{0}$. 
Then, for each $\ha{w} \in \ha{W}$, 
$\ha{w}(\ha{\mu})$ and $\ha{w}(\ha{\nu})$ are 
contained in the same Weyl chamber 
with respect to the simple coroots 
$\ha{h}_{j}$, $j \in \ha{I}_{0}$. 
Moreover, for each $\ha{w} \in \ha{W}$, 
we have either
$(\ha{w}(\ha{\mu}))(\ha{h}_{0}) \le 0$ and 
$(\ha{w}(\ha{\nu}))(\ha{h}_{0}) \le 0$, or 
$(\ha{w}(\ha{\mu}))(\ha{h}_{0}) \ge 0$ and 
$(\ha{w}(\ha{\nu}))(\ha{h}_{0}) \ge 0$.
\end{lem}

\begin{proof}
Recall that the Weyl group 
$\ha{W}$ of $\ha{\Fg}$ decomposes into 
the semidirect product $\ha{W}_{\ha{I}_{0}} 
\ltimes \ha{T}$ of the Weyl group 
$\ha{W}_{\ha{I}_{0}}:=
\langle \ha{r}_{j} \mid j \in \ha{I}_{0}\rangle$ 
(of finite type) and the abelian group 
$\ha{T}$ of translations. 
Hence, for each $\ha{w} \in \ha{W}$, 
there exists $\ha{z} \in \ha{W}_{\ha{I}_{0}}$ and 
$\ha{t} \in \ha{T}$ such that $\ha{w}=\ha{z}\,\ha{t}$. 
Then, since 
$\ha{\mu},\,\ha{\nu} \in (\ha{P}_{\cl})_{0}$ 
by the assumption, it follows from 
\cite[Chap.~6, formula (6.5.5)]{Kac} 
that $\ha{t}(\ha{\mu})=\ha{\mu}$, 
$\ha{t}(\ha{\nu})=\ha{\nu}$, 
and hence that 
$\ha{w}(\ha{\mu})=\ha{z}(\ha{\mu})$, 
$\ha{w}(\ha{\nu})=\ha{z}(\ha{\nu})$. 
Consequently, for each 
$\ha{w} \in \ha{W}$, $\ha{w}(\ha{\mu})$ and 
$\ha{w}(\ha{\nu})$ are contained 
in the same Weyl chamber 
with respect to the simple coroots 
$\ha{h}_{j}$, $j \in \ha{I}_{0}$.

Now, we fix an arbitrary $\ha{w} \in \ha{W}$. 
Note that $\ha{h}_{0}=\ha{c}-\ha{\theta}^{\vee}$, 
where $\ha{\theta}^{\vee}$ is the highest coroot 
of the dual root system of $\ha{\Fg}_{\ha{I}_{0}}$. 
Since $\ha{w}(\ha{\mu})$ and 
$\ha{w}(\ha{\nu})$ are contained in 
$(\ha{P}_{\cl})_{0}$, we have 
%
%
\begin{equation} \label{eq:S2-01-1}
\begin{array}{l} 
(\ha{w}(\ha{\mu}))(\ha{h}_{0})
  =(\ha{w}(\ha{\mu}))(\ha{c}-\ha{\theta}^{\vee})
  =-(\ha{w}(\ha{\mu}))(\ha{\theta}^{\vee}), \\[1.5mm]
(\ha{w}(\ha{\nu}))(\ha{h}_{0})
  =(\ha{w}(\ha{\nu}))(\ha{c}-\ha{\theta}^{\vee})
  =-(\ha{w}(\ha{\nu}))(\ha{\theta}^{\vee}).
\end{array}
\end{equation}
Because $\ha{w}(\ha{\mu})$ and $\ha{w}(\ha{\nu})$ 
are contained in the same Weyl chamber as shown above, 
we conclude from \eqref{eq:S2-01-1} that either 
$(\ha{w}(\ha{\mu}))(\ha{h}_{0}) \le 0$ and 
$(\ha{w}(\ha{\nu}))(\ha{h}_{0}) \le 0$, or 
$(\ha{w}(\ha{\mu}))(\ha{h}_{0}) \ge 0$ and 
$(\ha{w}(\ha{\nu}))(\ha{h}_{0}) \ge 0$. 
This proves the lemma. 
\end{proof}

%
%
\begin{lem} \label{lem:S2-02}
Let $\ha{\mu},\,\ha{\nu}$ be 
elements of $(\ha{P}_{\cl})_{0}$ 
that are $\ha{I}_{0}$-dominant, 
and assume that $\ha{\mu} \ne \ha{\nu}$. Then, 
there exists some $\ha{z} \in \ha{W}_{\ha{I}_{0}}$ 
such that $(\ha{z}(\ha{\mu}))(\ha{h}_{0}) \ne 
(\ha{z}(\ha{\nu}))(\ha{h}_{0})$, and 
$(\ha{z}(\ha{\mu}))(\ha{h}_{0}) \le 0$, 
$(\ha{z}(\ha{\nu}))(\ha{h}_{0}) \le 0$.
\end{lem}

\begin{proof}
First, let us show that 
there exists some 
$\ha{w} \in \ha{W}$ such that 
$(\ha{w}(\ha{\mu}))(\ha{h}_{0}) \ne 
(\ha{w}(\ha{\nu}))(\ha{h}_{0})$, and such that 
$(\ha{w}(\ha{\mu}))(\ha{h}_{0}) \le 0$, 
$(\ha{w}(\ha{\nu}))(\ha{h}_{0}) \le 0$.
Suppose that $(\ha{w}(\ha{\mu}))(\ha{h}_{0}) =
(\ha{w}(\ha{\nu}))(\ha{h}_{0})$ 
for all $\ha{w} \in \ha{W}$, or equivalently 
\begin{equation} \label{eq:S2-02-01}
\ha{\mu}(\ha{h}) =\ha{\nu}(\ha{h}) \quad 
\text{for all $\ha{h} \in \ha{W}\ha{h}_{0}$}.
\end{equation}
It is easy to show (cf. the argument 
in the hint for \cite[Exercise~6.9]{Kac}) 
that $\ha{\Fh}_{\cl}$ is generated by 
$\ha{W}\ha{h}_{0}$ as a $\BC$-vector space: 
$\ha{\Fh}_{\cl}=
 \sum_{\ha{h} \in \ha{W}\ha{h}_{0}} \BC \ha{h}$.
Therefore, it immediately follows from \eqref{eq:S2-02-01} 
that $\ha{\mu}(\ha{h}) =\ha{\nu}(\ha{h})$ for all 
$\ha{h} \in \ha{\Fh}_{\cl}$, and hence that 
$\ha{\mu}=\ha{\nu}$, which contradicts the assumption. 
This shows that 
$(\ha{w}(\ha{\mu}))(\ha{h}_{0}) \ne 
 (\ha{w}(\ha{\nu}))(\ha{h}_{0})$ 
for some $\ha{w} \in \ha{W}$. 
Furthermore, for this $\ha{w} \in \ha{W}$, 
we see from Lemma~\ref{lem:S2-01} that either 
$(\ha{w}(\ha{\mu}))(\ha{h}_{0}) \le 0$ and 
$(\ha{w}(\ha{\nu}))(\ha{h}_{0}) \le 0$, or 
$(\ha{w}(\ha{\mu}))(\ha{h}_{0}) \ge 0$ and 
$(\ha{w}(\ha{\nu}))(\ha{h}_{0}) \ge 0$. 
If $(\ha{w}(\ha{\mu}))(\ha{h}_{0}) \ge 0$ and 
$(\ha{w}(\ha{\nu}))(\ha{h}_{0}) \ge 0$, then 
$(\ha{r}_{0}\ha{w}(\ha{\mu}))(\ha{h}_{0}) \le 0$ and 
$(\ha{r}_{0}\ha{w}(\ha{\nu}))(\ha{h}_{0}) \le 0$, with 
$(\ha{r}_{0}\ha{w}(\ha{\mu}))(\ha{h}_{0}) \ne 
(\ha{r}_{0}\ha{w}(\ha{\nu}))(\ha{h}_{0})$.
Thus, we have obtained an element $\ha{w} \in \ha{W}$ 
such that $(\ha{w}(\ha{\mu}))(\ha{h}_{0}) \ne 
(\ha{w}(\ha{\nu}))(\ha{h}_{0})$, and 
$(\ha{w}(\ha{\mu}))(\ha{h}_{0}) \le 0$, 
$(\ha{w}(\ha{\nu}))(\ha{h}_{0}) \le 0$.

If we write the $\ha{w} \in \ha{W}$ in the form 
$\ha{w}=\ha{z}\,\ha{t}$, with $\ha{z} \in \ha{W}_{\ha{I}_{0}}$ 
and $\ha{t} \in \ha{T}$, 
then we have $\ha{w}(\ha{\mu})=\ha{z}(\ha{\mu})$ and 
$\ha{w}(\ha{\nu})=\ha{z}(\ha{\mu})$ 
(see the proof of Lemma~\ref{lem:S2-01}), and hence that 
$(\ha{z}(\ha{\mu}))(\ha{h}_{0}) \ne 
 (\ha{z}(\ha{\nu}))(\ha{h}_{0})$, and 
$(\ha{z}(\ha{\mu}))(\ha{h}_{0}) \le 0$, 
$(\ha{z}(\ha{\nu}))(\ha{h}_{0}) \le 0$. 
This proves the lemma. 
\end{proof}

\begin{proof}[Proof of Proposition~\ref{prop:S2}]
Recall from \S\ref{subsec:perfixed} that 
$\ti{\CB}^{i,s}$ is a perfect 
$U_{q}^{\prime}(\Fg)$-crystal isomorphic to 
the crystal base of a tensor product of 
finite-dimensional irreducible $U_{q}^{\prime}(\Fg)$-modules. 
By Proposition~\ref{prop:perfect}\,(2), we have 
an energy function $H$ on 
the tensor product $U_{q}^{\prime}(\Fg)$-crystal 
$\ti{\CB}^{i,s} \otimes \ti{\CB}^{i,s}$. 
We show that if $b$ were not equal to $\ti{u}$, then 
the energy function $H$ would attain infinitely many different 
values on the finite set 
$\ti{\CB}^{i,s} \otimes \ti{\CB}^{i,s}$. 
This proves the proposition by contradiction.

Now, suppose that $b \ne \ti{u}$. 
Since $\ti{u} \in \ha{\CB}^{i,s}$ 
by Lemma~\ref{lem:tiu}, and $b \in \ha{\CB}^{i,s}$ by assumption, 
it follows that $\mu:=\wt \ti{u} \in (P_{\cl})_{0} \cap \sw$ and 
$\nu:=\wt b \in (P_{\cl})_{0} \cap \sw$. 
Note that $\mu \ne \nu$ since 
$(\ti{\CB}^{i,s})_{\mu}=
\bigl\{\ti{u}\bigr\}$ by 
condition (S3) of Definition~\ref{dfn:simple}. 
If we set $\ha{\mu}:=(\pos)^{-1}(\mu)$ and 
$\ha{\nu}:=(\pos)^{-1}(\nu)$, then 
by Proposition~\ref{prop:S1}, we have 
$\ha{\mu},\,\ha{\nu} \in (\ha{P}_{\cl})_{0}$. 
Let us take
$\ha{z} \in \ha{W}_{\ha{I}_{0}}$ 
such that 
$(\ha{z}(\ha{\mu}))(\ha{h}_{0}) \ne 
(\ha{z}(\ha{\nu}))(\ha{h}_{0})$, and 
such that 
$(\ha{z}(\ha{\mu}))(\ha{h}_{0}) \le 0$, 
$(\ha{z}(\ha{\nu}))(\ha{h}_{0}) \le 0$ 
(see Lemma~\ref{lem:S2-02}). 

\paragraph{Step 1.} 
Define an action of 
$\omega:\ti{\CB}^{i,s} \otimes 
\ti{\CB}^{i,s} \rightarrow 
\ti{\CB}^{i,s} \otimes 
\ti{\CB}^{i,s}$ of the diagram automorphism $\omega$ by: 
$\omega(b_{1} \otimes b_{2}) = 
\omega(b_{1}) \otimes \omega(b_{2})$ for 
$b_{1} \otimes b_{2} \in \ti{\CB}^{i,s} \otimes 
\ti{\CB}^{i,s}$, and let
$(\ti{\CB}^{i,s} \otimes 
  \ti{\CB}^{i,s})^{\omega}$ be the fixed point 
subset of $\ti{\CB}^{i,s} \otimes 
\ti{\CB}^{i,s}$ under this action of $\omega$ 
as in \S\ref{subsec:fixed-tensor}. 
Obviously, 
$\ti{u} \otimes b$ is contained in 
the fixed point subset 
$(\ti{\CB}^{i,s} \otimes 
  \ti{\CB}^{i,s})^{\omega}$. 
Furthermore, 
we deduce from Lemma~\ref{lem:ak02} that 
$\ti{u} \otimes b$ is a $\ha{W}$-extremal element 
of the regular $U_{q}^{\prime}(\ha{\Fg})$-crystal 
$(\ti{\CB}^{i,s} \otimes 
  \ti{\CB}^{i,s})^{\omega}$ 
equipped with $\omega$-Kashiwara operators, 
and that 
$\ha{S}_{\ha{z}} (\ti{u} \otimes b) = 
 \ha{S}_{\ha{z}}\ti{u} \otimes \ha{S}_{\ha{z}}b$. 
We set 
$p:=-(\ha{z}(\ha{\mu}))(\ha{h}_{0})$, and 
$q:=-(\ha{z}(\ha{\nu}))(\ha{h}_{0})$; note 
that $p \ne q$, and $p \ge 0$, $q \ge 0$. 
Then, because 
\begin{equation*}
\Bigl(
 \ha{\wt}\bigl(
 \ha{S}_{\ha{z}}\ti{u} \otimes
 \ha{S}_{\ha{z}}b\bigr)
\Bigr)(\ha{h}_{0}) = 
\bigl(\ha{z}(\ha{\mu})+\ha{z}(\ha{\nu})\bigr)
(\ha{h}_{0}) = -(p+q) \le 0, 
\end{equation*}
we have by the definition of $\ha{S}_{0}$, 
\begin{equation*}
\ha{S}_{0} \ha{S}_{\ha{z}} (\ti{u} \otimes b) 
= (\ti{e}_{0})^{p+q}(
  \ha{S}_{\ha{z}} \ti{u} \otimes \ha{S}_{\ha{z}}b) \\
= e_{0}^{p+q}(
  \ha{S}_{\ha{z}} \ti{u} \otimes \ha{S}_{\ha{z}}b)
\end{equation*}
(note that $\ti{e}_{0}=e_{0}$). 
In addition, 
since $\ha{S}_{\ha{z}}\ti{u}$ and $\ha{S}_{\ha{z}}b$ are 
$\ha{W}$-extremal elements of 
the $U_{q}^{\prime}(\ha{\Fg})$-crystal $\ha{\CB}^{i,s}$, 
and since 
$(\ha{\wt}\,\ha{S}_{\ha{z}}\ti{u})(\ha{h}_{0})=
(\ha{z}(\ha{\mu}))(\ha{h}_{0})=-p \le 0$, 
$(\ha{\wt}\,\ha{S}_{\ha{z}}b)(\ha{h}_{0})=
(\ha{z}(\ha{\nu}))(\ha{h}_{0})=-q \le 0$, 
we obtain that 
\begin{align*}
& 
\ve_{0}(\ha{S}_{\ha{z}}\ti{u})=
\ha{\ve}_{0}(\ha{S}_{\ha{z}}\ti{u}) = 0 
\quad \text{and} \quad 
\ve_{0}(\ha{S}_{\ha{z}}b)=
\ha{\ve}_{0}(\ha{S}_{\ha{z}}b) = 0, \\ 
& 
\vp_{0}(\ha{S}_{\ha{z}}\ti{u})=
\ha{\vp}_{0}(\ha{S}_{\ha{z}}\ti{u}) = p
\quad \text{and} \quad 
\vp_{0}(\ha{S}_{\ha{z}}b)=
\ha{\vp}_{0}(\ha{S}_{\ha{z}}b) = q 
\end{align*}
(note that $\ti{e}_{0}=e_{0}$, $\ti{f}_{0}=f_{0}$). 
Using these equalities, 
we can easily show by the tensor product rule 
for the $U_{q}^{\prime}(\Fg)$-crystal 
$\ti{\CB}^{i,s} \otimes \ti{\CB}^{i,s}$ that 
%
%
\begin{equation} \label{eq:S2-00}
e_{0}^{l}
  \bigl(
  \ha{S}_{\ha{z}} \ti{u} \otimes \ha{S}_{\ha{z}} b
  \bigr) 
  = 
  \begin{cases}
  \ha{S}_{\ha{z}} \ti{u} 
  \otimes 
  e_{0}^{l}(\ha{S}_{\ha{z}} b)
  & \text{for $0 \le l \le q$}, \\[3mm]
  e_{0}^{l-q}(\ha{S}_{\ha{z}} \ti{u})
  \otimes 
  e_{0}^{q}(\ha{S}_{\ha{z}}b)
  & \text{for $q \le l \le p+q$}.
  \end{cases}
\end{equation}
Thus, by using \eqref{eq:energy-e} 
(applied to the case $\CB=\ti{\CB}^{i,s}$) successively, 
we obtain that 
%
%
\begin{equation} \label{eq:S2-01}
H \Bigl(
  \ha{S}_{0} \ha{S}_{\ha{z}} 
  (\ti{u} \otimes b) 
  \Bigr)=
 H \Bigl(\ha{S}_{\ha{z}}\ti{u}
   \otimes 
   e_{0}^{q}(\ha{S}_{\ha{z}}b)
   \Bigr)+p. 
\end{equation}
Indeed, we deduce that 
\begin{align*}
& 
H \Bigl(
  \ha{S}_{0} \ha{S}_{\ha{z}} 
  (\ti{u} \otimes b) 
  \Bigr)
= 
H \Bigl(
  e_{0}^{p+q} 
  (\ha{S}_{\ha{z}}\ti{u} \otimes 
  \ha{S}_{\ha{z}} b) 
  \Bigr) \\[3mm]
& \hspace*{5mm} = 
 H \Bigl(
   e_{0}\bigl(
        e_{0}^{p-1}(\ha{S}_{\ha{z}}\ti{u})
        \otimes 
        e_{0}^{q}(\ha{S}_{\ha{z}}b)
        \bigr)
   \Bigr) \quad 
   \text{by \eqref{eq:S2-00}} \\[3mm]
& \hspace*{5mm} = 
 H \Bigl(
   e_{0}^{p-1}(\ha{S}_{\ha{z}}\ti{u})
   \otimes 
   e_{0}^{q}(\ha{S}_{\ha{z}}b)
   \Bigr)+1 \quad 
   \text{since $\ve_{0}
   \bigl(e_{0}^{q}
    (\ha{S}_{\ha{z}}b)
   \bigr)=0$}
   \\[3mm]
& \hspace*{5mm} = 
 H \Bigl(
   e_{0}\bigl(
        e_{0}^{p-2}(\ha{S}_{\ha{z}}\ti{u})
        \otimes 
        e_{0}^{q}(\ha{S}_{\ha{z}}b)
        \bigr)
   \Bigr)+1 \quad 
   \text{by \eqref{eq:S2-00}} \\[3mm]
& \hspace*{5mm} = 
 H \Bigl(
   e_{0}^{p-2}(\ha{S}_{\ha{z}}\ti{u})
   \otimes 
   e_{0}^{q}(\ha{S}_{\ha{z}}b)
   \Bigr)+2 \quad 
   \text{since $\ve_{0}
   \bigl(e_{0}^{q}
    (\ha{S}_{\ha{z}}b)
   \bigr)=0$}. 
\end{align*}
Continuing in this way, we finally obtain \eqref{eq:S2-01}. 
Furthermore, again by using \eqref{eq:energy-e} 
successively, we obtain that 
%
%
\begin{equation} \label{eq:S2-02}
 H \Bigl(\ha{S}_{\ha{z}}\ti{u}
   \otimes 
   e_{0}^{q}(\ha{S}_{\ha{z}}b)
   \Bigr) = 
H \bigl(
   \ha{S}_{\ha{z}}\ti{u}
   \otimes 
   \ha{S}_{\ha{z}}b
   \bigr) - q.
\end{equation}
Indeed, we deduce that 
\begin{align*}
& 
H \Bigl(
   \ha{S}_{\ha{z}}\ti{u}
   \otimes 
   e_{0}^{q}(\ha{S}_{\ha{z}}b)
   \Bigr) = 
H \Bigl(
   e_{0}\bigl(
      \ha{S}_{\ha{z}}\ti{u}
      \otimes 
      e_{0}^{q-1}(\ha{S}_{\ha{z}}b)
        \bigr)
   \Bigr)  \quad 
   \text{by \eqref{eq:S2-00}} \\[3mm]
& \hspace*{5mm} = 
H \Bigl(
      \ha{S}_{\ha{z}}\ti{u}
      \otimes 
      e_{0}^{q-1}(\ha{S}_{\ha{z}}b)
   \Bigr) -1 \quad 
\text{since $\vp_{0}
  \bigl(
   \ha{S}_{\ha{z}}\ti{u}
  \bigr) = 0$ and 
   $\ve_{0}
  \bigl(
   e_{0}^{q-1}(\ha{S}_{\ha{z}}b)
  \bigr) > 0$}
  \\[3mm]
& \hspace*{5mm} = 
H \Bigl(
   e_{0}\bigl(
      \ha{S}_{\ha{z}}\ti{u}
      \otimes 
      e_{0}^{q-2}(\ha{S}_{\ha{z}}b)
        \bigr)
   \Bigr)  \quad 
   \text{by \eqref{eq:S2-00}} \\[3mm]
& \hspace*{5mm} = 
H \Bigl(
      \ha{S}_{\ha{z}}\ti{u}
      \otimes 
      e_{0}^{q-2}(\ha{S}_{\ha{z}}b)
   \Bigr) -2 \quad 
  \text{since $\vp_{0}
   \bigl(
   \ha{S}_{\ha{z}}\ti{u}
   \bigr) = 0$ and 
   $\ve_{0}
   \bigl(
   e_{0}^{q-2}(\ha{S}_{\ha{z}}b)
   \bigr) > 0$}. 
\end{align*}
Continuing in this way, we finally obtain \eqref{eq:S2-02}. 
Hence, by combining \eqref{eq:S2-01} and \eqref{eq:S2-02}, 
we have
\begin{align*}
H \bigl(
  \ha{S}_{0} \ha{S}_{\ha{z}} 
  (\ti{u} \otimes b) 
  \bigr) & = 
H \bigl(
   \ha{S}_{\ha{z}}\ti{u}
   \otimes 
   \ha{S}_{\ha{z}}b
   \bigr) + p-q. 
\end{align*}
Also, we deduce from \eqref{eq:energy-e} and 
\eqref{eq:energy-f} (applied to the case $\CB=\ti{\CB}^{i,s}$) 
that 
\begin{equation*}
H \Bigl(
   \ha{S}_{\ha{z}}\ti{u}
   \otimes 
   \ha{S}_{\ha{z}}b
   \Bigr) = 
H \Bigl(
   \ha{S}_{\ha{z}} (\ti{u}
   \otimes b)
   \Bigr) = 
H(\ti{u} \otimes b), 
\end{equation*}
since $\ha{S}_{\ha{z}}$ is defined by using only 
$\ti{e}_{j}$ and $\ti{f}_{j}$, $j \in \ha{I}_{0}$, 
and hence by using only $e_{j}$ and $f_{j}$, $j \in I_{0}$. 
Thus we conclude that 
%
%
\begin{equation} \label{eq:S2-03}
H \Bigl(
  \ha{S}_{0} \ha{S}_{\ha{z}} 
  (\ti{u} \otimes b) 
  \Bigr)
= H(\ti{u} \otimes b) + p-q. 
\end{equation}

\paragraph{Step 2.}
We write $\ha{r}_{0} \in \ha{W}$ in the form 
$\ha{r}_{0}=\ha{z}^{\prime}\,\ha{t}$, with 
$\ha{z}^{\prime} \in \ha{W}_{\ha{I}_{0}}$ and 
$\ha{t} \in \ha{T}$. Then we obtain that 
%
%
\begin{align} \label{eq:S2-04}
\ha{z}^{\prime}(\ha{r}_{0}\ha{z}(\ha{\mu}))
 & =\ha{z}^{\prime}\,\ha{t}\,(\ha{r}_{0}\ha{z}(\ha{\mu}))
 \quad \text{since $\ha{r}_{0}\ha{z}(\ha{\mu}) \in 
 (\ha{P}_{\cl})_{0}$} \nonumber \\
 & = \ha{r}_{0}(\ha{r}_{0}\ha{z}(\ha{\mu}))=\ha{z}(\ha{\mu}), 
\end{align}
and similarly that 
$\ha{z}^{\prime}\,\ha{r}_{0}\ha{z}(\ha{\nu})=
\ha{z}(\ha{\nu})$. 
Set $b_{1} \otimes b_{2}:=
\ha{S}_{\ha{z}^{\prime}}\ha{S}_{0} \ha{S}_{\ha{z}}
(\ti{u} \otimes b)$. 
It follows from Lemma~\ref{lem:ak02} that 
$b_{1}=\ha{S}_{\ha{z}^{\prime}} \ha{S}_{0} \ha{S}_{\ha{z}}(\ti{u})$ and 
$b_{2}=\ha{S}_{\ha{z}^{\prime}} \ha{S}_{0} \ha{S}_{\ha{z}}(b)$. 
Hence, by \eqref{eq:S2-04}, we have 
$\ha{\wt}\,b_{1}=\ha{z}(\ha{\mu})$ and 
$\ha{\wt}\,b_{2}=\ha{z}(\ha{\nu})$. 
Now, by the same argument as in Step 1 above, 
we can deduce that 
\begin{equation*}
H(\ha{S}_{0}(b_{1} \otimes b_{2}))=
H(b_{1} \otimes b_{2}) + p-q.
\end{equation*}
Note that 
since $\ha{z}^{\prime} \in \ha{W}_{\ha{I}_{0}}$, 
the $\ha{S}_{\ha{z}^{\prime}}$ is defined 
by using only 
$\ti{e}_{j}$, $\ti{f}_{j}$, $j \in \ha{I}_{0}$, 
and hence by using 
only $e_{j}$, $f_{j}$, $j \in I_{0}$. 
Therefore, 
by \eqref{eq:energy-e} and \eqref{eq:energy-f}, 
we obtain that 
\begin{align*}
H(\ha{S}_{0}(b_{1} \otimes b_{2})) & 
= H(b_{1} \otimes b_{2}) + p-q
= H\bigl(
  \ha{S}_{\ha{z}^{\prime}} \ha{S}_{0} \ha{S}_{\ha{z}} 
  (\ti{u} \otimes b)\bigr) + p-q \\[1.5mm]
& = H\bigl(\ha{S}_{0} \ha{S}_{\ha{z}} 
  (\ti{u} \otimes b)\bigr) + p-q \\[1.5mm]
& = H(\ti{u} \otimes b) + 2(p-q) 
  \quad \text{by \eqref{eq:S2-03}}. 
\end{align*}
Repeating this argument, we can show that 
for every $k \in \BZ_{\ge 0}$, 
there exists some $x \in 
(\ti{\CB}^{i,s} \otimes \ti{\CB}^{i,s})^{\omega} 
\subset \ti{\CB}^{i,s} \otimes \ti{\CB}^{i,s}$ 
such that $H(x)=H(\ti{u} \otimes b) + 2k(p-q)$. 
This contradicts the fact that 
$\ti{\CB}^{i,s} \otimes \ti{\CB}^{i,s}$ 
is a finite set. Thus we have proved the proposition. 
\end{proof}
Finally, we show that the $\ha{\CB}^{i,s}$ 
satisfies condition (S3) of Definition~\ref{dfn:simple}. 
%
%
\begin{prop} \label{prop:S3}
Let $b \in \ha{\CB}^{i,s}$ be a $\ha{W}$-extremal element 
of the $U_{q}^{\prime}(\ha{\Fg})$-crystal $\ha{\CB}^{i,s}$, 
and set $\ha{\mu}:=\ha{\wt}\,b$. Then, the subset 
$(\ha{\CB}^{i,s})_{\ha{\mu}} \subset 
\ha{\CB}^{i,s}$ of all elements of weight 
$\ha{\mu}$ consists only of the element $b$, i.e., 
$(\ha{\CB}^{i,s})_{\ha{\mu}}=\bigl\{b\bigr\}$. 
\end{prop}

\begin{proof}
By Proposition~\ref{prop:S2}, together with 
the comment before it, we see that 
$b$ is contained in the $\ha{W}$-orbit of 
$\ti{u} \in \ha{\CB}^{i,s}$: 
$b=\ha{S}_{\ha{w}}\ti{u}$ 
for some $\ha{w} \in \ha{W}$. 
Furthermore, 
it follows from Lemma~\ref{lem:Weyl} that 
$b=S_{\Theta(\ha{w})}\ti{u}$, 
and hence that $b$ is contained also in the $W$-orbit of 
$\ti{u}$. 
Note that since $\ti{u} \in \ti{\CB}^{i,s}$ 
is a $W$-extremal element
of the simple $U_{q}^{\prime}(\Fg)$-crystal 
$\ti{\CB}^{i,s}$, so is the $b \in \ti{\CB}^{i,s}$. 
In particular, 
by condition (S3) of Definition~\ref{dfn:simple}, 
we have $(\ti{\CB}^{i,s})_{\wt b}=\bigl\{b\bigr\}$. 
Therefore, we obtain that $(\ha{\CB}^{i,s})_{\ha{\mu}}=
\bigl\{b\bigr\}$, as desired. 
\end{proof}

Combining Propositions~\ref{prop:S1}, \ref{prop:S2}, 
and \ref{prop:S3}, we now establish the simplicity of 
the $U_{q}^{\prime}(\ha{\Fg})$-crystal $\ha{\CB}^{i,s}$. 
%
%
%
%
\subsection{Proof of bijectivity.}
\label{subsec:prf-bij}
%
%
%
%
\begin{prop} \label{prop:bij}
The level of the simple 
$U_{q}^{\prime}(\ha{\Fg})$-crystal 
$\ha{\CB}^{i,s}$ is equal to $s$. Moreover, 
the maps $\ha{\ve}, \ha{\vp}:(\ha{\CB}^{i,s})_{\min} 
\rightarrow (\ha{P}_{\cl}^{+})_{s}$ are bijective. 
\end{prop}

\begin{proof}
We see from Lemma~\ref{lem:vevp} 
applied to the case where $\CB=\ti{\CB}^{i,s}$ and 
$\CB^{\omega}=\ha{\CB}^{i,s}$ that 
$\pos(\ha{\ve}(b))=\ve(b)$ 
for all $b \in \ha{\CB}^{i,s}$. 
Therefore, we deduce that 
\begin{align*}
(\ha{\ve}(b))(\ha{c}\,) 
 & = \bigl((\pos)^{-1}(\ve(b))\bigr)(\ha{c}\,) \nonumber \\
 & = (\ve(b))(P_{\omega}^{-1}(\ha{c}\,)) 
   \quad \text{(see the comment after \eqref{eq:pos-aff})} 
   \nonumber \\
 & = (\ve(b))(c) 
   \quad \text{by \eqref{eq:pos-aff01}},
\end{align*}
and hence 
that $(\ha{\ve}(b))(\ha{c}\,)=(\ve(b))(c) \ge s$ 
for all $b \in \ha{\CB}^{i,s} \subset \ti{\CB}^{i,s}$  
since $\ti{\CB}^{i,s}$ is a perfect 
$U_{q}^{\prime}(\Fg)$-crystal of level $s$. 
This implies that 
the level of $\ha{\CB}^{i,s}$ is 
greater than or equal to $s$. 
So, to prove that the level of 
$\ha{\CB}^{i,s}$ is equal to $s$, 
it suffices to show that 
there exists some $b \in \ha{\CB}^{i,s}$ 
such that $(\ha{\ve}(b))(\ha{c}\,)=s$. 
Let $\ha{\mu} \in (\ha{P}_{\cl}^{+})_{s}$; 
note that $(\ha{P}_{\cl}^{+})_{s} \ne \emptyset$ for all 
$s \in \BZ_{\ge 1}$ since $s \ha{\Lambda}_{0} 
\in (\ha{P}_{\cl}^{+})_{s}$. 
Then we deduce that 
$(\pos(\ha{\mu}))(c) =\ha{\mu}(P_{\omega}(c))= 
 \ha{\mu}(\ha{c}\,)=s$,
and hence that 
$\pos(\ha{\mu}) \in (P_{\cl}^{+})_{s}$. 
Because $\ti{\CB}^{i,s}$ is 
a perfect $U_{q}^{\prime}(\Fg)$-crystal 
of level $s$, there exists 
a (unique) $b \in (\ti{\CB}^{i,s})_{\min}$ such that 
$\ve(b)=\pos(\ha{\mu}) \in (P_{\cl}^{+})_{s}$. 
It follows that this $b$ is contained 
in $\ha{\CB}^{i,s}$, i.e., that $\omega(b)=b$. 
Indeed, we see from \eqref{eq:omega} that 
$\ve_{j}(\omega(b))=\ve_{\omega^{-1}(j)}(b)$ 
for all $j \in I$, and hence that
$\ve(\omega(b))=\omega^{\ast}(\ve(b))$. 
But, since $\ve(b)=\pos(\ha{\mu}) \in 
(P_{\cl}^{+})_{s} \cap \sw$, 
we have $\omega^{\ast}(\ve(b))=\ve(b)$. 
So, we have $\ve(\omega(b))=\ve(b)$. 
Because $\ve:(\ti{\CB}^{i,s})_{\min} 
 \rightarrow (P_{\cl}^{+})_{s}$ is bijective, 
we conclude that $\omega(b)=b$. 
Also, it follows from the formulas
$\pos(\ha{\ve}(b))=\ve(b)$ and 
$\pos(\ha{\mu})=\ve(b)$
that $\ha{\ve}(b)=\ha{\mu}$. 
Consequently, we have $(\ha{\ve}(b))(\ha{c})=
\ha{\mu}(\ha{c})=s$ 
since $\ha{\mu} \in (\ha{P}_{\cl}^{+})_{s}$. 
Thus, we have shown 
that $\lev \ha{\CB}^{i,s}=s$. 

Next, we prove that the map 
$\ha{\ve}:(\ha{\CB}^{i,s})_{\min} \rightarrow 
 (\ha{P}_{\cl}^{+})_{s}$ is bijective. 
The argument above shows that 
for each $\ha{\mu} \in (\ha{P}_{\cl}^{+})_{s}$, there 
exists some $b \in (\ha{\CB}^{i,s})_{\min}$ such that 
$\ha{\ve}(b)=\ha{\mu}$, 
which means that the map 
$\ha{\ve}:(\ha{\CB}^{i,s})_{\min} 
\rightarrow (\ha{P}_{\cl}^{+})_{s}$ is surjective. 
Let us show that 
$\ha{\ve}:(\ha{\CB}^{i,s})_{\min} 
\rightarrow (\ha{P}_{\cl}^{+})_{s}$ is injective. 
Assume that $\ha{\ve}(b)=\ha{\ve}(b^{\prime})$ 
for $b,\,b^{\prime} \in (\ha{\CB}^{i,s})_{\min}$. 
Note that, by Lemma~\ref{lem:vevp}, 
$\ve(b)=\pos(\ha{\ve}(b))=
 \pos(\ha{\ve}(b^{\prime}))=\ve(b^{\prime})$, 
and hence that 
$(\ve(b))(c)=(\ha{\ve}(b))(\ha{c}\,)=s$, 
$(\ve(b^{\prime}))(c)=(\ha{\ve}(b^{\prime}))(\ha{c}\,)=s$, i.e., 
that $b,\,b^{\prime} \in (\ti{\CB}^{i,s})_{\min}$. 
Therefore, we conclude from the equality 
$\ve(b)=\ve(b^{\prime})$ that $b=b^{\prime}$ since
$\ve:(\ti{\CB}^{i,s})_{\min} \rightarrow 
(P_{\cl}^{+})_{s}$ is bijective. 
Thus, we have shown that 
the map 
$\ha{\ve}:(\ha{\CB}^{i,s})_{\min} \rightarrow 
 (\ha{P}_{\cl}^{+})_{s}$ is injective, and hence bijective. 
Similarly, we can show that the map 
$\ha{\vp}:(\ha{\CB}^{i,s})_{\min} \rightarrow 
(\ha{P}_{\cl}^{+})_{s}$ is bijective. 
This proves the proposition. 
\end{proof}

By Proposition~\ref{prop:bij}, we can conclude that 
$\ha{\CB}^{i,s}$ is a perfect 
$U_{q}^{\prime}(\Fg)$-crystal of level $s$, 
thereby completing the proof of 
Theorem~\ref{thm:main}. 
%
%
%
%
\subsection{Relation to virtual crystals.}
\label{subsec:virtual}
Let $\ti{u}$ be the unique $W$-extremal element of 
the simple $U_{q}^{\prime}(\Fg)$-crystal
$\ti{\CB}^{i,s}$ whose weight is $I_{0}$-dominant. 
Then, by Lemma~\ref{lem:tiu}, 
the element $\ti{u}$ is contained in the fixed point subset
$\ha{\CB}^{i,s}$. 
In addition, 
since $\ha{\CB}^{i,s}$ is a perfect 
(and hence simple) $U_{q}^{\prime}(\ha{\Fg})$-crystal 
by Theorem~\ref{thm:main}, 
we see by Proposition~\ref{prop:simple}\,(1) that 
$\ha{\CB}^{i,s}$ is equal to 
the set of elements of $\ti{\CB}^{i,s}$ obtained 
by applying the $\omega$-Kashiwara operators 
$\ti{e}_{j}$, $\ti{f}_{j}$, $j \in \ha{I}$, 
successively to $\ti{u}$. 
Therefore, if the pair $(\Fg,\omega)$ is 
in Case (a) (resp., (c), (d), (e)) in \S\ref{subsec:diag-aff}, 
then $(\ha{\CB}^{i,s},\ti{\CB}^{i,s})$ 
agrees with the virtual crystal for 
$(D_{n+1}^{(2)}, A_{2n-1}^{(1)})$ 
defined in \cite[\S6.7]{OSS1} (resp., 
$(A_{2n-1}^{(2)}, D_{n+1}^{(1)})$, 
$(D_{4}^{(3)}, D_{4}^{(1)})$,
$(E_{6}^{(2)}, E_{6}^{(1)})$ defined in 
\cite[Definition~2.6]{OSS2}). 

As for combinatorial $R$-matrices, we have 
the following proposition 
(see also \cite[Definition--Conjecture~3.4]{OSS2} 
and the comment after it). 
\begin{prop}
Let us fix $i_{1} \in \ha{I}_{0}$, $s_{1} \in \BZ_{\ge 1}$, 
and $i_{2} \in \ha{I}_{0}$, $s_{2} \in \BZ_{\ge 1}$, 
for which Assumption~\ref{ass} is satisfied. 
Then, there exists 
an isomorphism $\ha{R}: 
\ha{\CB}^{i_{1},s_{1}} \ha{\otimes}\, \ha{\CB}^{i_{2},s_{2}} 
\rightarrow 
\ha{\CB}^{i_{2},s_{2}} \ha{\otimes}\, \ha{\CB}^{i_{1},s_{1}}$ of 
$U_{q}^{\prime}(\ha{\Fg})$-crystals.
\end{prop}

\begin{proof}
As in \S\ref{subsec:fixed-tensor}, 
we define an action
$\omega:
\ti{\CB}^{i_{1},s_{1}} \otimes \ti{\CB}^{i_{2},s_{2}}
\rightarrow 
\ti{\CB}^{i_{1},s_{1}} \otimes \ti{\CB}^{i_{2},s_{2}}$
of the diagram automorphism $\omega$ by: 
$\omega(b_{1} \otimes b_{2})=
 \omega(b_{1}) \otimes \omega(b_{2})$ 
 for $b_{1} \otimes b_{2} \in 
 \ti{\CB}^{i_{1},s_{1}} \otimes \ti{\CB}^{i_{2},s_{2}}$, 
and let 
$(\ti{\CB}^{i_{1},s_{1}} \otimes \ti{\CB}^{i_{2},s_{2}})^{\omega}$ 
be the fixed point subset under this action of $\omega$. 
Similarly, we define an action 
$\omega:
\ti{\CB}^{i_{2},s_{2}} \otimes \ti{\CB}^{i_{1},s_{1}}
\rightarrow 
\ti{\CB}^{i_{2},s_{2}} \otimes \ti{\CB}^{i_{1},s_{1}}$ of 
$\omega$ as above, and let 
$(\ti{\CB}^{i_{2},s_{2}} \otimes \ti{\CB}^{i_{1},s_{1}})^{\omega}$ 
be the fixed point subset under this action of $\omega$. 

\begin{claim}
Let $R:\ti{\CB}^{i_{1},s_{1}} \otimes \ti{\CB}^{i_{2},s_{2}}
\rightarrow
\ti{\CB}^{i_{2},s_{2}} \otimes \ti{\CB}^{i_{1},s_{1}}$ be the 
combinatorial $R$-matrix in 
Proposition~\ref{prop:perfect}\,{\rm(}1\,{\rm)}. 
Then we have $R\bigl((\ti{\CB}^{i_{1},s_{1}} \otimes 
\ti{\CB}^{i_{2},s_{2}})^{\omega}\bigr) \subset 
(\ti{\CB}^{i_{2},s_{2}} \otimes \ti{\CB}^{i_{1},s_{1}})^{\omega}$.
\end{claim}

\noindent
{\it Proof of Claim.}
Let $\ti{u}_{1}$ (resp., $\ti{u}_{2}$) be 
the unique $W$-extremal element of 
the simple $U_{q}^{\prime}(\Fg)$-crystal 
$\ti{\CB}^{i_{1},s_{1}}$
(resp., $\ti{\CB}^{i_{2},s_{2}}$) whose weight is 
$I_{0}$-dominant (see Lemma~\ref{lem:dom}). 
Then we see from Lemma~\ref{lem:tiu} and 
\eqref{eq:fixed-ten01} that 
$\ti{u}_{1} \otimes \ti{u}_{2} \in 
(\ti{\CB}^{i_{1},s_{1}} \otimes 
 \ti{\CB}^{i_{2},s_{2}})^{\omega}$ and 
$\ti{u}_{2} \otimes \ti{u}_{1} \in 
(\ti{\CB}^{i_{2},s_{2}} \otimes 
 \ti{\CB}^{i_{1},s_{1}})^{\omega}$. 
Now, let $b \in (\ti{\CB}^{i_{1},s_{1}} \otimes 
 \ti{\CB}^{i_{2},s_{2}})^{\omega}$.
Because $(\ti{\CB}^{i_{1},s_{1}} \otimes 
 \ti{\CB}^{i_{2},s_{2}})^{\omega}$ 
 is isomorphic to 
 $\ha{\CB}^{i_{1},s_{1}} \ha{\otimes}\, 
 \ha{\CB}^{i_{2},s_{2}}$
as a $U_{q}^{\prime}(\ha{\Fg})$-crystal
by Lemma~\ref{lem:fixed-tensor}, 
it follows from 
Theorem~\ref{thm:main} and Proposition~\ref{prop:simple}\,(2) 
that the $U_{q}^{\prime}(\ha{\Fg})$-crystal 
$(\ti{\CB}^{i_{1},s_{1}} \otimes 
 \ti{\CB}^{i_{2},s_{2}})^{\omega}$ is simple. 
Therefore, 
the $b \in (\ti{\CB}^{i_{1},s_{1}} \otimes 
 \ti{\CB}^{i_{2},s_{2}})^{\omega}$ can be written as: 
$b = \ti{x}_{j_{1}} \ti{x}_{j_{2}} \cdots \ti{x}_{j_{l}}
(\ti{u}_{1} \otimes \ti{u}_{2})$, where $\ti{x}_{j}$ is 
either $\ti{e}_{j}$ or $\ti{f}_{j}$ for each $j \in \ha{I}$.
Since $R$ is an isomorphism of $U_{q}^{\prime}(\Fg)$-crystals, 
we see from the definition of the $\omega$-Kashiwara operators 
$\ti{e}_{j}$ and $\ti{f}_{j}$, $j \in \ha{I}$, 
that $R(b)=\ti{x}_{j_{1}} \ti{x}_{j_{2}} \cdots \ti{x}_{j_{l}} 
R(\ti{u}_{1} \otimes \ti{u}_{2})$. 
Furthermore, it follows 
from Lemma~\ref{lem:ak} that $\ti{u}_{1} \otimes \ti{u}_{2}$
(resp., $\ti{u}_{2} \otimes \ti{u}_{1}$) is a $W$-extremal element of 
the simple $U_{q}^{\prime}(\Fg)$-crystal
$\ti{\CB}^{i_{1},s_{1}} \otimes \ti{\CB}^{i_{2},s_{2}}$ 
(resp. $\ti{\CB}^{i_{2},s_{2}} \otimes \ti{\CB}^{i_{1},s_{1}}$) 
whose weight is $I_{0}$-dominant. 
Because such an element of 
the simple $U_{q}^{\prime}(\Fg)$-crystal is unique 
by Lemma~\ref{lem:dom}, we conclude that 
$R(\ti{u}_{1} \otimes \ti{u}_{2})=\ti{u}_{2} \otimes \ti{u}_{1}$, 
and hence that 
$R(b)=\ti{x}_{j_{1}} \ti{x}_{j_{2}} \cdots \ti{x}_{j_{l}} 
(\ti{u}_{2} \otimes \ti{u}_{1})$. In addition, we deduce from 
Proposition~\ref{prop:fixed-reg} applied to the case 
$\CB=\ti{\CB}^{i_{2},s_{2}} \otimes \ti{\CB}^{i_{1},s_{1}}$ that 
$R(b) \in (\ti{\CB}^{i_{2},s_{2}} \otimes 
\ti{\CB}^{i_{1},s_{1}})^{\omega}$ since 
$\ti{u}_{2} \otimes \ti{u}_{1} 
\in (\ti{\CB}^{i_{2},s_{2}} \otimes 
\ti{\CB}^{i_{1},s_{1}})^{\omega}$.
This proves the claim. 

\vspace{3mm}

Now we define a map $\ha{R}:
\ha{\CB}^{i_{1},s_{1}} \ha{\otimes}\, \ha{\CB}^{i_{2},s_{2}}
\rightarrow 
\ha{\CB}^{i_{2},s_{2}} \ha{\otimes}\, \ha{\CB}^{i_{1},s_{1}}$ 
by the following commutative diagram: 
\begin{equation}
\begin{CD}
\ha{\CB}^{i_{1},s_{1}} \ha{\otimes}\, \ha{\CB}^{i_{2},s_{2}}
@>{\sim}>{\Phi_{1}}>
(\ti{\CB}^{i_{1},s_{1}} \otimes \ti{\CB}^{i_{2},s_{2}})^{\omega} \\
@V{\ha{R}}VV @VV{R}V \\
\ha{\CB}^{i_{2},s_{2}} \ha{\otimes}\, \ha{\CB}^{i_{1},s_{1}}
@>{\sim}>{\Phi_{2}}>
(\ti{\CB}^{i_{2},s_{2}} \otimes \ti{\CB}^{i_{1},s_{1}})^{\omega}, 
\end{CD}
\end{equation}
where the isomorphism $\Phi_{1}$ (resp., $\Phi_{2}$) of 
$U_{q}^{\prime}(\ha{\Fg})$-crystals on the top 
(resp., on the bottom) is 
given by Lemma~\ref{lem:fixed-tensor}. 
It immediately follows from this commutative diagram that 
$\ha{R}:
\ha{\CB}^{i_{1},s_{1}} \ha{\otimes}\, \ha{\CB}^{i_{2},s_{2}} 
\rightarrow 
\ha{\CB}^{i_{2},s_{2}} \ha{\otimes}\, \ha{\CB}^{i_{1},s_{1}}$ 
is an embedding of $U_{q}^{\prime}(\ha{\Fg})$-crystals. 
In addition, since the $U_{q}^{\prime}(\ha{\Fg})$-crystals 
$\ha{\CB}^{i_{1},s_{1}}$ and $\ha{\CB}^{i_{2},s_{2}}$ are 
perfect (and hence simple) by Theorem~\ref{thm:main}, 
we see from Proposition~\ref{prop:simple}\,(2) that 
$\ha{\CB}^{i_{2},s_{2}} \ha{\otimes}\, \ha{\CB}^{i_{1},s_{1}}$ 
is a simple $U_{q}^{\prime}(\ha{\Fg})$-crystal, and that 
$\ha{\CB}^{i_{2},s_{2}} \ha{\otimes}\, \ha{\CB}^{i_{1},s_{1}}$ 
is connected. 
Consequently, we deduce that 
$\ha{R}: 
\ha{\CB}^{i_{1},s_{1}} \ha{\otimes}\, \ha{\CB}^{i_{2},s_{2}} 
\rightarrow 
\ha{\CB}^{i_{2},s_{2}} \ha{\otimes}\, \ha{\CB}^{i_{1},s_{1}}$ is 
surjective, and hence bijective. 
Thus we have proved the proposition. 
\end{proof}

Similarly, using the isomorphism 
$\Phi:\ha{\CB}^{i,s} \ha{\otimes}\, \ha{\CB}^{i,s} 
\stackrel{\sim}{\rightarrow}
 (\ti{\CB}^{i,s} \otimes \ti{\CB}^{i,s})^{\omega} 
 \subset \ti{\CB}^{i,s} \otimes \ti{\CB}^{i,s}$
of $U_{q}^{\prime}(\ha{\Fg})$-crystals in 
Lemma~\ref{lem:fixed-tensor}, 
we can prove the following proposition. 

\begin{prop}
We define a $\BZ$-valued function 
$\ha{H}:\ha{\CB}^{i,s} \ha{\otimes}\, \ha{\CB}^{i,s} 
\rightarrow \BZ$ by 
\begin{equation}
\ha{H}(b_{1} \ha{\otimes}\, b_{2}):=
H(\Phi(b_{1} \ha{\otimes}\, b_{2}))=
H(b_{1} \otimes b_{2}) 
\quad \text{\rm for $b_{1} \ha{\otimes}\, b_{2} \in 
\ha{\CB}^{i,s} \ha{\otimes}\, \ha{\CB}^{i,s}$},
\end{equation}
where $H:\ti{\CB}^{i,s} \otimes \ti{\CB}^{i,s} \rightarrow \BZ$ 
is the energy function in 
Proposition~\ref{prop:perfect}\,{\rm(}2\,{\rm)}. 
Then, the function $\ha{H}$ enjoys 
the following property\,{\rm:}
\begin{equation}
\ha{H}(\ha{e}_{j}(b_{1} \ha{\otimes}\, b_{2})) = 
\begin{cases}
\ha{H}(b_{1} \ha{\otimes}\, b_{2})+1 
 & \text{\rm if $j=0$ and 
   $\ha{\vp}_{0}(b_{1}) \ge 
    \ha{\ve}_{0}(b_{2})$,} \\[1.5mm]
\ha{H}(b_{1} \ha{\otimes}\, b_{2})-1 
 & \text{\rm if $j=0$ and 
   $\ha{\vp}_{0}(b_{1}) < 
    \ha{\ve}_{0}(b_{2})$,} \\[1.5mm]
\ha{H}(b_{1} \ha{\otimes}\, b_{2})
 & \text{\rm if $j \ne 0$},
\end{cases}
\end{equation}
for all $j \in \ha{I}$ and $b_{1} \ha{\otimes}\, b_{2} \in 
\ha{\CB}^{i,s} \ha{\otimes}\, \ha{\CB}^{i,s}$ such that 
$\ha{e}_{j}(b_{1} \ha{\otimes}\, b_{2}) \ne \theta$. 
\end{prop}
%
%
%
%
%
%
\section{Branching rules 
with respect to $U_{q}(\ha{\Fg}_{\ha{I}_{0}})$.}
\label{sec:branch}
In this section, we use the notation of \S\ref{subsec:perfixed} and 
keep Assumption~\ref{ass} for the (arbitrarily) fixed 
$i \in \ha{I}_{0}=\ha{I} \setminus \{0\}$ and $s \in \BZ_{\ge 1}$.
Since $\ha{\CB}^{i,s}$ is a perfect (and hence regular)
$U_{q}^{\prime}(\ha{\Fg})$-crystal for 
$i \in \ha{I}_{0}=\ha{I} \setminus \{0\}$ and $s \in \BZ_{\ge 1}$ 
by Theorem~\ref{thm:main}, it decomposes, 
as a $U_{q}(\ha{\Fg}_{\ha{I}_{0}})$-crystal, 
into a direct sum of the crystal bases of 
integrable highest weight 
$U_{q}(\ha{\Fg}_{\ha{I}_{0}})$-modules. 
In this section, 
we explicitly describe the branching rule, i.e., 
how the $\ha{\CB}^{i,s}$ decomposes, with respect to 
the restriction to $U_{q}(\ha{\Fg}_{\ha{I}_{0}})$
for almost all $i \in \ha{I}_{0}$ 
and $s \in \BZ_{\ge 1}$. 
%
%
\subsection{Preliminary results.}
\label{subsec:pre}
Let us set
\begin{equation}
(\ha{\CB}^{i,s})_{\hw}:=
\bigl\{
 b \in \ha{\CB}^{i,s} \mid \ti{e}_{j}b = \theta 
 \text{\rm\  for all } j \in \ha{I}_{0}\bigr\}
\end{equation}
for $i \in \ha{I}_{0}$ and $s \in \BZ_{\ge 1}$. 
Then, by the regularity of 
the $U_{q}^{\prime}(\ha{\Fg})$-crystal $\ha{\CB}^{i,s}$, 
we have 
\begin{equation}
\ha{\CB}^{i,s} \cong 
 \bigoplus_{b \in (\ha{\CB}^{i,s})_{\hw}}
 \ha{\CB}_{\ha{I}_{0}}(\ha{\wt}\,b) 
 \quad \text{as $U_{q}(\ha{\Fg}_{\ha{I}_{0}})$-crystals}, 
\end{equation}
where $\ha{\CB}_{\ha{I}_{0}}(\ha{\lambda})$ denotes
the crystal base of the integrable highest weight 
$U_{q}(\ha{\Fg}_{\ha{I}_{0}})$-module of 
($\ha{I}_{0}$-dominant) 
highest weight $\ha{\lambda} \in \ha{P}_{\cl}$. 
Similarly, we set 
\begin{equation}
(\ti{\CB}^{i,s})_{\hw}:=
\bigl\{
 b \in \ti{\CB}^{i,s} \mid e_{j}b = \theta 
 \text{\rm\  for all } j \in I_{0}\bigr\}
\end{equation}
for $i \in \ha{I}_{0}$ and $s \in \BZ_{\ge 1}$. 
Then, by the regularity of 
the $U_{q}^{\prime}(\Fg)$-crystal 
$\ti{\CB}^{i,s}$, we have 
%
%
\begin{equation} \label{eq:br01}
\ti{\CB}^{i,s} \cong 
 \bigoplus_{b \in (\ti{\CB}^{i,s})_{\hw}}
 \CB_{I_{0}}(\wt b)
 \quad \text{as $U_{q}(\Fg_{I_{0}})$-crystals},
\end{equation}
where $\CB_{I_{0}}(\lambda)$ denotes 
the crystal base of the integrable highest weight 
$U_{q}(\Fg_{I_{0}})$-module of ($I_{0}$-dominant) 
highest weight $\lambda \in P_{\cl}$. 
%
%
\begin{lem} \label{lem:branch01}
We have
\begin{equation}
(\ha{\CB}^{i,s})_{\hw} 
=
(\ti{\CB}^{i,s})_{\hw} \cap \ha{\CB}^{i,s}
= 
\bigl\{
 b \in (\ti{\CB}^{i,s})_{\hw} \mid \omega(b)=b
 \bigr\}.
\end{equation}
\end{lem}

\begin{proof}
The inclusion 
$(\ha{\CB}^{i,s})_{\hw} \supset 
 (\ti{\CB}^{i,s})_{\hw} \cap 
 \ha{\CB}^{i,s}$ is obvious 
from the definition of the raising $\omega$-Kashiwara 
operators $\ti{e}_{j}$, $j \in \ha{I}$. 
For the reverse inclusion, 
let $b \in (\ha{\CB}^{i,s})_{\hw}$ 
(note that $\omega(b)=b$ by definition). 
Then it follows that 
$\ha{\ve}_{j}(b)=0$ for all $j \in \ha{I}_{0}$, and hence 
by Lemma~\ref{lem:vevp} that $\ve_{j}(b)=0$ 
for all $j \in I_{0}$. 
This implies that 
$e_{j}b = \theta$ for all $j \in I_{0}$. 
Hence we have 
$b \in (\ti{\CB}^{i,s})_{\hw} \cap 
 \ha{\CB}^{i,s}$. 
This proves the lemma. 
\end{proof}
%
%
\begin{lem} \label{lem:branch02}
Assume that the decomposition \eqref{eq:br01} 
of $\ti{\CB}^{i,s}$ 
as a $U_{q}(\Fg_{I_{0}})$-crystal is 
multiplicity-free. Then we have 
\begin{equation}
(\ha{\CB}^{i,s})_{\hw} = 
\bigl\{
 b \in (\ti{\CB}^{i,s})_{\hw} \mid 
 \omega^{\ast}(\wt b) =\wt b
 \bigr\}.
\end{equation}
\end{lem}

\begin{proof}
The inclusion $\subset$ is obvious 
by Lemma~\ref{lem:branch01}. 
For the reverse inclusion, 
let $b \in (\ti{\CB}^{i,s})_{\hw}$ be such that 
$\omega^{\ast}(\wt b) =\wt b$.
It immediately follows from 
\eqref{eq:omega} that $\omega(b)$ is also contained in 
$(\ti{\CB}^{i,s})_{\hw}$ 
(note that $\omega(I_{0})=I_{0}$). 
In addition, we have $\wt(\omega(b))=
\omega^{\ast}(\wt b)=\wt b$. 
Therefore, we deduce that $\omega(b)=b$, since 
the decomposition \eqref{eq:br01} of $\ti{\CB}^{i,s}$ 
as a $U_{q}(\Fg_{I_{0}})$-crystal is 
multiplicity-free by the assumption.
Thus, the $b \in (\ti{\CB}^{i,s})_{\hw}$ is 
contained in the set 
$(\ti{\CB}^{i,s})_{\hw} \cap 
 \ha{\CB}^{i,s}$, and 
hence in the set $(\ha{\CB}^{i,s})_{\hw}$ 
by Lemma~\ref{lem:branch01}. 
This proves the lemma. 
\end{proof}

In \S\ref{subsec:case-a} -- \S\ref{subsec:case-e} below,   
by using Lemma~\ref{lem:branch02}, 
we give an explicit description of 
the branching rule for $\ha{\CB}^{i,s}$ for 
$i \in \ha{I}_{0}$ and $s \in \BZ_{\ge 1}$
with respect to the restriction 
to $U_{q}(\ha{\Fg}_{\ha{I}_{0}})$, 
except a few cases where 
the decomposition \eqref{eq:br01} of 
$\ti{\CB}^{i,s}$ 
as a $U_{q}(\Fg_{I_{0}})$-crystal is not
multiplicity-free. In the following, we use 
the notation:
\begin{align}
& \vpi_{i}:=\Lambda_{i}-a_{i}^{\vee}\Lambda_{0} \in P_{\cl}
\quad \text{for $i \in I_{0}$}, \\[1.5mm]
& 
\ti{\vpi}_{i}:=
\sum_{k=0}^{N_{i}-1} \vpi_{\omega^{k}(i)} \in P_{\cl}
\quad \text{for $i \in \ha{I}_{0}$}, 
\quad \text{and} \quad 
\ti{\vpi}_{0}:=0, \\[1.5mm]
&
\ha{\vpi}_{i}:=\ha{\Lambda}_{i} - 
\ha{a}^{\vee}_{i} \ha{\Lambda}_{0} \in \ha{P}_{\cl}
\quad \text{for $i \in \ha{I}_{0}$}, 
\quad \text{and} \quad 
\ha{\vpi}_{0}:=0.
\end{align}
Note that $\pos(\ha{\vpi}_{i})=\ti{\vpi}_{i}$ 
for all $i \in \ha{I}_{0}$, and also for $i=0$. 

\begin{rem} \label{rem:dom03}
For $i \in \ha{I}_{0}$ and $s \in \BZ_{\ge 1}$, we set 
$\ti{u}_{i,s}=u_{i,s} \otimes u_{\omega(i),s} \otimes 
\cdots \otimes u_{\omega^{N_{i}-1}(i),s}$, where 
the $u_{\omega^{k}(i),s} \in \CB^{\omega^{k}(i),s}$, 
$0 \le k \le N_{i}-1$, are the $W$-extremal elements 
in Remark~\ref{rem:dom02}. Then, by Lemma~\ref{lem:ak}, 
this $\ti{u}_{i,s}$ is the $W$-extremal element 
of the simple $U_{q}^{\prime}(\Fg)$-crystal 
$\ti{\CB}^{i,s}$ such that $(\wt \ti{u}_{i,s})(h_{j}) \ge 0$ 
for all $j \in I_{0}$. In fact, we have 
$\wt \ti{u}_{i,s} = s \ti{\vpi}_{i}$. Therefore, it follows from 
Remark~\ref{rem:dom} that 
the weights of elements of $\ti{\CB}^{i,s}$ are 
all contained in the set 
$s\ti{\vpi}_{i} - \sum_{j \in I_{0}} \BZ_{\ge 0}\,\alpha_{j}$.
Also, the weights of elements of $\ha{\CB}^{i,s}$ are 
all contained in the set 
$s\ha{\vpi}_{i} - \sum_{j \in \ha{I}_{0}} 
\BZ_{\ge 0}\,\ha{\alpha}_{j}$.
\end{rem}

%
%
\subsection{Branching rule for Case (a).}
\label{subsec:case-a}
We know from \cite{KMN} that 
%
%
\begin{equation} \label{eq:case-a01}
\CB^{i,s} \cong \CB_{I_{0}}(s\vpi_{i}) \quad 
 \text{as $U_{q}(\Fg_{I_{0}})$-crystals} 
\end{equation}
for $i \in I_{0}$ and $s \in \BZ_{\ge 1}$. 
Here we recall the (well-known) fact 
that the integrable highest weight 
$U_{q}(\Fg_{I_{0}})$-module $V_{I_{0}}(\lambda)$ of 
$I_{0}$-dominant highest weight $\lambda \in P_{\cl}$ has
the same character as the integrable highest weight 
$\Fg_{I_{0}}$-module of the same highest weight, and this 
character is equal to 
$\sum_{b \in \CB_{I_{0}}(\lambda)} e^{\wt b}$. 
By using this fact, 
we deduce from Decomposition Rule on page~145 
of \cite{L} that if $1 \le i \le n-1$, then 
as $U_{q}(\Fg_{I_{0}})$-crystals, 
\begin{align*}
\ti{\CB}^{i,s} 
 & = \CB^{i,s} \otimes \CB^{\omega(i),s}
   \cong \CB_{I_{0}}(s\vpi_{i}) \otimes 
     \CB_{I_{0}}(s\vpi_{\omega(i)}) 
 \quad \text{by \eqref{eq:case-a01}} \\[3mm]
 & \cong \bigoplus_{
     \begin{subarray}{c}
     s_{0} + s_{1} + \cdots + s_{i} = s \\[1mm]
     s_{0},\,s_{1},\,\dots,\,s_{i} \in \BZ_{\ge 0}
     \end{subarray}}
     \CB_{I_{0}}(s_{0}\ti{\vpi}_{0} + s_{1}\ti{\vpi}_{1} + 
     \cdots + s_{i}\ti{\vpi}_{i})\\[3mm]
 & = \bigoplus_{
     \begin{subarray}{c}
     s_{1} + \cdots + s_{i} \le s \\[1mm]
     s_{1},\,\dots,\,s_{i} \in \BZ_{\ge 0}
     \end{subarray}}
     \CB_{I_{0}}(s_{1}\ti{\vpi}_{1} + 
     \cdots + s_{i}\ti{\vpi}_{i})
\end{align*}
for $s \in \BZ_{\ge 1}$. 
If $i=n$, then 
$\ti{\CB}^{n,s} = \CB^{n,s} \cong 
 \CB_{I_{0}}(s\vpi_{n})$ 
 as $U_{q}(\Fg_{I_{0}})$-crystals 
 for $s \in \BZ_{\ge 1}$.
Note that for every 
$i \in \ha{I}_{0}$ and $s \in \BZ_{\ge 1}$, 
the decomposition of $\ti{\CB}^{i,s}$ 
above is multiplicity-free, 
and that the highest weights appearing 
in the decomposition are all 
fixed by $\omega^{\ast}$. 
Consequently, by using Lemma~\ref{lem:branch02}, 
we obtain the branching rule as follows. 
%
%
\begin{prop} \label{prop:case-a}
For $i \in \ha{I}_{0}$ and $s \in \BZ_{\ge 1}$, we have
%
%
\begin{equation} \label{eq:br-a}
\ha{\CB}^{i,s} \cong 
\begin{cases}
\displaystyle{
\bigoplus_{
     \begin{subarray}{c}
     s_{1} + \cdots + s_{i} \le s \\[1mm]
     s_{1},\,\dots,\,s_{i} \in \BZ_{\ge 0}
     \end{subarray}}}
     \ha{\CB}_{\ha{I}_{0}}(s_{1}\ha{\vpi}_{1} + 
     \cdots + s_{i}\ha{\vpi}_{i})
     & \text{\rm if \, } 1 \le i \le n-1, \\[12mm]
\ha{\CB}_{\ha{I}_{0}}(s\ha{\vpi}_{n}) 
     & \text{\rm if \, } i=n,
\end{cases}
\end{equation}
as $U_{q}(\ha{\Fg}_{\ha{I}_{0}})$-crystals.
\end{prop}

%
%
\subsection{Branching rule for Case (b).}
\label{subsec:case-b}
In exactly the same way as in Case (a), 
we have the following proposition. 
%
%
\begin{prop} \label{prop:case-b}
For $i \in \ha{I}_{0}$ and 
$s \in \BZ_{\ge 1}$, we have
%
%
\begin{equation} \label{eq:br-b}
\ha{\CB}^{i,s} \cong 
\bigoplus_{
     \begin{subarray}{c}
     s_{1} + \cdots + s_{i} \le s \\[1mm]
     s_{1},\,\dots,\,s_{i} \in \BZ_{\ge 0}
     \end{subarray}}
     \ha{\CB}_{\ha{I}_{0}}(s_{1}\ha{\vpi}_{1} + 
     \cdots + s_{i}\ha{\vpi}_{i})
\end{equation}
as $U_{q}(\ha{\Fg}_{\ha{I}_{0}})$-crystals.
\end{prop}

%
%
\subsection{Branching rule for Case (c).}
\label{subsec:case-c}
We know from \cite{C} that the KR module 
$W^{(i)}_{s}(\zeta^{(i)}_{s})$ over $U_{q}^{\prime}(\Fg)$
for $i \in I_{0}$, $s \in \BZ_{\ge 1}$, and 
$\zeta^{(i)}_{s} \in \BC(q)^{\times}$, 
decomposes under the restriction 
to $U_{q}(\Fg_{I_{0}})$ as follows: 
\begin{align} 
& W^{(i)}_{s}(\zeta^{(i)}_{s}) \cong \nonumber \\
& \label{eq:case-c-w}
\begin{cases}
\displaystyle{
     \bigoplus_{
     \begin{subarray}{c}
     s_{p_{i}} + s_{p_{i}+2} + \cdots + s_{i} = s \\[0.5mm]
     s_{p_{i}},\,s_{p_{i}+2},\,\dots,\,s_{i} \in \BZ_{\ge 0}
     \end{subarray}}
     }
     V_{I_{0}}(s_{p_{i}}\vpi_{p_{i}} + 
     s_{p_{i}+2}\vpi_{p_{i}+2} + \cdots + s_{i}\vpi_{i}) & 
\text{if \, } 1 \le i \le n-1, \\[12mm]
V_{I_{0}}(s\vpi_{i}) & 
\text{if \, } i =n,\,n+1,
\end{cases}
\end{align}
as $U_{q}(\Fg_{I_{0}})$-modules, 
where the $p_{i} \in \bigl\{0,\,1\bigr\}$ for 
$1 \le i \le n-1$ is defined to be $0$ (resp., $1$) 
if $i$ is even (resp., odd), and $V_{I_{0}}(\lambda)$ 
is the integrable highest weight 
$U_{q}(\Fg_{I_{0}})$-module of ($I_{0}$-dominant) 
highest weight $\lambda \in P_{\cl}$. 
Accordingly, from \eqref{eq:case-c-w}, we obtain the following 
decomposition of the crystal base $\CB^{i,s}$ of 
$W^{(i)}_{s}(\zeta^{(i)}_{s})$, 
regarded as a $U_{q}(\Fg_{I_{0}})$-crystal
by restriction: 
\begin{align} 
& \CB^{i,s} \cong \nonumber \\
& \label{eq:case-c}
\begin{cases}
\displaystyle{
     \bigoplus_{
     \begin{subarray}{c}
     s_{p_{i}} + s_{p_{i}+2} + \cdots + s_{i} = s \\[0.5mm]
     s_{p_{i}},\,s_{p_{i}+2},\,\dots,\,s_{i} \in \BZ_{\ge 0}
     \end{subarray}}
     }
     \CB_{I_{0}}(s_{p_{i}}\vpi_{p_{i}} + 
     s_{p_{i}+2}\vpi_{p_{i}+2} + \cdots + s_{i}\vpi_{i}) & 
\text{if \, } 1 \le i \le n-1, \\[12mm]
\CB_{I_{0}}(s\vpi_{i}) & 
\text{if \, } i =n,\,n+1,
\end{cases}
\end{align}
for $s \in \BZ_{\ge 1}$. 
Consequently, as in Case (a), we deduce from 
Decomposition Rule on page~145 of \cite{L} that 
as $U_{q}(\Fg_{I_{0}})$-crystals, 
\begin{align*}
\ti{\CB}^{n,s}
 & = \CB^{n,s} \otimes \CB^{n+1,s} 
   \cong \CB_{I_{0}}(s\vpi_{n}) \otimes 
     \CB_{I_{0}}(s\vpi_{n+1}) 
 \quad \text{by \eqref{eq:case-c}} \\[3mm]
 & \cong \displaystyle{
     \bigoplus_{
     \begin{subarray}{c}
     s_{p_{n}} + s_{p_{n}+2} + \cdots + s_{n} = s \\[0.5mm]
     s_{p_{n}},\,s_{p_{n}+2},\,\dots,\,s_{n} \in \BZ_{\ge 0}
     \end{subarray}}
     }
     \CB_{I_{0}}(s_{p_{n}}\ti{\vpi}_{p_{n}} + 
     s_{p_{n}+2}\ti{\vpi}_{p_{n}+2} + \cdots + 
     s_{n}\ti{\vpi}_{n})
\end{align*}
for $s \in \BZ_{\ge 1}$. 
Also, if $1 \le i \le n-1$, then we have 
\begin{align*}
\ti{\CB}^{i,s}=\CB^{i,s} 
& \cong 
\displaystyle{
     \bigoplus_{
     \begin{subarray}{c}
     s_{p_{i}} + s_{p_{i}+2} + \cdots + s_{i} = s \\[0.5mm]
     s_{p_{i}},\,s_{p_{i}+2},\,\dots,\,s_{i} \in \BZ_{\ge 0}
     \end{subarray}}
     }
     \CB_{I_{0}}(s_{p_{i}}\vpi_{p_{i}} + 
     s_{p_{i}+2}\vpi_{p_{i}+2} + \cdots + s_{i}\vpi_{i})
 \quad \text{by \eqref{eq:case-c}} \\[3mm]
& = 
\displaystyle{
     \bigoplus_{
     \begin{subarray}{c}
     s_{p_{i}} + s_{p_{i}+2} + \cdots + s_{i} = s \\[0.5mm]
     s_{p_{i}},\,s_{p_{i}+2},\,\dots,\,s_{i} \in \BZ_{\ge 0}
     \end{subarray}}
     }
     \CB_{I_{0}}(s_{p_{i}}\ti{\vpi}_{p_{i}} + 
     s_{p_{i}+2}\ti{\vpi}_{p_{i}+2} + \cdots + s_{i}\ti{\vpi}_{i}).
\end{align*}
Because in all cases, the decomposition of 
$\ti{\CB}^{i,s}$ above is multiplicity-free, 
and the highest weights appearing in the decomposition 
are all fixed by $\omega^{\ast}$, 
we obtain the following proposition by using 
Lemma~\ref{lem:branch02}. 
%
%
\begin{prop} \label{prop:case-c}
For $i \in \ha{I}_{0}$ and 
$s \in \BZ_{\ge 1}$, we have
%
%
\begin{equation} \label{eq:br-c}
\ha{\CB}^{i,s} \cong 
     \bigoplus_{
     \begin{subarray}{c}
     s_{p_{i}} + s_{p_{i}+2} + \cdots + s_{i} = s \\[0.5mm]
     s_{p_{i}},\,s_{p_{i}+2},\,\dots,\,s_{i} \in \BZ_{\ge 0}
     \end{subarray}}
     \ha{\CB}_{\ha{I}_{0}}(s_{p_{i}}\ha{\vpi}_{p_{i}} + 
     s_{p_{i}+2}\ha{\vpi}_{p_{i}+2} + \cdots + 
     s_{i}\ha{\vpi}_{i})
\end{equation}
as $U_{q}(\ha{\Fg}_{\ha{I}_{0}})$-crystals. 
\end{prop}
%
%
\subsection{Branching rule for Case (d).}
\label{subsec:case-d}
First, we should remark that our numbering 
of the index set $I$ is different from 
the ones in \cite{HKOTY} and \cite{L}. 
As in Case (c), we obtain that  
as $U_{q}(\Fg_{I_{0}})$-crystals, 
\begin{align}
& \ti{\CB}^{1,s} \cong 
  \bigoplus_{0 \le s_{1} \le s}
  \CB_{I_{0}}(s_{1}\vpi_{1}), \label{eq:case-d1} \\[1.5mm]
& \ti{\CB}^{2,s} \cong 
  \CB_{I_{0}}(s\vpi_{2}) \otimes 
  \CB_{I_{0}}(s\vpi_{3}) \otimes 
  \CB_{I_{0}}(s\vpi_{4}), 
\end{align}
for $s \in \BZ_{\ge 1}$.
Note that in general, the irreducible decomposition of 
the tensor product $V_{I_{0}}(s\vpi_{2}) \otimes 
  V_{I_{0}}(s\vpi_{3}) \otimes 
  V_{I_{0}}(s\vpi_{4})$ as a $U_{q}(\Fg_{I_{0}})$-module 
is not multiplicity-free, and hence that 
the decomposition \eqref{eq:br01} of $\ti{\CB}^{2,s}$, 
regarded as a $U_{q}(\Fg_{I_{0}})$-crystal by restriction, 
is not multiplicity-free.
Hence we exclude this case, i.e., the case where $i=2$.
Then, by using Lemma~\ref{lem:branch02} as above, 
we obtain the following proposition. 
%
%
\begin{prop} \label{prop:case-d}
For $s \in \BZ_{\ge 1}$, we have
%
%
\begin{equation} \label{eq:br-d}
\ha{\CB}^{1,s} \cong 
  \bigoplus_{0 \le s_{1} \le s}
  \ha{\CB}_{\ha{I}_{0}}(s_{1}\ha{\vpi}_{1})
\end{equation}
as $U_{q}(\ha{\Fg}_{\ha{I}_{0}})$-crystals. 
\end{prop}
%
%
\subsection{Branching rule for Case (e).}
\label{subsec:case-e}
First, we should remark that our numbering 
of the index set $I$ is different from the ones in 
\cite{HKOTY} and \cite{L}. 
Assume that $i \ne 2$ 
(note that 
the decomposition of $\CB^{2,s}$ is known not to be 
multiplicity-free in general; 
cf.~the formula for $\CW^{(3)}_{s}$ 
in the case where $X_{n}=E_{6}$ on page 278 of 
\cite[Appendix A]{HKOTY}). 
Then, as in Case (c), we obtain from \cite{C}
the following decomposition of $\CB^{i,s}$, 
regarded as a $U_{q}(\Fg_{I_{0}})$-crystal by restriction: 
\begin{equation*}
\CB^{i,s} \cong 
\begin{cases}
\CB_{I_{0}}(s\vpi_{i})
 & \text{if $i = 4,\,6$}, \\[3mm]
{\displaystyle\bigoplus_{0 \le s_{1} \le s}}
\CB_{I_{0}}(s_{1}\vpi_{1})
 & \text{if $i=1$}, \\[7mm]
{\displaystyle\bigoplus_{
  \begin{subarray}{c}
  s_{3}+s_{6}=s \\
  s_{3},s_{6} \ge 0
  \end{subarray}}}
\CB_{I_{0}}(s_{3}\vpi_{3}+s_{6}\vpi_{6})
 & \text{if $i=3$}, \\[10mm]
{\displaystyle\bigoplus_{
  \begin{subarray}{c}
  s_{4}+s_{5}=s \\
  s_{4},s_{5} \ge 0
  \end{subarray}}}
\CB_{I_{0}}(s_{4}\vpi_{4}+s_{5}\vpi_{5})
 & \text{if $i=5$}, 
\end{cases}
\end{equation*}
for $s \in \BZ_{\ge 1}$. 
Consequently, as in Cases (a), (c), we deduce from 
Decomposition Rule on page~145 of \cite{L} that as 
$U_{q}(\Fg_{I_{0}})$-crystals, 
\begin{align}
\ti{\CB}^{1,s} & = \CB^{1,s} \cong 
\bigoplus_{0 \le s_{1} \le s} 
\CB_{I_{0}}(s_{1}\vpi_{1}), \label{eq:case-e1} \\[3mm]
\ti{\CB}^{4,s} & 
= \CB^{4,s} \otimes \CB^{6,s} \cong 
 \CB_{I_{0}}(s\vpi_{4}) \otimes \CB_{I_{0}}(s\vpi_{6}) 
\nonumber \\[3mm]
& \cong \bigoplus_{
    \begin{subarray}{c} 
    s_{1}+s_{4} \le s \\[1mm]
    s_{1},\,s_{4} \in \BZ_{\ge 0} 
    \end{subarray}} 
    \CB_{I_{0}}(s_{1}\vpi_{1}+s_{4}\ti{\vpi}_{4}),
\label{eq:case-e2}
\end{align}
for $s \in \BZ_{\ge 1}$. 
Here we have exclude the case $i=3$,
since the decomposition of $\ti{\CB}^{3,s}=
\CB^{3,s} \otimes \CB^{5,s}$, 
regarded as a $U_{q}(\Fg_{I_{0}})$-crystal
by restriction, is known not to be multiplicity-free in general.
Then, by using Lemma~\ref{lem:branch02} as above, 
we obtain the following proposition. 
%
%
\begin{prop} \label{prop:case-e}
For $s \in \BZ_{\ge 1}$, we have
%
%
\begin{equation} \label{eq:br-e}
\ha{\CB}^{1,s} \cong 
\bigoplus_{0 \le s_{1} \le s}
\ha{\CB}_{\ha{I}_{0}}(s_{1}\ha{\vpi}_{1}), 
\qquad 
\ha{\CB}^{4,s} \cong 
\bigoplus_{
    \begin{subarray}{c} 
    s_{1}+s_{4} \le s \\[1mm]
    s_{1},\,s_{4} \in \BZ_{\ge 0} 
    \end{subarray}}
    \ha{\CB}_{\ha{I}_{0}}
    (s_{1}\ha{\vpi}_{1}+s_{4}\ha{\vpi}_{4}), 
\end{equation}
as $U_{q}(\ha{\Fg}_{\ha{I}_{0}})$-crystals. 
\end{prop}
%
%
%
\subsection{Comments.}
\label{subsec:comments}
In addition to Assumption~\ref{ass} for $\Fg$, let us assume 
that Conjecture~\ref{conj:hkott} with $\Fg$ replaced by $\ha{\Fg}$ 
holds for the (fixed) $i \in \ha{I}_{0}$ and $s \in \BZ_{\ge 1}$. 
Denote the KR module over $U_{q}^{\prime}(\ha{\Fg})$ having 
a crystal base by $\ha{W}^{(i)}_{s}(\ha{\zeta}^{(i)}_{s})$, where 
$\ha{\zeta}^{(i)}_{s} \in \BC(q)^{\times}$. Then we make the 
following observations. 

\begin{obs}
By comparing \eqref{eq:br-a} (resp., 
\eqref{eq:br-b}, \eqref{eq:br-c}) 
with the formula 
for $\CW_{s}^{(a)}$ in the case where 
$\Fg=D_{n+1}^{(2)}$ (resp., $A_{2n}^{(2)}$, 
$A_{2n-1}^{(2)}$) on 
page 247 (resp., page 247, page 246) of 
\cite[Appendix~A]{HKOTT} 
(with ``$a$'' replaced by ``$i$'', 
``$\ba{\Lambda}_{i}$'' replaced by 
``$\ha{\vpi}_{i}$'', and ``$q$'' equal to ``$1$''), 
we can check that they indeed agree. 
Similarly, 
by comparing \eqref{eq:br-d} 
(resp., \eqref{eq:br-e}) 
with the formula 
for $\CW_{s}^{(1)}$ 
(resp., $\CW_{s}^{(1)}$ and $\CW_{s}^{(4)}$) 
in the case where 
$\Fg=D_{4}^{(3)}$ (resp., $E_{6}^{(2)}$) 
on page 249 (resp., page 248) of 
\cite[Appendix~A]{HKOTT} 
(with ``$\ba{\Lambda}_{i}$'' replaced by 
``$\ha{\vpi}_{i}$'', and ``$q$'' equal to ``$1$''), 
we can check that they also agree. Here the formula 
for $\CW^{(a)}_{s}$ (with $q=1$) describes the branching 
rule for the KR module 
$\ha{W}^{(a)}_{s}(\ha{\zeta}^{(a)}_{s})$ over 
$U_{q}^{\prime}(\ha{\Fg})$ with respect to the restriction 
to $U_{q}(\ha{\Fg}_{\ha{I}_{0}})$. 
\end{obs}

Motivated by the observations above, 
we propose the following conjecture. 
%
%
%
\begin{conj} \label{conj:kr}
Let us fix {\rm(}arbitrarily{\rm)} $i \in \ha{I}_{0}$ 
and $s \in \BZ_{\ge 1}$, and keep Assumption~\ref{ass}. 
Then, the perfect $U_{q}^{\prime}(\ha{\Fg})$-crystal 
$\ha{\CB}^{i,s}$ is isomorphic to the 
{\rm(}conjectural\,{\rm)} crystal base of the KR module 
$\ha{W}^{(i)}_{s}(\ha{\zeta}^{(i)}_{s})$ over the 
twisted quantum affine algebra $U_{q}^{\prime}(\ha{\Fg})$. 
\end{conj}
%
%

{\small \setlength{\baselineskip}{12pt}
\renewcommand{\refname}{References}

}

\end{document}